\documentclass[fontsize=12pt,a4paper,headings=normal,
twoside=false,leqno,parskip=half-,abstract=true]{scrartcl}
\usepackage[english]{babel}
\usepackage[utf8]{inputenc}
\usepackage{csquotes}
\usepackage[hyphens]{url}
\usepackage{hyperref}

\hypersetup{
 pdftitle={scalar blow-up},
 pdfauthor={Bernold Fiedler},
 colorlinks=true,
 linkcolor=blue,
 citecolor=blue,
 filecolor=blue,
 urlcolor=blue}

 \usepackage{afterpage}

\usepackage{graphicx}
\usepackage{wrapfig} 
\usepackage[format=plain,labelfont=bf,font=small]{caption}
\usepackage{xcolor}
\usepackage[arrow, matrix, curve]{xy}
\usepackage{float}

\usepackage{caption}
\usepackage{tabulary}
\usepackage{array}
\newcolumntype{N}[1]{>{\centering\arraybackslash}m{#1}}

\usepackage{amsmath,amsthm}
\swapnumbers 
\usepackage{amssymb, eurosym} 

\makeatletter
\newcommand{\tpitchfork}{%
  \vbox{
    \baselineskip\z@skip
    \lineskip-.52ex
    \lineskiplimit\maxdimen
    \m@th
    \ialign{##\crcr\hidewidth\smash{$-$}\hidewidth\crcr$\pitchfork$\crcr}
  }%
}
\makeatother
\usepackage{latexsym}
\usepackage{enumerate}

\usepackage[notref,notcite,color,final 
]{showkeys}

\definecolor{refkey}{rgb}{1,0,0}
\definecolor{labelkey}{rgb}{1,0,0}

\usepackage{cancel}

\usepackage{tikz}

  \mathchardef\ordinarycolon\mathcode`\:
  \mathcode`\:=\string"8000
  \begingroup \catcode`\:=\active
    \gdef:{\mathrel{\mathop\ordinarycolon}}
  \endgroup

\theoremstyle{plain}
\newtheorem{thm}{Theorem}[section]
\newtheorem{lem}[thm]{Lemma}

\newtheorem{cor}[thm]{Corollary}

\newcommand\eps{\varepsilon}
\newcommand\mi{\mathrm{i}}
\renewcommand\theta{\vartheta}
\renewcommand\rho{\varrho}
\renewcommand\Re{\mathrm{Re}\,}
\renewcommand\Im{\mathrm{Im}\,}

\begin{document}

\title{\LARGE{Scalar polynomial vector fields\\ 
in real and complex time }}
\vspace{1cm}
{\subtitle{}
	\vspace{1ex}
	{}}\vspace{1ex}

\author{
 \\
\emph{-- In memoriam Leonid P. Shilnikov --}\\
{~}\\
Bernold Fiedler* 
\\
\vspace{2cm}}

\date{\small{version of \today}}
\maketitle
\thispagestyle{empty}

\vfill

*\\
Institut für Mathematik\\
Freie Universität Berlin\\
Arnimallee 3\\ 
14195 Berlin, Germany


\newpage
\pagestyle{plain}
\pagenumbering{roman}
\setcounter{page}{1}

\begin{abstract}
\noindent
Recent PDE studies address global boundedness versus finite-time blow-up in equations like the quadratic parabolic heat equation versus the nonconservative quadratic Schrödinger equation.
The two equations are related by passage from real to purely imaginary time.
Renewed interest in pioneering work by Masuda, in particular, has further explored the option to circumnavigate blow-up in real time, by a detour in complex time.

\smallskip\noindent
In the present paper, the simplest scalar ODE case is studied for polynomials
\begin{equation*}
\label{*}
\dot{w}=f(w)=(w-e_0)\cdot\ldots\cdot(w-e_{d-1})\,,  \tag{*}
\end{equation*}
of degree $d$ with $d$ simple complex zeros.
The explicit solution by separation of variables and explicit integration is an almost trivial matter.

\smallskip\noindent
In a classical spirit, indeed, we describe the complex Riemann surface $\mathcal{R}$ of the global nontrivial solution $(w(t),t)$ in complex time, as an unbranched cover of the punctured Riemann sphere $w\in\widehat{\mathbb{C}}_d:=\widehat{\mathbb{C}}\setminus\{e_0,\ldots,e_{d-1}\}$\,.
The flow property, however, fails at $w=\infty\in\widehat{\mathbb{C}}_d$.
The global consequences depend on the period map of the residues $2\pi\mi/f'(e_j)$ of $1/f$ at the punctures, in detail.
We therefore show that polynomials $f$ exist for arbitrarily prescribed residues with zero sum.
This result is not covered by standard interpolation theory.

\smallskip\noindent
Motivated by the PDE case, we also classify the planar \emph{real-time} phase portraits of \eqref{*}.
Here we prefer a Poincaré compactification of $w\in\mathbb{C}=\mathbb{R}^2$ by the closed unit disk. This regularizes $w=\infty$ by $2(d-1)$ equilibria, alternatingly stable and unstable within the invariant circle boundary at infinity.
In structurally stable hyperbolic cases of nonvanishing real parts $\Re f'(e_j)\neq 0$, for the linearizations at all equilibria $e_j$\,, and in absence of saddle-saddle heteroclinic orbits, we classify all compactified phase portraits, up to orientation preserving orbit equivalence and time reversal. 
Combinatorially, their source/sink connection graphs correspond to the planar trees of $d$ vertices or, dually, the circle diagrams with $d-1$ non-intersecting chords.
The correspondence provides an explicit count of the above equivalence classes of ODE \eqref{*}, in real time.

\smallskip\noindent
We conclude with a discussion of some higher-dimensional problems.
Not least, we offer a 1,000 \euro\ reward for the discovery, or refutation, of complex entire homoclinic orbits.

\end{abstract}

\newpage
\tableofcontents


\newpage
\pagenumbering{arabic}
\setcounter{page}{1}

\section{Introduction and main results} \label{Int}

\numberwithin{equation}{section}
\numberwithin{figure}{section}
\numberwithin{table}{section}

\subsection{PDE motivation}\label{Para}

Why study complex time?

Some PDE motivation arises from a comparison between the nonlinear heat equation
\begin{equation}
\label{PDEw}
w_r=w_{xx}+w^2-1\,,
\end{equation}
on the one hand, and, on the other hand, the nonlinear and nonconservative Schrödinger equation
\begin{equation}
\label{PDEpsi}
\mi\psi_s=\psi_{xx}+\psi^2-1\,.  
\end{equation}
Indices indicate partial derivatives, and $r,s$ stand for real time in the two PDEs, respectively.
To be completely specific, we consider consider both equations on an interval $0<x<\ell$ of length $\ell$, and under Neumann boundary conditions $u_x=\psi_x=0$, at $x=0,\ell$.

Traditionally, solutions $w$ are thought of as real-valued.
Complex-valued solutions $w=u+\mi v$ would then be considered as solutions of the system of two coupled reaction-diffusion equations which arises when we split \eqref{PDEw}  into real and imaginary parts.
For some relevant literature in this context, see for example \cite{Masuda1,Masuda2,LiSinai,COS,Yanagida,Kevrekidis,Stukediss,Stukearxiv,Jaquetteqp,JaquetteHet, JaquetteStuke,JaquetteMasuda, Fasondini23, Fasondini24}.
See also \eqref{ODEuv2} below.

Unlike general Schrödinger solutions $\psi$, solutions $w$ of the nonlinear heat equation \eqref{PDEw} are often real analytic in time $r$.
Complex-valued extensions of real (or complex) solutions $w(r,x)$ to complex times $t=r+\mi s$ then provide Schrödinger solutions
\begin{equation}
\label{psiw}
\psi(s,x):=w(r_0-\mathrm{i}s,x)
\end{equation}
Here $r=r_0$ is any fixed real part $\Re t$ of complex time $t$, and the Schrödinger solution proceeds along the vertical imaginary ``time'' direction $\mi s$.
Conversely, any complex analytic Schrödinger solutions $t\mapsto\psi(t,x)$ of \eqref{PDEpsi} define complex solutions \begin{equation}
\label{wpsi}
r\mapsto w(r,x):=\psi(-s_0+\mi r)\,, 
\end{equation}
for any fixed $s_0=-\Re t$.
In conclusion, complex time extension provides families of Schrödinger solutions $\psi$, from a single heat solution $w$, and vice versa. 
Moreover, the families of one type are related by a (semi)flow of the other type.
In fact, the two semiflows commute, locally and on analytic solutions.

Here and below, the terms \emph{``analytic''} and \emph{``analyticity''} emphasize local expansions by convergent power series.
\emph{``Holomorphic''} emphasizes complex differentiability and, therefore, Cauchy-Riemann equations.
For continuously real differentiable functions, e.g., the two notions coincide.
\emph{``Entire''} functions are globally analytic, e.g. for all complex time arguments $t\in\mathbb{C}$.

One PDE peculiarity of the heat equation \eqref{PDEw} is the frequent appearance of heteroclinic orbits, in real time $r$.
\emph{Heteroclinic} solutions $\Gamma(r)$ in real time $r\in\mathbb{R}$ connect time-independent equilibria $w=W_\pm$\,, i.e.
\begin{equation}
\label{heteroclinic}
\Gamma(r) \rightarrow W_\pm\,, \quad\textrm{for}\quad r\rightarrow\pm\infty\,,
\end{equation}
in suitable Banach spaces of solutions.
Depending on context, they are also called \emph{connecting orbits}, \emph{traveling fronts} (or backs), or \emph{solitons}.
We abbreviate heteroclinicity as $\Gamma:W_-\leadsto W_+$\,.
Unless specified otherwise, explicitly, we subsume the \emph{homoclinic} case $W_+=W_-$ of non-constant $\Gamma(r)$ under the heteroclinic label.

It turns out that any globally bounded, real-valued solution $w(r,x), \ r\in\mathbb{R}$ of \eqref{PDEw} is in fact real heteroclinic, $w=\Gamma:W_-\leadsto W_+$ between distinct equilibria $W_\pm$.
This is due to a decreasing Lyapunov functional on real solutions of \eqref{PDEw}.
See for example \cite{brfi88,brfi89,firoSFB,firoFusco,LappicyBlowup}, the introduction of \cite{fiestu24}, and the many references there.

In \cite{fiestu24} we have presented a preliminary analytical study of the complex heat-Schrödinger correspondence between \eqref{PDEw} and \eqref{PDEpsi}.
The interval length $\ell$, essentially, serves as a bifurcation parameter.
Our results hold for most lengths $\ell>0$, with mostly just discrete sets of exceptions.
We establish that the extensions of heteroclinic orbits $w(r,\cdot)=\Gamma(r):W_-\leadsto W_+$\,, from real time $r$ to complex time $t=r+\mi s$, cannot be complex entire for all $t\in\mathbb{C}$.
Failure to extend holomorphically to all $t\in\mathbb{C}$ entails blow-up (and blow-down) of the corresponding Schrödinger solutions \eqref{psiw}, for suitably fixed real $r_0$.
Here and below, \emph{blow-up} means that solutions become unbounded in some finite positive time; here for $0\leq s\nearrow s^*<\infty$.
Similarly, the term \emph{blow-down} is used to describe the same phenomenon in reverse, finite negative time.

Assumptions actually limit the above results to the case where the target $W_+=-1$ is the unique stable equilibrium, and hence spatially homogeneous.
The source equilibrium $W_-$ is required to be unstable hyperbolic, and of unstable dimension (alias Morse index) $i(W_-)$ not exceeding 22.
Only in the homogeneous case $W_-=+1$ were we able to drop that Morse limitation.
Alternatively, Schrödinger blow-up also occurs for $\Gamma$ which emanate from the fast unstable manifold at $W_-$ of dimension less than $1+i(W_-)/\sqrt{2}$.

The results are based on local Poincaré linearization \cite{ArnoldODE,Ilyashenko} in the locally analytic, finite-dimensional, (fast) unstable manifold of the heteroclinic source $W_-$\,.
The prerequisite spectral non-resonances impose the above restrictions on our blow-up results.
By Sturm-Liouville theory, the spectra of linearizations consist of simple real eigenvalues.
Let us fix the real time parameter $r_0$ sufficiently negative, so that $\Gamma(r_0)$ falls into the neighborhood of the source $W_-$ where Poincaré linearization rules.
For the Schrödinger solutions \eqref{psiw} associated to real heteroclinic orbits $w(r,\cdot)=\Gamma(r):W_-\leadsto W_+$\,, this implies global quasi-periodicity for all $s\in\mathbb{R}$, rather than blow-up for finite $s$.
Again, see \cite{fiestu24} for complete details, and a discussion of some further literature on PDEs in complex time.

Forty years ago, in contrast, pioneering work \cite{Masuda1,Masuda2} by Kyûya Masuda focused on orbits $w$ of the purely quadratic variant $w_r=w_{xx}+w^2$ of \eqref{PDEw}, in any space dimension.
Under Neumann boundary conditions, elliptic maximum principles imply that the spatially homogeneous solution $W=0$ is the only real equilibrium.
Real solutions $w(r,x)$ starting at initial profiles $w(0,x)=w_0(x)>0$, for $r=0$, then blow up in finite real time $0\leq r\nearrow r^*<\infty$.
For almost homogeneous real initial conditions $ w_0$\,, Masuda was then able to circumnavigate blow-up, at real time $r=r^*(w_0)$, via sectorial detours venturing into complex time $t=r+\mi s$.
See \cite{Masuda2} for proofs of the earlier announcement \cite{Masuda1}.
Notably, Masuda also cautioned that the detours via positive and negative imaginary parts $s$ agree in their real-time overlap after blow-up, if \emph{and only if} $w_0$ is spatially homogeneous.

\subsection{ODE setting}\label{ODEset}

In our present ODE study we explore the very special role of quadratic nonlinearities like $w^2$ or $w^2-1$ in \eqref{PDEw}.
Instead of the exceedingly demanding subject of PDEs in complex time, we limit ourselves to a detailed study of the ODE case of homogeneous solutions $w$ which do not depend on the spatial variable $0<x<\ell$, at all.
The diffusion terms $w_{xx}\,,\psi_{xx}$ in \eqref{PDEw}, \eqref{PDEpsi} then drop out, and we obtain quadratic ODEs, in real and imaginary time; see \eqref{ODEw2}, \eqref{ODEpsi2} below.
More generally, we consider scalar complex ODEs with polynomial nonlinearities $f$ of degree $d$, i.e.
\begin{equation}
\label{ODEw}
\dot{w}=f(w)=(w-e_0)\cdot\ldots\cdot(w-e_{d-1})=w^d+\ldots+f_0\,,
\end{equation}
for $w=w(t)$ in complex time $t=r+\mi s$.
Any polynomial ODE of degree $d$ may be brought to this univariate form with normalized top coefficient, of course, by a suitable complex scaling of $w$.
We assume the $d$ complex zeros $e_j$ of $f$ to be simple.
In other words, $W=e_j$ are the $d$ distinct equilibria of \eqref{ODEw}.

Although some of our questions (and answers) seem to be new, much of our analysis proceeds at almost textbook level.
For the convenience of our readers we summarize some elementary background from complex analysis.
The novelty of our questions arises because we will distinguish a real time direction, as motivated by their PDE background.

For an expert presentation of much background on complex analysis of ODEs see the beautiful book by Ilyashenko and Yakovenko \cite{Ilyashenko}.
For given initial conditions $w(0)=w_0$\,, let $\Phi^t(w_0):=w(t)$ denote the local solution flow of the (not necessarily scalar or polynomial) ODEs $\dot w(t)=f(w)$, with holomorphic vector fields $f$.
Then the flow property 
\begin{equation}
\label{flow}
\Phi^{t_2}\circ\Phi^{t_1}= \Phi^{t_1+t_2}, \qquad \Phi^0=\mathrm{id},
\end{equation}
holds, for any argument $w_0$ and, locally, for all $t_1,t_2\in\mathbb{C}$ such that the (small) closed complex parallelogram spanned  by $t_1,t_2$ in $\mathbb{C}$ is contained in the domain of existence of the local semiflow $\Phi^t(w_0)$.
This follows from the ODE and Cauchy's theorem.
Locally, and for fixed $t$, \eqref{flow} implies that the flow maps $w_0\mapsto \Phi^t(w_0)$ are biholomorphic.
Indeed, the holomorphic local inverse of $\Phi^t$ is $\Phi^{-t}$.
Moreover, and for example, the local flow $\Phi^t$ in real time $t=t_1=r$ commutes with the local flow in imaginary (Schrödinger) time $t=t_2=\mi s$.

Explicit solution of \eqref{ODEw}, by standard separation of variables, makes the ODE look deceptively innocent and completely classical.
However, the precise details of blow-up $|w|\nearrow\infty$ in finite complex time $t$, and ``circumnavigation of blow-up'' in the Masuda sense, are less classical topics.
Our PDE motivation, however, prods us to keep track of the real axis, as a distinguished direction of time.
The PDE paradigms, for example, should have convinced us to pay attention to the global dynamics of complex-time solutions $w(t)$ along real lines $r\mapsto t=r+\mi s_0$, even in the ODE case.
Along these lines, we do not favor complex foliations, which so elegantly blur the distinction between real and imaginary time.
By closer scrutiny, we hope to convince the reader that even a ``trivial'' ODE like \eqref{ODEw} does hold some interest, today.

\subsection{The quadratic ODE}\label{ODE2}
Let us first explore the simplest nonlinear case $d=2$ in \eqref{ODEw}.
To be specific, consider $e_0:=1,\ e_1:=-1$.
We then obtain the scalar quadratic ODEs
\begin{align}
    \label{ODEw2}
    \dot w &=f(w)=w^2-1, \\
    \label{ODEpsi2}
    \mi\dot\psi &=f(\psi)= \psi^2-1\,.
\end{align}
Here local analyticity of $w,\psi$ in complex time is obvious.
Note \emph{time reversibility} of the Schrödinger variant \eqref{ODEpsi2}:  $\psi(-s)$ is a solution, if and only if the complex conjugate $\overline\psi(s)$ is.

\begin{figure}[t]
\centering \includegraphics[width=0.86\textwidth]{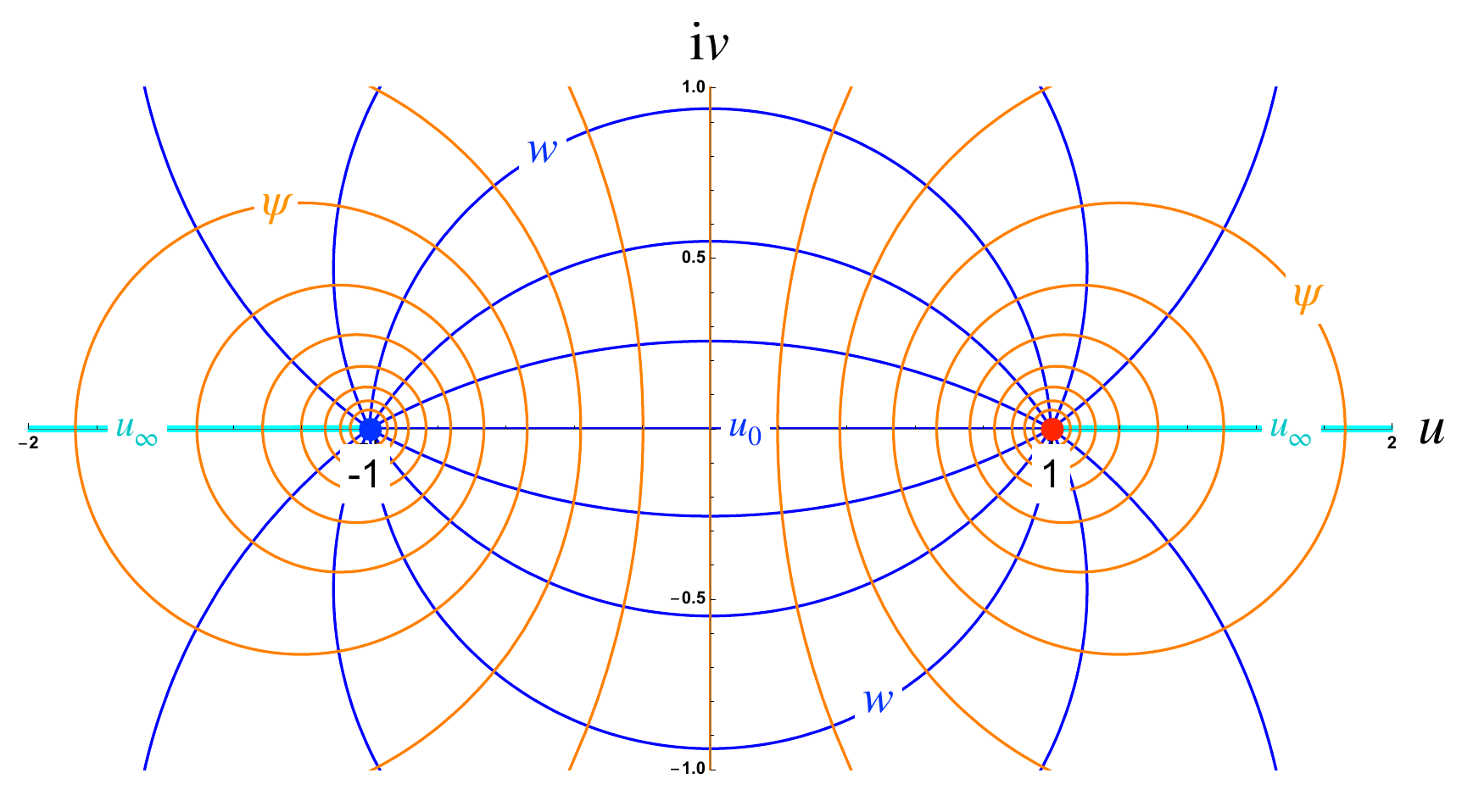}
\caption{\emph{
Phase portrait of quadratic complex ODEs \eqref{ODEw2} -- \eqref{ODEuv2} with equilibria $w=W_\mp=\pm1$ and the heteroclinic orbit $\Gamma:W_-\leadsto W_+$\,.
Orbits $r\mapsto\Gamma(r+\mi s)=u+\mi v$ of \eqref{ODEuv2}, \eqref{ODEw2}, in real time $r$, are circular arcs (blue). 
They are heteroclinic from source  $W_-=+1$ (red dot) to sink $W_+=-1$ (blue dot).
The $\pi$-periodic orbits $s\mapsto\psi(s)=\Gamma(r-\mi s)=u+\mi v$ of \eqref{ODEpsi2}, in contrast, are full circles (orange) in $\mathbb{C}$, each surrounding one of the two equilibria $W_\mp=\pm 1$.
By \eqref{flow}, the flows in real and imaginary time commute.
Therefore the orange circles also serve as isochrones which globally synchronize the blue heteroclinics.
Conversely, the blue circle segments are isochrones which globally synchronize the orange periodic solutions.
For any fixed time $t$, the holomorphic flow map $w_0\mapsto \Phi^t(w_0)$ is conformal, i.e. angle-preserving. 
In particular, the blue and orange circle families are mutually orthogonal.\\
The imaginary $v$-axis (orange) blows up in finite time $s^*=\pm\pi/2$; see \eqref{s*}.
On the Riemann sphere $\widehat{\mathbb{C}}$\,, this is just a longitude circle of period $\pi$ through the South Pole at $w=0$ and the North Pole at $w=\infty$.
The real $u$-axis features the real heteroclinic orbit $u_0: +1 \leadsto -1$, in blue, as well as the two cyan parts of the blow-up heteroclinic $u_\infty:+1\leadsto\infty$ and the blow-down heteroclinic $u_\infty:\infty\leadsto -1$.
On the Riemann sphere $\widehat{\mathbb{C}}$\,, these three line segments combine to another longitude circle, perpendicular to the first at the polar intersections 0 and $\infty$.
In conclusion, the blue bounded heteroclinic orbit $u_0$ on the real $u$-axis, in real time, gives rise to finite time blow-up and blow-down on the orange imaginary $v$-axis, in imaginary time, when started at $u=v=0$.
The complex viewpoint also reveals how the blue real heteroclinicity $u_0$ and the two cyan real blow-up/blow-down segments $u_\infty$\,, in real time, become three aspects of one and the same underlying trajectory $\Gamma$, in complex time $t=r+\mi s$.\\
}}
\label{fig1}
\end{figure}

For $w=u+\mi v$ we obtain the equivalent real system
\begin{equation}
\label{ODEuv2}
\begin{aligned}
    \dot{u} &= u^2-v^2-1\,,  \\
    \dot{v} &= \quad 2uv\,.
\end{aligned}
\end{equation}
See figure \ref{fig1} for a phase portrait.
The equilibria are the attractor, or \emph{sink}, equilibrium $W_+=e_1=-1$ (blue dot) and the repellor, or \emph{source}, $W_-=e_0=+1$ (red dot), with linearizations $f'(W_\mp)=\pm2$.
All nonstationary orbits (blue) are heteroclinic $\Gamma: W_-\leadsto W_+$.
The real $u$-axis $\{v=0\}$ is invariant and contains the monotonically decreasing heteroclinic orbit $u_0:1\leadsto -1$.
The two unbounded real orbits $u_\infty$ on the $u$-axis (cyan) look like an exception, at first: they blow up or blow down at $u=\pm\infty$ in finite positive or negative time, just as in the purely quadratic Masuda case $e_0=e_1=0$.

Unlike the general PDE case \eqref{PDEw}, the scalar ODE \eqref{ODEw2} can be regularized on the \emph{Riemann sphere} $w\in\widehat{\mathbb{C}}:=\mathbb{C}\cup\{\infty\}$\,.
Indeed, just note equivariance of \eqref{ODEw2} under the involutive automorphism $w\mapsto 1/w$.
In other words, $w(t)$ is a solution, if and only if $1/w(t)$ is.
In particular, the two blow-up pieces of $u_\infty$ can be joined to indicate yet another single heteroclinic orbit, on the Riemann sphere.
Indeed, cyan $u_\infty$ is just the involutive copy of blue $u_0$\,, on $\widehat{\mathbb{C}}$\,.

Because $\dim w = 1$, Poincaré linearization near $W_-=1$ shows that the complex heteroclinic orbit $s\mapsto\Gamma(r+\mi s)$ is periodic of minimal period $p=\pi$ in imaginary time, for any fixed $r\in\mathbb{R}$.
In other words, the meromorphic heteroclinic orbit $t\mapsto\Gamma(t)$ provides a bi-holomorphic map between complex cylinders,
\begin{equation}
\label{ODEhol}
\Gamma:\ \mathbb{C}/\pi\mi\mathbb{Z} \ \rightarrow\  \widehat{\mathbb{C}}_2:=\widehat{\mathbb{C}}\setminus\{\pm 1\}\,.
\end{equation}

By the PDE results of \cite{fiestu24}, blow-up to $\psi(s)=\infty$ must occur in the Schrödinger variant \eqref{psiw}, \eqref{ODEpsi2}, for some $r=r_0$ and at some $s=s^*$.
Bi-holomorphy \eqref{ODEhol} implies that $r_0\in\mathbb{R}$ and $s^*\in\mathbb{R}/\pi\mathbb{Z}$ are unique.
Here real $\Gamma(r):=u_0(r)\in (-1,1)$ has to track the blue bounded decreasing real heteroclinic orbit. 
By time reversibility of \eqref{ODEpsi2} in $s$, blow-up must therefore occur at 
\begin{equation}
\label{s*}
s^*\equiv\pi/2 \mod\pi,
\end{equation}
for some $\psi_0=\Gamma(r_0)$.
Also, $\Gamma(r\pm\pi\mi/2)\in\mathbb{R}$, except at the blow-up point $r=r_0$.

Without loss of generality, let us fix $\Gamma(0)=0$.
Explicitly, $v=0$ in \eqref{ODEuv2} then identifies the biholomorphic cylinder map \eqref{ODEhol} as
\begin{equation}
\label{Gammaexplicit}
\Gamma(t)=-\tanh t \ \in\  \widehat{\mathbb{C}}\setminus\{\pm 1\}\,,
\end{equation}
first for real $t=r\in\mathbb{R}$, i.e. for $u_0:1\leadsto -1$, and then for all $t\in\mathbb{C}$, by maximal analytic continuation.
This identifies $r_0=0$ and the Schrödinger blow-up solution
\begin{equation}
\label{psiexplicit}
\psi(s)=-\mi\tan s\,.
\end{equation}
Complex time shifts of $\Gamma$, in \eqref{ODEhol}, by the flow $\Phi^t$ of \eqref{flow}, induce biholomorphic automorphisms of the Riemann sphere $\widehat{\mathbb{C}}$.
The maps $\Phi^t$ are therefore fractional linear Möbius transformations, which fix the equilibria $\pm 1$. 
In \eqref{Gammaexplicit}, commutativity \eqref{flow} on the heteroclinic orbit $\Gamma$ therefore amounts to the elementary addition theorems for the hyperbolic tangent.

In passing, we note that the real and imaginary parts of the inverse function $t+\mi s=-\mathrm{arctanh}\,w$ of the solution \eqref{Gammaexplicit} are first integrals of \eqref{ODEpsi2} and \eqref{ODEw2}, respectively.
In figure \ref{fig1} their level curves are colored orange and blue/cyan.
Specifically, the blue/cyan heteroclinic orbits of \eqref{ODEuv2}, \eqref{ODEw2}, in real time $\mathbb{R}\ni r\mapsto t=r+\mi s$ and for fixed $s$, foliate the cylinder $\widehat{\mathbb{C}}\setminus\{\pm 1\}$ into segments of invariant Euclidean circles through $w=\pm1$, with centers on the imaginary $v$-axis.
The orange $\pi$-periodic orbits of \eqref{ODEuv2}, in the Schrödinger variant \eqref{ODEpsi2} of imaginary time $t=\mi s$, define a perpendicular foliation into circles with centers $w=a\in\mathbb{R},\ 1<|a|\leq\infty$ on the real $u$-axis.
In other words, their radii $0<\sqrt{a^2-1}\leq\infty$ are defined by the tangents from $a\in\mathbb{R}$ to the complex unit circle.
The straight axes themselves are "circles with centers at infinity", as already Cusanus has remarked \cite{Cusanus}.

In electrostatics, the blue/cyan and orange circle families of figure \ref{fig1} are known as field lines and potential levels of an electrical dipole with charges $\pm1$ located at $W=\pm1$.
Similarly they illustrate flow lines and potential of a static, incompressible, irrotational planar fluid flow with source and sink at $W=\pm1$. 
See for example \cite{Marsden}, ch 5.3.

We can also trivialize the cylinder map \eqref{ODEhol} by the fractional linear, biholomorphic automorphism $w\mapsto (z+1)/(z-1)$ of the Riemann sphere $\widehat{\mathbb{C}}$\,.  
This Möbius transformation maps the equilibria $W=\pm 1$ to $\infty$ and 0, respectively. 
It also globally linearizes the quadratic complex ODE $\dot{w}=w^2-1$, alias \eqref{ODEuv2}, to become $\dot{z}=-2z$.
In $z$, the heteroclinic orbits of \eqref{ODEw2} are inward radial, and the periodic orbits of \eqref{ODEpsi2} are circles of constant radius $|z|$.
Circles and lines in figure \ref{fig1} illustrate how Möbius transformations, just as circle inversions, map collections of circles and lines to circles and lines, preserving angles of intersection.
For general scalar ODEs $\dot w=f(w)$, we continue this discussion in section \ref{RRes}.

\subsection{Prescribing linearizations}\label{f'Res}

In general, we have assumed all zeros $e_j$ of $f$ to be simple.
Standard separation of variables then provides solutions $w(t)$ of the polynomial ODE \eqref{ODEw} with initial condition $w(0)=w_0$ by explicit integration of the differential form $\omega=dw/f(w)$:
\begin{equation}
\label{sovw}
t \,=\, \int_{w_0}^{w(t)}\,\frac{dw}{f(w)} \,=\, c(w_0) \,+ \, \sum_{j=0}^{d-1}\,\tfrac{1}{f'(e_j)} \,\log(w(t)-e_j)\,.
\end{equation}
The integration constant $c(w_0)$ annihilates the sum at $t=0$, of course.
The second equality follows from the partial fraction decomposition
\begin{equation}
\label{partialfrac}
\frac{1}{f(w)} \,=\, \sum_{j=0}^{d-1}\,\eta_j\cdot\frac{1}{w-e_j}\;,\qquad\mathrm{with} \ \eta_j:=\frac{1}{f'(e_j)}\,.
\end{equation}
Nonuniqueness of the logarithms by addition of $2\pi\mi\mathbb{Z}$ corresponds to different choices of complex integration paths from $w_0$ to $w(t)$.
Indeed, the \emph{residue} of $1/f(w)$ at the simple pole $w=e_j$ coincides with the reciprocal $\eta_j$ of the linearization $f'(e_j)$ at the equilibrium $e_j$\,.
To understand the $d$ resulting \emph{periods} 
\begin{equation}
\label{Tj}
T_j:=2\pi\mi\eta_j\,,
\end{equation}
of $w(t)$, we therefore attempt to prescribe the coefficients $\eta_j$\,.

For any degree $d\geq2$, one constraint arises from holomorphic integrability of $1/f$ in $z:=1/w$ near $w=\infty,\ z=0$, for $w$ in the punctured Riemann sphere
\begin{equation}
\label{Cd}
\widehat{\mathbb{C}}_d:=\widehat{\mathbb{C}}\setminus\{e_0,\ldots,e_{d-1}\}\,.
\end{equation}
Indeed, we may fix time $t=0$ at $w_0=\infty$ and solve ODE \eqref{ODEw} for $z=z(t)=1/w(t),\ z(0)=0$ by separation of variables, analogously to \eqref{sovw}. 
In other words, we rewrite the differential form $\omega=dw/f(w)$ in the new coordinate $z=1/w$ as
\begin{equation}
\label{sowom}
\omega \,=\, dw/f(w) \,=\, -\frac{z^{d-2}}{z^df(1/z)}\,dz
  \,=\, - (1+\ldots+f_0z^d)^{-1}\,z^{d-2} dz\,.
\end{equation}
In particular, the differential form $\omega$ is locally holomorphic at $z=0$, alias $w=\infty$, with locally convergent power series
\begin{equation}
\label{sovz}
\begin{aligned}
   t &\,=\, \int_\infty^{w}\,\omega\,=\, -\tfrac{1}{d-1}z^{d-1}(1+c_1z+\ldots)
\,,
\end{aligned}
\end{equation}
for some coefficients $c_1,\ldots$ .
Holomorphic integration \eqref{sovz} also implies the constraint
\begin{equation}
\label{sumf'eta}
\sum_{j=0}^{d-1}\,\frac{1}{f'(e_j)}\,=\,\sum_{j=0}^{d-1}\,\eta_j\,=\,0\,.
\end{equation}
Indeed, consider any left oriented, closed Jordan curve $\gamma\subset\mathbb{C}$ which is large enough to contain all zeros $e_j$ of $f$ in its interior.
Up to a factor of $2\pi\mi$, the residue theorem then identifies the sums as the integral of $\omega=dw/f(w)$ over $\gamma$.  
Shrinking the loop $\gamma$ to $z=0$, as in \eqref{sovz}, then proves constraint \eqref{sumf'eta}.

Observing this constraint we can now identify possible assignments of $f'(e_j)$.

\begin{thm}\label{f'thm}
Let $d\geq2$. For $0\leq j<d$, prescribe values $0\neq\eta_j\in\mathbb{C}$.
Assume
\begin{equation}
\label{sumseta}
\sum_{j=0}^{d-1}\,\eta_j\,=\,0\,,\quad\mathrm{but\ also}\quad
 \sum_{j\in J}\,\eta_j\,\neq\,0
\end{equation}
for any nonempty subset $\emptyset\neq J \subsetneq\{0,\ldots,d-1\}$.\\
Then there exists some univariate complex polynomial $f$ of degree $d$, as in \eqref{ODEw}, with derivatives $f'(e_j)=1/\eta_j\neq 0$ prescribed by $\eta_j$\,, at all $d$ simple zeros $e_j$\,.
\end{thm}

Theorem \ref{f'thm} will be proved in section \ref{f'Pf}.
Unlike common interpolation problems from Lagrange to Hermite-Birkhoff \cite{Interpol1, Interpol2, Interpol3}, this result is highly nonlinear.
The locations of the zeros $e_j$ are not given.
Rather, they result from the prescribed reciprocal derivatives $\eta_j$\,.
The resulting locations are in fact not unique.
For example, we may translate all zeros $e_j$ by a fixed complex constant $c$, without affecting the derivatives.

The sum over all $0\leq j<d$ in \eqref{sumseta} must vanish, by \eqref{sumf'eta}.
The additional nondegeneracy condition  on subsets $J$ is therefore equivalent to that same condition on all partitions of $0\leq j<d$ into complementary nonempty subsets $J,J^c$.
For degrees $d=2,3$ and simple zeros, in particular, nondegeneracy holds automatically.
For degrees $d\geq 4$, it is violated by certain polynomials.
Degeneracy at degree $d=4$ occurs, e.g., for the cyclotomic polynomial $f(w)=w^4-1$ and $J,J^c$ designating the pairs of equilibria $W=\pm 1$ and $W=\pm\mi$.
See section \ref{Cyclotomic} and figure \ref{fig2} for further discussion of this cyclotomic case.

\subsection{The period map}\label{PeriodT}

We begin to address global aspects of complex polynomial ODEs $\dot{w}=f(w)$ of degree $d\geq2$; see  \eqref{ODEw}.
Our presentation paraphrases \cite{Forster}, §§ 6--10.
See also section \ref{Apx} for a summary of terminology.

Puncturing the Riemann sphere $\widehat{\mathbb{C}}$ at the simple equilibria  $e_0,\ldots,e_{d-1}$, we obtain the domain $w\in\widehat{\mathbb{C}}_d$\,. See \eqref{Cd}.
Separation of variables \eqref{sovw} and regularity \eqref{sovz} of $t$ at $w=\infty\in\mathbb{C}_d$ then allow us to integrate the ODE via the holomorphic differential form $\omega=dw/f(w)$ on $w\in\widehat{\mathbb{C}}_d$\,.
This rewrites \eqref{sovw}, \eqref{sovz} more concisely as 
\begin{equation}
\label{tw}
t=\int_\infty^w\, \omega\,.
\end{equation}
However, the integral of $\omega$ is multivalued.
In fact, the integral depends on the complex integration path which we choose in $\widehat{\mathbb{C}}_d$, from $\infty$ to $w$.
To capture this ambiguity, we define the \emph{period map}
\begin{equation}
\label{T}
\begin{aligned}
\mathcal{P}: \pi_1(\widehat{\mathbb{C}}_d) &\rightarrow \mathbb{C} \\
							\gamma &\mapsto \int_\gamma\,\omega
\end{aligned}
\end{equation}
on the fundamental group $\pi_1(\widehat{\mathbb{C}}_d)$.
Here $\gamma$ denotes any closed loop through $w=\infty$ in $\widehat{\mathbb{C}}_d$\,,
which need not be an ODE trajectory.
The integral only depends on the homotopy class of $\gamma$ in $\pi_1$\,.
Since punctured spheres remain path connected, we may replace the base point $w=\infty\in\widehat{\mathbb{C}}_d$ by any other base point.
By \eqref{pi1Cd}, in fact, the fundamental group $\pi_1(\widehat{\mathbb{C}}_d)\cong\mathbb{F}_{d-1}$ coincides with the free group on $d-1$ generators.

To capture the ambiguity of the integration path, we can now rewrite \eqref{tw} more appropriately as 
\begin{equation}
\label{alltw}
t\ \in\ \mathrm{range}\,\mathcal{P}+ \int_\infty^w\, \omega
\end{equation}
along arbitrary paths in $\widehat{\mathbb{C}}_d$ from $\infty$ to $w$.
Any such $t$ can in fact be realized, by an appropriate choice of the integration path.

To determine the structure of $\mathrm{range}\,\mathcal{P}$, we first note that $\mathcal{P}$ is a group homomorphism into $(\mathbb{C},+)$, just by concatenation of cycles $\gamma$.
The kernel of $\mathcal{P}$ identifies those cycles $\gamma$ which do not affect integration \eqref{tw}.
The commutator subgroup $\pi_1'$ of $\pi_1$\,, i.e. the normal subgroup of $\pi_1$ with maximal Abelian factor $\pi_1/\pi_1'$\,, is generated by the commutator elements $\gamma_1\gamma_2\gamma_1^{-1}\gamma_2^{-1}$ of $\pi_1$\,.
Additivity of the integral \eqref{T} therefore implies
\begin{equation}
\label{kerTpi'}
\ker \mathcal{P} \geq \pi_1' \,.
\end{equation}
The commutator $\pi_1'$ is also a free group, but of infinite rank for $d\geq 3$.

The homomorphism theorem for $\mathcal{P}$, on the other hand, implies 
\begin{equation}
\label{Thom}
\mathrm{range}\,\mathcal{P}\, \cong\, \pi_1/\ker \mathcal{P} \,\leq\, \pi_1/\pi_1' \,\cong\, H_1(\widehat{\mathbb{C}}_d,\mathbb{Z}) \,=\, (\mathbb{Z}^{d-1},+)\,,
\end{equation}
for some $0<k<d$.
Indeed, the first homology $H_1$ coincides with the Abelianization $\pi_1/\pi_1'\cong (\mathbb{Z}^{d-1},+)$ of the fundamental group $\pi_1(\widehat{\mathbb{C}}_d)\cong\mathbb{F}_{d-1}$.

The subgroup $\mathrm{range}\,\mathcal{P}$ of the free Abelian group $\mathbb{Z}^{d-1}$ is again a free $\mathbb{Z}$-module, because the ring of integers is torsion free:
\begin{equation}
\label{Trange}
\mathrm{range}\,\mathcal{P}\,=\,\langle T_1,\ldots,T_{d-1}\rangle_\mathbb{Z}\,\cong\,(\mathbb{Z}^k,+)\,,
\end{equation}
for some $0<k<d$.
In case $k<d-1$, however, note that the \emph{elementary periods}
\begin{equation}
\label{Tjres}
T_j\,:=\,\int_{\gamma_j}\omega= 2\pi\mi\cdot \mathrm{Res}\,(1/f,e_j)\,=\, 2\pi\mi/f'(e_j)\,=\, 2\pi\mi\,\eta_j
\end{equation}
become linearly dependent over $\mathbb{Z}$ or, equivalently, over $\mathbb{Q}$; see also \eqref{sovw}, \eqref{Tj}.
The range of the period map $\mathcal{P}$, of course, just captures the ambiguity \eqref{alltw} generated by the multi-valued logarithms in the explicit solution \eqref{sovw}.
For the topological closure $\mathbb{B}$ of $\mathrm{range}\,\mathcal{P}$, see the options \eqref{closT} below.

\subsection{The Riemann surface of the flow}\label{RRes}

In this section we further explore classical global aspects of the explicit solutions \eqref{sovw}, \eqref{sovz} for the complex polynomial ODE \eqref{ODEw}.
At $w=\infty,\ z=1/w=0$ in the punctured Riemann sphere $w\in\widehat{\mathbb{C}}_d:=\widehat{\mathbb{C}}\setminus \{e_0,\ldots,e_{d-1}\}$, the expansion \eqref{sovz} has preserved local analyticity.
In general, we denote the solution set of the multivalued relation \eqref{alltw} as 
\begin{equation}
\label{Rdef}
\mathcal{R} := \{(w,t)\in \widehat{\mathbb{C}}_d\times\mathbb{C}\,\vert\, \eqref{alltw}\ \mathrm{holds}\, \}\,.
\end{equation}

To describe the complex structure of $\mathcal{R}$, we freely use some terminology on Riemann surfaces in the following theorem.
For details, following \cite{Forster, Hartshorne, Jost, Lamotke} for example, see the appendix in section \ref{Apx}, and our proof in section \ref{RPf}.

\begin{thm}\label{Rthm}
Let $f$ denote any univariate complex polynomial of degree $d\geq 3$ with simple zeros $e_0,\ldots,e_{d-1}$\,; see ODE \eqref{ODEw}.
Then the following four statements hold true.
\begin{enumerate}[(i)]
\item The set $\mathcal{R}$ defined in \eqref{Rdef} coincides, as a set, with the Riemann surface of the maximal analytic continuation of the integral $t$ of the differential form $\omega=dw/f(w)$ on $w\in\widehat{\mathbb{C}}_d$ as defined in \eqref{tw}; see also \eqref{alltw}. 
The starting germ of the analytic continuation is the local Taylor series \eqref{sovz} for $t$ as a function of $z=1/w$, at $w=\infty,\ t=0$.
For any $w_0\in\widehat{\mathbb{C}}_d$\,, the local Taylor series $t=t(w;w_0)$ at $w_0$ coincide, except for their constant term $t_0$ with $(w_0,t_0)\in\mathcal{R}$.
In other words, the constant difference $T=t_1(w;w_0)-t_2(w;w_0)$ of any two local expansions $t_1,t_2$ of $t$, at the same $w_0\in\widehat{\mathbb{C}}_d$\,, is a period $T\in\mathrm{range}\,\mathcal{P}=\langle T_1,\ldots,T_{d-1}\rangle_\mathbb{Z}\cong\mathbb{Z}^k$ of the period map $\mathcal{P}$; see \eqref{Trange}.
The fibers $t_0+\mathrm{range}\,\mathcal{P}$ over each $w_0$ are endowed with the discrete topology, on $\mathcal{R}$.

\item The two Riemann surfaces $(w,t)\in\mathcal{R}$ and $w\in\widehat{\mathbb{C}}_d$ are hyperbolic, with the complex upper half plane $\mathbb{H}$ as the universal cover of both:
\begin{equation}
\label{HRCd}
\mathbb{H}\xrightarrow{\mathbf{p}}\mathcal{R}\xrightarrow{\mathbf{q}_w}\widehat{\mathbb{C}}_d\,.
\end{equation}
The projection $\mathbf{q}_w(w,t):=w$ and the universal coverings $\mathbf{p},\,\mathbf{p}\circ \mathbf{q}_w$ are unbranched, unlimited, normal covering maps.
In particular, the topology of any fiber $\mathbf{q}_w^{-1}(w_0)$ of $\mathcal{R}$ is discrete.
The corresponding deck transformation groups are
\begin{align}
\label{deckHRCd}
    \mathrm{deck}(\mathbf{p}\circ \mathbf{q}_w)\,&\cong\, \pi_1(\widehat{\mathbb{C}}_d) \,\cong\, \mathbb{F}_{d-1}\,;   \\
\label{deckHR}
    \mathrm{deck}(\mathbf{p})\,\,\   &\cong\,\, \pi_1(\mathcal{R}) \;\cong\,\ker\,\mathcal{P}\,\trianglelefteq\, \mathbb{F}_{d-1}\,;   \\
\label{deckRCd}
    \mathrm{deck}(\mathbf{q}_w)\; &\cong\, \mathbb{F}_{d-1}/\ker\,\mathcal{P} \,\cong\, \mathrm{range} \,\mathcal{P}\,=\,\langle T_1,\ldots,T_{d-1}\rangle_\mathbb{Z}\,\cong\,\mathbb{Z}^k\,,
\end{align}
for some $0<k<d$.
Note how kernel and range of the period map $\mathcal{P}$ from \eqref{T} describe the deck groups, alias the fibers, of the first and second covering map in \eqref{HRCd}, respectively.

\item Global integration \eqref{alltw} amounts to the covering projection $\mathbf{q}_t(w,t):=t$,
\begin{equation}
\label{RC}
\mathbf{q}_t:\ \mathcal{R}\rightarrow\mathbb{C}\,.
\end{equation}
Branching points are all periods $t=T\in\mathrm{range}\,\mathcal{P}$, each of branching multiplicity $d-1$, with associated ramification points $(\infty,t)\in\mathcal{R}$.
The topological closure $\mathbb{B}$ of the branching set $\mathrm{range}\,\mathcal{P}\subset\mathbb{C}$ depends on the generators $T_j=2\pi\mi\eta_j$\,. 
This makes the closure $\mathbb{B}$ real linear equivalent to one the following five options
\begin{equation}
\label{closT}
\mathbb{Z},\quad \mathbb{Z}\times \mathbb{Z},\quad \mathbb{R},\quad \mathbb{R}\times \mathbb{Z},\quad \mathbb{C}\cong \mathbb{R}^2\,.
\end{equation}
With $\eta_j:=1/f'(e_j)$ realizable at least as prescribed in theorem \ref{f'thm}, all five cases actually arise for suitable polynomial nonlinearities $f$.
In the topology of \,$\widehat{\mathbb{C}}_d\times\mathbb{C}$, the description \eqref{Rdef} is an embedding of the Riemann surface $\mathcal{R}$ in the first two cases, only.
These are also the only cases where the set of branching points is discrete and, therefore, where the projection $\mathbf{q}_t$ qualifies as an unlimited branched covering. 
In the remaining three cases, branching points are densely accumulating in $\mathbb{B}$, and we just obtain an immersion $(w,t)\in\mathcal{R}\subset\widehat{\mathbb{C}}_d\times\mathbb{C}$ of the densely accumulating local $t$-sheets defined by $t=t(w;w_0)+T$, for periods $T\in\mathrm{range}\,\mathcal{P}$.

\item The local flow $\Phi^t$ of ODE \eqref{ODEw} on $\widehat{\mathbb{C}}_d\setminus\{\infty\}$ lifts to local flows on $\mathcal{R}':=\mathcal{R}\setminus(\{\infty\}\times\ker\mathcal{P})$ and $\mathbb{H}':=\mathbf{p}^{-1}(\mathcal{R}')\subsetneq\mathbb{H}$ by the covering maps $\mathbf{q}_w$ and $\mathbf{p}\circ\mathbf{q}_w$\,, respectively.
The lifted flows are equivariant under the deck groups $\mathrm{deck}(\mathbf{q}_w)$ and $\mathrm{deck}(\mathbf{p}\circ \mathbf{q}_w)$ of $\mathbf{q}_w$ and $\mathbf{p}\circ\mathbf{q}_w$\,; see \eqref{deckHRCd}, \eqref{deckRCd}.
\end{enumerate}
\end{thm}

For $d\geq 3$, we cannot lift the local flow $\Phi^t$ of our original ODE \eqref{ODEw} to include the branch point $w=\infty$.
On $t\in\mathbb{C}$ in fact, i.e. after projection by $\mathbf{q}_t$\,, that flow simply amounts to a complex time shift.
At the ramification points $(w=\infty,t_0)\in\mathcal{R}$ of branch points $t_0$ of multiplicity $d-1$, in theorem \ref{Rthm}\emph{(iii)}, however, expansion \eqref{sovz} shows how $d-1$ different real-time trajectories $r\mapsto z(t_0+r)$ cross each other, at $r=0$. 
See also section \ref{Cyclotomic} and figure \ref{fig2} below.
In conclusion, indeed, real time shift in $t$ cannot lift to a real-time flow on the Riemann surface $\mathcal{R}$, by this construction.
We will resort to smooth (but not holomorphic) regularization and Poincaré compactification of the real-time flow, instead, in the next section.

In the quadratic case $d=2$ of section \ref{ODE2}, we have already seen the relevant modifications of theorem \ref{Rthm} at work. 
For the parabolic cylinder $w\in\widehat{\mathbb{C}}_2$\,, we only have to substitute the hyperbolic upper half plane $\mathbb{H}$ by the universal cover $\mathbb{C}$.
This implies the Abelian fundamental group $\pi_1(\widehat{\mathbb{C}}_2)\cong\mathbb{F}_1\cong\mathbb{Z}\cong H_1(\widehat{\mathbb{C}}_2)$; see \eqref{deckHRCd}.
The period map features a single purely imaginary period $T_1=2\pi\mi/f'(e_1)=-\pi\mi$. 
See \eqref{Tjres}, \eqref{deckRCd}.
In particular, $\pi_1(\mathcal{R})\cong\ker \mathcal{P}=\{0\}$ is trivial.
Therefore $\mathcal{R}$ is biholomorphically equivalent to the simply connected universal cover $\mathbb{C}$, and the unbranched projection $\mathbf{q}_t: \mathcal{R}\rightarrow\mathbb{C}$ of \eqref{RC} is biholomorphic as well.
The composition $\mathbf{q}_w\circ\mathbf{q}_t^{-1}$ of \eqref{HRCd} and \eqref{RC} provides a solution $\mathbb{C}\ni t\mapsto w(t)$ of ODE \eqref{ODEw}. The initial condition is $w=\infty$ at time $t=0$, and the purely imaginary time period is $T_1$\,.

For all degrees $d\geq 2$, of course, these results confirm and extend the PDE result of \cite{fiestu24} in the simplest ODE case of spatially homogeneous solutions; see section \ref{Para}.
Indeed, any real-time heteroclinic orbit $w=\Gamma:W_-\leadsto W_+$ between any two hyperbolic equilibria $e_j$ has to be contained in $\widehat{\mathbb{C}}_d$\,.
In the case $d\geq 3$ of theorem \ref{Rthm}, let $t^*=r_0+\mi s^*$ denote any branch point of the branched covering $\mathbf{q}_t:\mathcal{R}\rightarrow\mathbb{C}$ such that $\lvert s^*\rvert>0$ is minimal, for that choice of $r_0$\,. 
In the unbranched case $d=2$, we simply pick $t^*$ from $\mathbf{q}_t(\mathbf{p}_w^{-1}(\infty))$.
Then $s\mapsto \Gamma(r_0+\mi s)$ starts at $\Gamma(r_0)$ and blows up (or down) to the ramification point $(w,t)=(\infty,t^*)$, at imaginary Schrödinger time $s=s^*$.

Only in the quadratic Masuda case $d=2$, however, the isolated blow-up at $t=t^*$ can be circumnavigated in complex time.
For $d\geq 3$, unique continuation by circumnavigation fails due to branching. 
In fact, we will see in the next sections, how two basic options for continuation will lead to two different target equilibria $W_+$\,, after circumnavigation, once real time flow is resumed.

\subsection{Poincaré compactification}\label{PoinRes}

\begin{figure}[t]
\centering \includegraphics[width=\textwidth]{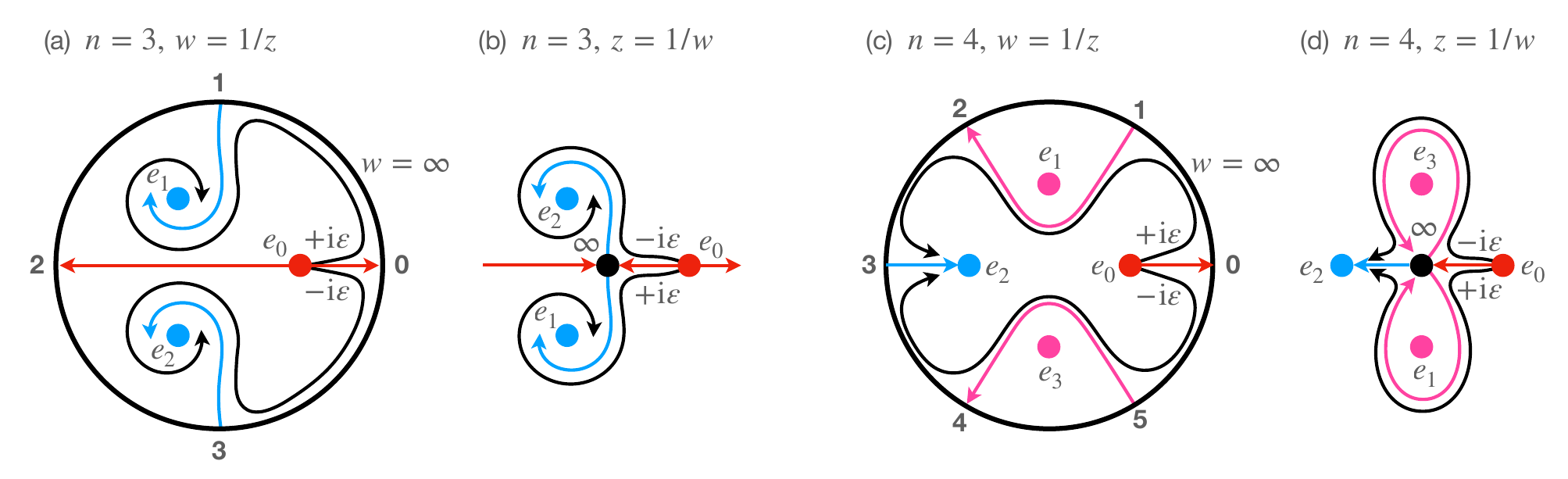}
\caption{\emph{
Schematic phase portraits, in real time, of complex-valued ODEs \eqref{ODEw}, \eqref{gpolyzEuler} -- \eqref{alpha} for cyclotomic vector fields $\dot{w}=w^d-1$; see \eqref{cyclotomic}.
For $d=3$ see $w$ in (a), and $z=1/w$ in (b).
Similarly, $w$ in (c) and $z=1/w$ in (d) refer to $d=4$.
The invariant circle $\rho=|z|=0,\ \alpha\in\mathbb{S}^1$ of the Poincaré compactification at $w=\infty$ is marked black in (a), (c).
Interior equilibria $e_j\in\mathbb{D}$ are stable sinks (blue), Lyapunov centers (purple) surrounded by periodic orbits, or unstable sources (red).
Unstable blow-down separatrices (blue) emanate from $w=\infty$ at odd-labeled vertices $\mathbf{k}=\mathbf{1},\mathbf{3},\mathbf{5}$.  
Stable blow-up separatrices (red) run towards the even-labeled saddles $\mathbf{k}=\mathbf{0},\mathbf{2},\mathbf{4}$.
In (c), two pairs of interior saddle separatrices coincide (purple) in each of the interior saddle-connections $\mathbf{1}\leadsto \mathbf{2}$ and $\mathbf{5}\leadsto \mathbf{4}$.
In (d), the two purple separatrices become homoclinic to $w=\infty$.\\
The black trajectories marked $\pm\mi\eps$ are complex perturbations of the red real blow-up separatrices $e_0\leadsto \mathbf{0}$, say with initial conditions $w(0)=2\pm\mi\eps$.
For $d=3$ in (a), they closely follow the heteroclinic chains $e_0\leadsto \mathbf{0}\leadsto \mathbf{1}\leadsto e_1$ and $e_0\leadsto \mathbf{0}\leadsto \mathbf{3}\leadsto e_2$, respectively.
The intermediate boundary connections $\mathbf{0}\leadsto \mathbf{1}$ and $\mathbf{0}\leadsto \mathbf{3}$ within the invariant boundary circle $\mathbb{S}^1$ concatenate initial red blow-up to terminal blue blow-down.
In the polar view (b), centered at $w=\infty,\ z=0$, the boundary parts are conflated into $z=0$.
Note the \emph{markedly distinct limits of the two perturbations}, for $\eps\searrow 0$, given by the two distinct remaining blow-up-down concatenations $e_0\leadsto \infty\leadsto e_1$ and $e_0\leadsto \infty\leadsto e_2$\,.\\
In case $d=4$ (d), the initial red blow-up and terminal blue blow-down limits of both perturbations $\pm\eps$ coincide.
Each limit $e_0\leadsto \infty\leadsto\infty\leadsto e_3$ contains an additional interior blow-down-up separatrix $\infty\leadsto\infty$ (purple) which is homoclinic to $z=0$.
The two counter-rotating purple homoclinic separatrices, however, are markedly distinct: their lobes surround the distinct centers $e_1$ and $e_3$, respectively, in opposite direction.
In (c), this is manifest by the two purple interior saddle-connections $\mathbf{1}\leadsto \mathbf{2}$ and $\mathbf{5}\leadsto \mathbf{4}$.
Together with the heteroclinic boundary connections within the black invariant boundary circle $\mathbb{S}^1$ which run between the same saddles, in opposite direction, we obtain two heteroclinic cycles. 
Their interiors are properly foliated by families of nested, synchronously iso-periodic orbits of minimal periods $\mp\pi/2$ around the Lyapunov centers $e_1$ and $e_3$\,; see lemma \ref{perlem} (iv),(v).
}}
\label{fig2}
\end{figure}

Any de-singularization of \eqref{ODEw}, in real time $t=r$, should provide some global real flow on the Riemann sphere $\widehat{\mathbb{C}}$, such that the orbits contain the orbits of \eqref{ODEw} on $w\in\mathbb{C}$.
However, we have already mentioned how branching of multiplicity $d-1$ in theorem \ref{Rthm}\emph{(iii)} leads to trajectory crossings at $w=\infty$, for polynomial degrees $d\geq3$.
For examples of degrees $d=3,4$ see figure \ref{fig2} (b),(d).
To obtain a regular flow near $w=\infty$, and to preserve real time direction, we therefore have to slow down the flow by an Euler multiplier which introduces an artificial equilibrium at $w=\infty$.
However, the Euler multiplier cannot be holomorphic on $\widehat{\mathbb{C}}$ because compactness would render it constant.
In fact, the Euler multiplier will have to be real and nonzero on $w\in\mathbb{C}$, to preserve real-time orbits.

By the Lefschetz theorem on the sphere $\widehat{\mathbb{C}}$ of Euler characteristic $\chi=2$, the local Brouwer degree of the artificial equilibrium $w=\infty$ has to be $\chi-d=-(d-2)<0$, for $d\geq3$.
In particular, the equilibrium $w=\infty$ will be of hyperbolic type, and degenerate for $d>3$.
Therefore we will replace the point $w=\infty$ by a ``circle $\mathbb{S}^1$ at infinity''.
This nonholomorphic type of Poincarè compactification of the complex plane by a closed disk $\mathbf{D}$ corresponds to ``stereographic'' projection from the center of a sphere, rather than its North Pole, to the tangent plane $\mathbb{C}$ of the South Pole.
The closed disk $\mathbf{D}$ then represents the closed lower hemisphere, and its equatorial boundary circle $\mathbb{S}^1$ represents the compactifying ``circle at infinity‘‘ of the complex plane  $\mathbb{C}$.

Specifically, we consider a smooth real Euler multiplier which coincides with $|z|^{2(d-2)}$ near $z=0$ and is strictly positive elsewhere.
Near $z=0$, the original ODE \eqref{ODEw} reads 
\begin{equation}
\label{gpolyz}
|z|^{2(d-2)}\,\dot{z}=-\bar{z}^{d-2}(1+\ldots+f_0 z^d)\,.
\end{equation}
Rescaling time by the Euler multiplier replaces \eqref{gpolyz} by the real analytic equation
\begin{equation}
\label{gpolyzEuler}
\dot{z}=\tilde{f}(z):=-\bar{z}^{d-2}(1+\ldots+f_0 z^d)\,,
\end{equation}
for small $|z|$.
As anticipated, we have lost complex analyticity, but at least we have retained real analyticity near $z=0$, on the right hand side.

For Masuda's quadratic case $d=2$, the vector field becomes regular nonzero at $z=0$.
For $d\geq 3$, in contrast, we introduce an equilibrium at $z=0$.
For $d=3$, the equilibrium $z=0$ is a hyperbolic saddle. 
The unstable and the stable manifolds separate the plane into four hyperbolic sectors, locally.
For $d>3$, the equilibrium is degenerately hyperbolic, in the terminology of \cite{Hartman}.
In the following, see figure \ref{fig2} for an illustration of the special cyclotomic cases $f(z)=z^d-1$ with  $d=3,4$.

For general $d\geq3$, we introduce polar coordinates $z=\rho\exp(\mi\alpha)$\,:
\begin{align}
\label{rho}
\dot{\rho}\ &=\rho\big(-\cos((d-1)\alpha)+\ldots \big)\,;   \\
\label{alpha}
\dot{\alpha}\ &=\phantom{\rho\big(-}\ \sin((d-1)\alpha)+\ldots\,.
\end{align}
Here we have vested \eqref{gpolyzEuler} with yet another local Euler multiplier $\rho^{-d+3}$, and we have omitted higher order terms in $\rho$.
This ``blow up'' at $z=0$, in the sense of singularity theory, replaces $w=\infty$ by the circle $\rho=0,\ \alpha\in\mathbb{S}^1$.
In other words, local polar coordinates for $z=1/w$ essentially compactify the plane $w\in\mathbb{C}=\mathbb{R}^2$ to the closed unit disc $\mathbf{D}$, instead of the Riemann sphere $\widehat{\mathbb{C}}$.
Indeed the boundary circle $\mathbb{S}^1$ of $\mathbf{D}$ corresponds to the circle $\alpha\in\mathbb{S}^1$ at $\rho=0$.

This is equivalent to Poincaré compactification of $w\in\mathbb{C}$.
See figure \ref{fig2} (a),(c).
On the invariant boundary circle $\rho=0$, a total of $2(d-1)$ hyperbolic saddle equilibria $\alpha_k=\pi k/(d-1)\in\mathbb{S}^1$ appear, for $k\,\textrm{mod}\,2(d-1)$.
In figure \ref{fig2} (a) and (c), we label the saddles at $\alpha_k$ by $\mathbf{k}$\,.
Within the boundary circle, they are alternatingly unstable, at even $\mathbf{k}$, and stable, at odd $\mathbf{k}$.
Their corresponding stable and unstable manifold counterparts, received from and sent into $|w|=|1/z|<\infty$, are marked red and blue, respectively.
In \eqref{ODEw}, they mark red blow-up and blue blow-down of $w(t)$, in finite real original time $t=r$.
In \eqref{gpolyzEuler}, they also delimit the $2(d-1)$ local hyperbolic sectors of the equilibrium $z=0$, according to the planar classification of degenerate saddles by Poincaré. 
See section VII.9 in \cite{Hartman}, and (b), (d) of figure \ref{fig2}.
In (b), i.e. for $d=3$, red and blue simply mark the $2(d-1)=2+2$ half branches of the stable and unstable manifolds at the hyperbolic equilibrium $z=0$ of \eqref{gpolyzEuler}.

Let us compare Poincaré compactification \eqref{rho}, \eqref{alpha} with the branching at $t=0$ associated to the ramification point $(w,t)=(\infty,0)$ of the branched covering projection $\mathbf{q}_t$\,; see theorem \ref{Rthm}\emph{(iii)} and our comments in section \ref{RRes}.
The $d-1$ branches converging towards the boundary circle (red) simply correspond to the $w$-projection $\mathbf{q}_w$ of the $d-1$ lifts of the negative real axis $t=r<0$ to the sheet of the ramification point in the Riemann surface $\mathcal{R}$.
Indeed, the $d-1$ lifted trajectories $(r,w(r))\in\mathcal{R}$ given by $z(r)=1/w(r)$ in \eqref{sovz} parametrize the red blow-up trajectories in figure \ref{fig2}, for $r\nearrow 0$.
Similarly, the $(d-1)$ branches emanating from the boundary circle (blue) correspond to real $t=r>0$ with blow-down for $r\searrow 0$.
In normal form \eqref{branchm} of branching, i.e. for coordinates $t=\zeta^{d-1}$, the lifts of the real axis become $d-1$ straight lines $\Im \zeta^{d-1} = 0$ in $\zeta$, at angles $k\pi/(d-1),\ 0\leq k<d-1$.
After biholomorphic transformation of $\zeta$ to $z=1/w$ of \eqref{sovz}, the alternating blue and red parts become $d-1$ real analytic curves in the complex $z$-plane which intersect, at $z=0$, under equal asymptotic angles $\pi/(d-1)$; see \eqref{alpha} and figure \ref{fig2}, again.

\textbf{Glossary}. We summarize this section with a brief glossary of terms and color codings concerning the real-time dynamics of Poincaré compactifications, for perusal in this paper.
See figure \ref{fig2} (a),(c) for an illustration of all terms.
Dynamics resides on the closed unit disk $\mathbf{D}$ with \emph{interior} $\mathbb{D}$ and invariant \emph{boundary circle} $\mathbb{S}^1$.
The $d$ \emph{interior equilibria} $e_j\in\mathbb{D},\ 0\leq j<d$, possess nonzero complex linearizations $f'(e_j)=1/\eta_j$ and are called (blue) \emph{sinks}, (purple)\emph{Hopf} or \emph{Lyapunov centers}, and (red) \emph{sources} in case $\Re f'(e_j)$ is strictly negative, zero, or strictly positive, respectively.
Heteroclinic orbits $e_j\leadsto e_k$ from sources $e_j$ to sinks $e_k$ in $\mathbb{D}$ (black) foliate open regions and are called \emph{interior source/sink heteroclinics}.

The $2(d-1)$ \emph{boundary equilibria} $\mathbf{k}\in\mathbb{S}^1$ are equidistantly spaced, at angles $\alpha_k=\pi k/{(d-1)}\in\mathbb{S}^1,\ 0\leq k < 2(d-1)$.
They are all hyperbolic saddles of alternating stability and instability inside $\mathbb{S}^1$.
For even $k$, the unique (red) half branches of their stable manifolds (i.e., separatrices) emanating into $\mathbb{D}$ 
are called  \emph{blow-up orbits}.
Their unique (blue) unstable manifold separatrix counterparts, arriving from $\mathbb{D}$ are called \emph{blow-down orbits}.
If two such separatrices happen to coincide, they form a (purple) saddle-saddle heteroclinic orbit between boundary saddles, which we call an \emph{interior saddle-connection} or \emph{blow-down-up} orbit.
In all other cases, blue blow-down orbits are heteroclinic orbits towards sinks, and red blow-up orbits are heteroclinic orbits emanating from sources.
All non-equilibrium orbits within the invariant boundary circle $\mathbb{S}^1$ are saddle-saddle heteroclinic orbits  between adjacent boundary saddles, which we call \emph{boundary saddle-connections}.

\subsection{Cyclotomic examples}\label{Cyclotomic}

Geometric analysis of the special cyclotomic case $f(w)=w^d-1$ also illustrates how Masuda's paradigm is limited to the quadratic case $d=2$.
It will be an easy exercise to obtain the global planar phase portraits of figure \ref{fig2} under real-time flows $t=r\in\mathbb{R}$, say for $d=3,4$.

For general $d\geq 2$, let us consider the cyclotomic examples 
\begin{equation}
\label{cyclotomic}
\dot{w}=f(w)=w^d-1\,.
\end{equation} 
The equilibria are $w=e_j=\exp(2\pi\mi j/d)$ with $f'(e_j)=de_{-j}$ and $\eta_j=\tfrac{1}{d}\,e_j$\,, for $j\,\mathrm{mod}\,d$.
In particular, the nonvanishing assumption \eqref{sumseta} of theorem \ref{f'thm} is violated for nonprime degrees $d$ with a proper divisor $d'\neq 1,d$.
Indeed let $J=\{k\, d' \,\mathrm{mod}\,d\}$ denote the integer multiples of $d'$.
Then
\begin{equation}
\label{sumsetaviol}
d\cdot\sum_{j\in J}\,\eta_j\quad =\sum_{k \,\mathrm{mod}\,d/d'}\,e_{k d'} \ =\ 0\,.
\end{equation}
This example demonstrates how the nonvanishing part of assumption \eqref{sumseta} may not be necessary, after all.
Technically, the assumption will prevent clustering of equilibria in the proof of closedness lemma \ref{fbclosed} below.

From \eqref{rho}, \eqref{alpha}, more dynamically, we recall invariance of the boundary circle $\rho=0,\ \alpha\in\mathbb{S}^1$ at $w=\infty$, with $2(d-1)$ alternating hyperbolic saddles $\mathbf{k}$\,.
Since the cyclotomic polynomial $f$ possesses real coefficients, complex conjugation provides reflection symmetry between the flows in the upper and lower half plane of $w=u+\mi v$.
Indeed $w(t)$ solves \eqref{ODEw}, if and only if $\overline{w}(t)$ does.
In particular, the horizontal real axis is invariant.
This determines the horizontal blow-up / blow-down orbits of the boundary equilibria $\mathbf{k}=\mathbf{0},\ \mathbf{d-1}$\,, and the real blow-up orbit, or orbits, of $e_0$\,.
For even $d$, the real heteroclinic orbit $e_0\leadsto e_{d/2}$ extends to, both, real blow-up $e_0\leadsto\mathbf{0}$ and to real blow-down $\mathbf{d-1}\leadsto e_{d/2}$ via $w=\infty$, in complex time.
For odd $d$, blow-up $e_0\leadsto\mathbf{0}$ and $e_0\leadsto\mathbf{d-1}$ occurs in both real time directions.
By reflection symmetry, it remains to study the upper half plane $v>0$ of $w=u+\mi v$.

Consider the case $d=3$ of figure \ref{fig2} (a), (b), first.
The absence of Lyapunov centers excludes periodic orbits; see lemma \ref{perlem}\emph{(iv)}.
The Poincaré-Bendixson theorem \cite{Hartman} therefore identifies the sink $e_1$ as the only possible $\boldsymbol{\omega}$-limit  set of the blue blow-down separatrix of $\mathbf{1}$ in the upper half plane.
All other trajectories $w$ in the open upper half plane are heteroclinic of type $e_0\leadsto e_1$.

In case $d=4$ of figure \ref{fig2} (c), as for any even $d$, we encounter time reversibility of \eqref{cyclotomic}.
The time reversor is horizontal reflection at the vertical imaginary axis.
For even $d$, indeed, $w(t)$ solves \eqref{cyclotomic}, whenever $-\bar{w}(-t)$ does.
Therefore, the Lyapunov centers $e_1, e_3 =\pm\mi$ are locally surrounded by periodic orbits, only, up to the purple interior saddle-connections $\mathbf{1}\leadsto\mathbf{2}$ and $\mathbf{5}\leadsto\mathbf{4}$.
In fact, the unstable blow-down separatrix of $\mathbf{1}$ (usually blue, here purple) has to cross the imaginary axis $\Re w=0<\Im w$. 
Indeed, the first quadrant $\{\Re w>0,\ \Im w>0\}$ does not contain equilibria and, therefore, cannot contain periodic orbits.
See lemma \ref{perlem}\emph{(ii)}.
Reflecting on the first crossing, say at time $t=0$, reversibility implies that the unstable blow-down separatrix has to coincide with the stable blow-up separatrix of $\mathbf{2}$, usually colored red.
We have therefore colored the resulting interior saddle-connection $\mathbf{1}\leadsto\mathbf{2}$ purple.
Similar arguments show that the resulting heteroclinic cycle between $\mathbf{1}$ and $\mathbf{2}$ is filled with periodic orbits around the center $e_1$.
(In absence of reversibility, see also lemma \ref{perlem}\emph{(v)}.)
In (d), i.e. upon identification of the boundary circle $\rho=|z|=0,\ \alpha\in\mathbb{S}^1$ with $w=\infty$, the heteroclinic cycles become counter-rotating homoclinic to $z=1/w=0$.
Symmetrically, the homoclinic lobes remain foliated by iso-periodic families of synchronously counter-rotating periodic orbits, as in lemma \ref{perlem}\emph{(v)}.
Indeed their minimal periods $\mp\pi/2$, in original time and prior to any rescaling, all coincide.
The periodic families are unbounded in $w=1/z$.
All remaining trajectories $w$ in the open upper half plane are again interior source/sink heteroclinic, of type $e_0\leadsto e_2$.

Similar real-time homoclinic lobes to $w=\infty$, alias interior saddle connections, foliated by counter-rotating iso-periodic orbits arise around the Lyapunov centers $e_{d/4}$ and $e_{3d/4}$\,, for all polynomial degrees $d \equiv 0 \mod 4$. 
In imaginary time and for $d \equiv 2 \mod 4$, the analogous phenomenon occurs around $e_0=1$ and $e_{d/2}$\,.

Figure \ref{fig2} also demonstrates the behavior of the global trajectories $w=\Gamma_{\pm\mi\eps}$ under slight perturbations $w(0)=2\pm\mi\eps,\ \eps\searrow 0$, of the real blow-up initial condition $w(0)=2$.
In the case $d=3$ of figure \ref{fig2} (a), the interior source/sink heteroclinic $w=\Gamma_{+\mi\eps}: e_0\leadsto e_1$ in the upper half plane is labeled by $+\mi\eps$ (black).
For small $\eps>0$, it closely follows the concatenated heteroclinic chain $e_0\leadsto \mathbf{0}\leadsto \mathbf{1}\leadsto e_1$. 
The interior heteroclinic $\Gamma_{-\mi\eps}: e_0\leadsto e_2$ in the lower half plane\,, in contrast, labeled $-\mi\eps$, closely follows the \emph{different} concatenation $e_0\leadsto \mathbf{0}\leadsto \mathbf{3}\leadsto e_2$.
In sharp contrast to the quadratic case $d=2$, therefore, the red blow-up cannot be circumnavigated in complex time. 

The polar view of (b), centered at $w=\infty$ alias $z=1/w=0$, conflates the boundary circle $\rho=|z|=0,\ \alpha\in\mathbb{S}^1$ to a single point.
This shortens the black trajectories $\Gamma_{\pm\mi\eps}$ to perturbations of the two \emph{distinct concatenations} $e_0\leadsto \infty\leadsto e_1$ and $e_0\leadsto \infty\leadsto e_2$\,.
In the limit $\eps\searrow 0$, this amounts to a shared red blow-up  $e_0\leadsto \infty$ followed by two different subsequent blue blow-down orbits $\infty\leadsto e_1$ and $\infty\leadsto e_2$\,.
Again, this is in sharp contrast to the case $d=2$ of figure \ref{fig1}, where both approximations to the cyan blow-up orbit $1\leadsto\infty$ continue along \emph{the same} cyan blow-down $\infty\leadsto-1$.
Since this happens in finite original time $w(t)$, the Masuda paradigm of complex near-recovery for near-homogeneous PDE solutions, after real blow-up, turns out to be a peculiarity which is limited to the quadratic nonlinearity $d=2$.

The first nonprime degree $d=4$ of figure \ref{fig2}, (d) leads to heteroclinic perturbations $w=\Gamma_{\pm\mi\eps}: e_0\leadsto e_2$ which, at least formally, seem to limit onto identical concatenations $e_0\leadsto \infty\leadsto\infty\leadsto e_2$, for $\eps\searrow 0$.
Both limits share the initial red blow-up part $e_0\leadsto\infty$ and the terminal blue blow-down part $\infty\leadsto e_2$\,.
The two purple homoclinic blow-down-up orbits $\infty\leadsto\infty$, however, remain distinct.
The limiting homoclinic loop part $+\mi\eps$ of $\Gamma_{+\mi\eps}$ contains the Lyapunov center $e_1$ in its interior, whereas $-\mi\eps$ of $\Gamma_{-\mi\eps}$ surrounds the Lyapunov center $e_3$\,, along with their mandatory homoclinic lobes, respectively filled by nested families of counter-rotating, synchronously iso-periodic orbits with constant minimal periods $\mp\pi/2$.
Again, these perturbations run against the quadratic Masuda paradigm.
The observed discrepancies between $\Gamma_{\pm\mi\eps}$ run deeper than ``mere technical'' PDE difficulties like the absence of a heat semiflow in reverse time $r=\Re t<0$.
In fact they already originate from non-unique complex continuation of homogeneous real-time blow-up itself.

\subsection{Classifying Poincaré compactifications, and counting}\label{ClassRes}

In this section we classify the Poincaré compactifications \eqref{rho}, \eqref{alpha} of flows to \eqref{ODEw}.
We consider degrees $d\geq2$ and assume that all $d$ finite equilibria $e_j$ of \eqref{ODEw} are hyperbolic, i.e. $\Re\,f'(e_j)\neq0$.
Recall that this makes $e_j$ a \emph{sink}, for $\Re\,f'(e_j)<0$, and a \emph{source} for $\Re\,f'(e_j)>0$.

We compare and count the compactified flows, in real time $t=r$, up to \emph{orientation preserving} $C^0$ \emph{orbit equivalence}.
In other words, two flows on the closed unit disk $\mathbf{D}$ are considered orbit equivalent, if there exists an orientation preserving disk homeomorphism $H: \mathbf{D}\rightarrow\mathbf{D}$ which maps real-time orbits to real-time orbits.
We do not require $H$ to conjugate the real-time flows, and we do allow $H$ to reverse the real time direction.
So, $H$ conjugates real-time orbits, as sets, but not necessarily their specific time parametrizations.

Abstractly, we compare orbit equivalence classes to equivalence classes of planar undirected trees with $d$ vertices and $d-1$ edges, say contained in the open unit disk $\mathbb{D}$.
A \emph{tree} is a finite, connected, undirected graph without cycles.
\emph{Planar trees} are embedded in the plane $\mathbb{R}^2$.
We consider two planar trees as \emph{equivalent}, if there exists an orientation preserving disk homeomorphism $H$ of $\mathbb{D}$ which acts as a graph isomorphism on the trees.
In other words, $H$ maps vertices to vertices, and edges to edges, preserving their adjacency relations and the left cyclic orderings of corresponding edges around each vertex.
The direction of edges may be reversed under $H$.
We do not require any distinguished vertices or edges to be mapped to each other, e.g. any vertices marked as ``roots'' or otherwise labeled. 

Specifically, as we will explain further in section \ref{StarQuad}, the vertices correspond to sources and sinks in the Poincaré compactified flows \eqref{ODEw}.
Edges indicate the existence of families of heteroclinic orbits between vertices.
Planar embedding aside, this identifies the tree as the \emph{connection graph} among sources and sinks in $\mathbb{D}$.
For many other applications of this concept, see for example \cite{Conley, brfi88, brfi89, Mischaikow, firoSFB, firoFusco, Yorke-a, Yorke-b} and the many references there.

A third view point are \emph{diagrams of} $d-1$ \emph{non-intersecting chords} of the closed unit disk $\mathbf{D}$, up to rotation.
Here the chords are closed undirected straight lines between their $2(d-1)$ end points on the disk boundary $\mathbb{S}^1$, spaced equidistantly at angles $\beta_k:=\pi (k+\tfrac{1}{2})/(d-1),\ 0\leq k<2(d-1)$.
The non-intersecting chords are neither allowed to share end points, nor are they allowed to cross each other.
We call two \emph{chord diagrams} equivalent, if they coincide up to a proper rotation of $\mathbf{D}$.
Chord diagrams are also called non-crossing (single-armed) \emph{handshakes} of $2(d-1)$ people on a round table, or partitions of $2(d-1)$ elements into non-crossing blocks of size 2.

\begin{thm}\label{Classthm}
Consider flows \eqref{ODEw} of univariate polynomials $f$ with $d$ equilibria $e_j$\,.
For the reciprocal linearizations $\eta_j=1/f'(e_j)$ at $e_j$ we assume the following \emph{strong nondegeneracy condition}
\begin{equation}
\label{sumsetaRe}
 \sum_{j\in J}\,\Re\,\eta_j\,\neq\,0
\end{equation}
to hold, for any nonempty subset $\emptyset\neq J \subsetneq\{0,\ldots,d-1\}$.
Then the following two statements hold true, in terms of the orientation preserving equivalence classes just described.
\begin{enumerate}[(i)]
\item  The real-time phase portraits of the Poincaré compactifications correspond, one-to-one, to certain unlabeled, unrooted, undirected, planar trees with $d$ vertices.
\item Equivalently, they correspond, one-to-one, to certain chord diagrams of $d-1$ unlabeled, non-intersecting chords of the unit circle.
\end{enumerate}
\end{thm}

Correspondence theorem \ref{Classthm} only asserts that each phase portrait corresponds to \emph{some} planar tree, alias circular handshake or non-intersecting chord diagram.
It does not assert that \emph{all} planar trees actually do occur.
In other words: the correspondence is injective, but surjectivity is still missing.
The following realization theorem closes this gap.

\begin{thm}\label{Realizethm}
Each unlabeled, unrooted, undirected, planar tree with $d$ vertices, alias each circular handshake or each chord diagram of $d-1$ non-intersecting chords, is realized by Poincaré compactifications of \eqref{ODEw}, for suitable polynomials $f$ of degree $d$.
\end{thm}

We postpone comments on realization theorem \ref{Realizethm} to section  \ref{RealizethmComm}.

For the construction of the precise correspondences, we refer to our  discussion of \emph{tree portraits}, in section \ref{StarQuad}, and to the proof of theorem \ref{Classthm}\emph{(i)}, in section \ref{OrbEqu}.
For $J=\{j\}$, strong nondegeneracy \eqref{sumsetaRe} excludes Lyapunov centers $e_j$\,, i.e.~equilibria with purely imaginary eigenvalues. 
Collaterally, this also excludes nonstationary periodic orbits in the open disk $\mathbb{D}=\mathbf{D}\setminus \mathbb{S}^1$.
See lemma \ref{perlem} for details.
General $J$ in assumption \eqref{sumsetaRe} forbid any interior saddle-connections.
All $2(d-1)$ saddles $\mathbf{0},\ldots,\mathbf{2d-3}$ sit on the boundary circle  $\mathbb{S}^1$.
Boundary saddle-connections within the flow-invariant boundary circle $\mathbb{S}^1$ itself, of course, cannot be excluded.
See sections \ref{PoinRes}, \ref{Cyclotomic}, and lemma \ref{sad2lem1}.

In absence of Lyapunov centers and interior saddle-connections, the source/sink connection graphs persist, under small perturbations of the polynomials $f$.
In particular, this recovers standard results on structural stability of Morse systems \cite{PalisSmale, Sotomayor, PalisdeMelo}: small perturbations of equilibria only produce $C^0$ orbit equivalent Poincaré compactifications.

Correspondence theorem \ref{Classthm} extends well beyond local perturbations.
See figure 59 in \cite{PalisdeMelo} for an example of planar phase portraits with planar trees which are isomorphic, in the sense of abstract graph theory.
The planar trees are not equivalent, however, by a planar homeomorphism.
The planar phase portraits, consequently, also fail to be $C^0$ orbit equivalent.

Explicit counts of planar trees have been provided by \cite{oeis} and in theorem 2 of \cite{Treecount1}, as follows.
See also $\omega_D^2$ in \cite{Treecount2}, equation (10) and table 1. 
By correspondence theorem \ref{Classthm} and realization theorem \ref{Realizethm}, we therefore obtain the following counts of real-time global phase portraits for polynomial ODEs \eqref{ODEw} with $d$ nondegenerate source/sink equilibria.

\begin{thm}\label{Classcor} \cite{Treecount1, oeis, Treecount2}
Up to orientation preserving equivalence, the number of unlabeled, unrooted, undirected, planar trees with $d$ vertices, alias chord diagrams of $d-1$ unlabeled, non-intersecting chords of the unit circle, is given by entry $A_{d-1}$ of sequence A002995 in the online encyclopedia of integer sequences \cite{oeis}.
The counts for $2\leq d\leq 16$ are
\begin{equation}
\label{Treecount}
1, 1, 2, 3, 6, 14, 34, 95, 280, 854, 2694, 8714, 28640, 95640, 323396.
\end{equation}
An explicit expression for the general counts $A_{d-1}$ in closed form is
\begin{equation}
\label{Treecount1}
\begin{aligned}
A_{d-1}\,=\,\tfrac{1}{2(d-1)d}\, \tbinom{2(d-1)}{d-1}\,&+\,\tfrac{1}{4(d-1)}\tbinom{d}{d/2}\,+\,\tfrac{1}{d-1}\,\phi(d-1) \,+\\
                  &+\,\tfrac{1}{2(d-1)}\,\sum_{k=2}^{d-2}\,\tbinom{2k}{k}\, \phi(\tfrac{d-1}{k}) \,.           
\end{aligned}
\end{equation}
Here $\phi$ denotes the Euler totient count of coprime elements, and the sum only runs over proper divisors $k$ of $d-1$.
For odd $d$, the second summand on the right is omitted.
\end{thm}

Since complex conjugation $e_j \mapsto \overline{e}_j$ of all equilibria reflects their phase portrait, orientation reversing homeomorphisms should be admitted, as well.
Alas, the corresponding counts of $A_{d-1}$, reduced by reflection symmetry, seem not to have been studied.

The single quadratic case $A_1=1$ of $d=2$ equilibria has been discussed in section \ref{ODE2}: indeed all cases are holomorphically equivalent, by a Möbius automorphism of $\widehat{\mathbb{C}}$.
The single cubic $A_2=1$ of $d=3$ is equivalent to the spatially homogeneous complex version of the famous Chafee-Infante attractor $f(w)=\lambda w(1-w^2)$; see \cite{chin74,brfi88,firoSFB, firoFusco} for the real case. 
Blow-up in real time occurs, e.g., for purely imaginary initial conditions.
Both cases $A_3=2$ with $d=4$ arise by perturbations of figure \ref{fig2} (c), (d) to hyperbolic cases with one or two sinks $e_m$, respectively.
See also our discussion in section \ref{ClassPf}.
A first pair of inequivalent hyperbolic examples with the same number of (two) sinks occurs among the $A_4=3$ cases with $d=5$ equilibria.
In that case, homotopies in the reciprocal linearizations $\eta_j=1/f'(e_j)$ connect all cases via Lyapunov centers at purely imaginary eigenvalues $f'(e_j)$, in leaves of the trees, and their collateral interior saddle-connections in the Poincaré compactification.
See section \ref{BifGraph} and theorem \ref{Bdthm} for further comments on the connectivity of the resulting bifurcation graph $\mathcal{B}_d$\,.

\subsection{Outline} \label{Out}

We proceed as follows.
In section \ref{f'Pf} we prove theorem \ref{f'thm} on prescribed linearizations $f'(e_j)$ at zeros $e_j$ of polynomials $f$.
Theorem \ref{Rthm} on the global Riemann surface of ODE solutions is proved in section \ref{RPf}.
Section \ref{ClassPf} establishes the equivalences of Poincaré compactifications collected in correspondence theorem \ref{Classthm}.
Section \ref{RealizethmComm} sketches a proof of realization theorem \ref{Realizethm}; for full details we refer to \cite{FiedlerYamaguti}.
We discuss some higher-dimensional aspects, in section \ref{Dis}.
Specifically we address blow-up in complex time, for homoclinic and heteroclinic orbits of entire nonlinearities $f$, rather than just polynomial ones; see section \ref{ODEhet}.
Section \ref{1000} recalls a 1000 \euro\ question concerning the possibility of complex entire homoclinic orbits for entire ODEs in higher-dimensional systems $w\in\mathbb{C}^N$.
The question was first raised in \cite{FiedlerClaudia}.
We conclude with a time-reversible real ODE of second order which features entire periodic orbits and, upon forcing, exponentially sharp Arnold tongues; see section \ref{EntirePer}.

\subsection{Acknowledgment}\label{Ack} 
This paper is dedicated to the memory of Leonid Pavlovich Shilnikov, distinguished leader of an internationally outstanding and warmly welcoming group of colleagues and friends at Nizhny, for decades.
We are also indebted to Vassili Gelfreich, for very patient explanations of his profound work on exponential splitting asymptotics, to Anatoly Neishtadt for lucid conversations on adiabatic elimination, and to Karsten Matthies, Carlos Rocha, Jürgen Scheurle, Hannes Stuke, and the late colleagues Marek Fila and Claudia Wulff, for their lasting interest in real and complex times.

\section{Proof of linearization}\label{f'Pf}

In this section, we prove theorem \ref{f'thm} on the realizability of reciprocal derivatives $\eta_j=1/f'(e_j)$. 
We recall the constraint $\sum_j \eta_j = 0$ for sums over all $0\leq j<d$, and the nonvanishing assumption $\sum_{j\in J} \eta_j\neq 0$ for any  $\emptyset\neq J \subsetneq\{0,\ldots,d-1\}$.
See \eqref{sumseta}.

For our proof, we reformulate the theorem in a technically more convenient manner.
First we note that translations $e_j\mapsto e_j+a$ do not affect the prescribed nonzero derivatives $f'$ at the zeros of $f$.
Indeed
\begin{equation}
\label{f'ej}
f'(e_j)=\prod_{k\neq j} (e_j-e_k)\,,
\end{equation}
with $0\leq k,j < d$.
This allows us to fix $e_0=0$ and $0\notin J$, without loss of generality.
We then rewrite $f_j:=f'(e_j) $ as
\begin{equation}
\label{fj}
f_j =e_j\, \prod_{k\neq j} (e_j-e_k)\,.
\end{equation}
Here and below we restrict the range of indices to $0<j,k<d$.
We collect those $f_j$ as a function $\mathbf{f}$ of $\mathbf{e}=(e_1,\ldots,e_{d-1})$,
\begin{equation}
\label{fbe}
\mathbf{f}(\mathbf{e}):=(f_1,\ldots,f_{d-1})\,.
\end{equation}
Given $e_0=0$, we also introduce the spaces
\begin{align}
\label{E}
 \mathbf{E}\,&:=\, \{\mathbf{e}\in\mathbb{C}^{d-1}\,|\ e_j\neq 0 \textrm{ are pairwise disjoint}\, \}\,;   \\
\label{Y}
 \mathbf{Y}\,&:=\, \{\mathbf{y}=(y_1,\ldots,y_{d-1})\in\mathbb{C}^{d-1}\,|\textrm{ all }y_j\neq 0\, \}\,;  \\
\label{Y'}
 \mathbf{Y}'\,&:=\, \{\mathbf{y} \in \mathbf{Y}\,|\ \sum_{j\in J}\ 1/y_j\,\neq\, 0\,, \textrm{ for all }\, \emptyset\neq J \subseteq\{1,\ldots,d-1\} \} \,.
\end{align}
If we think of $y_j$ as prescribed values of $f_j=f'(e_j)$, then $\emptyset\neq J \subseteq\{1,\ldots,d-1\}$ stands for the index set $J$ in nonvanishing assumption \eqref{sumseta}.
In other words, consider the vector polynomial map
\begin{equation}
\label{fbEY}
\mathbf{f}:\mathbf{E}\rightarrow \mathbf{Y}\,.
\end{equation}
The value $0\neq f_0=f'(e_0)=1/\eta_0$ at $0\in J^c$, which we have ignored, is then correctly prescribed by \eqref{sumseta}, automatically.
Our very nonlinear interpolation theorem \ref{f'thm} therefore amounts to solving the equation $\mathbf{f}(\mathbf{e}^*)=\mathbf{y}^*$ for $\mathbf{e}^*\in\mathbf{E}$, given $\mathbf{y}^*\in \mathbf{Y}'$.
This leads to the following equivalent formulation of theorem \ref{f'thm}.

\begin{thm}\label{fbthm}
In the above notation and under the assumptions of theorem \ref{f'thm},
\begin{equation}
\label{fbY'}
\mathbf{f}(\mathbf{E})\supseteq \mathbf{Y}'\,.
\end{equation}
\end{thm}

Although the equation  $\mathbf{f}(\mathbf{e}^*)=\mathbf{y}^*$ is a polynomial system in $\mathbf{e}^*$, we did not succeed to bring algebraic methods to bear, based on ideals and Gröbner bases.
Our ``topological'' proof will rely on the following two technical lemmata, instead, which we prove afterwards.

\begin{lem}\label{fbopen}
The map $\mathbf{f}$ of \eqref{fbEY} is open.
\end{lem}

\begin{lem}\label{fbclosed}
The set $\mathbf{f}(\mathbf{E})\cap \mathbf{Y}'$ is closed in $\mathbf{Y}'$.
\end{lem}

\emph{Proof of theorem \ref{fbthm}\,.}\\
We first note that the spaces $\mathbf{E},\mathbf{Y},\mathbf{Y}'$ are open and connected.
Indeed, this is obvious for any finite-dimensional complex vector space, even upon removal of any finite number of subspaces of strictly positive complex codimension.
This works because the complex unit circle is connected, unlike $\{\pm1\}\in\mathbb{R}$\,; the corresponding claim fails in real vector spaces.
To show connectedness of $\mathbf{Y}'$, the same argument applies when we pass to the reciprocals $\eta_j=1/y_j$, biholomorphically.

It remains to prove that the set $\mathbf{f}(\mathbf{E})\cap \mathbf{Y}'$ is nonempty, open, and closed in $\mathbf{Y}'$. 
Then connectedness of $\mathbf{Y}'$ implies $\mathbf{f}(\mathbf{E})\supseteq \mathbf{Y}'$\,, as claimed in \eqref{fbY'}.

By lemmata \ref{fbopen} and \ref{fbclosed}, the set $\mathbf{f}(\mathbf{E})\cap \mathbf{Y}'$ is open and closed in $\mathbf{Y}'$.
The open domain $\mathbf{E}$ of $\mathbf{f}$ is nonempty, the removed sets $\mathbf{Y}\setminus \mathbf{Y}'$ are of strictly positive complex codimension, and $\mathbf{f}$ is open.
Therefore $\mathbf{f}(\mathbf{E})\cap \mathbf{Y}'$ is also nonempty.

This proves theorems \ref{fbthm} and \ref{f'thm}, up to the two lemmata on the map $\mathbf{f}$.
\hfill$\bowtie$

\emph{Proof of lemma \ref{fbopen}\,.}\\
To show that the map $\mathbf{f}$ of \eqref{fbe}, \eqref{fbEY} is open, we calculate the determinant $D_{d-1}$ of the Jacobian $\mathbf{f}'(\mathbf{e})$ as
\begin{equation}
\label{Dd1fact}
D_{d-1}\,=\,(d-1)!\,\prod_{k\neq j}\,(e_j-e_k)\,\neq\,0\,,
\end{equation}
for indices $0< j,k < d$ in the product.

We first establish the general form
\begin{equation}
\label{Dd1a}
D_{d-1}\,=\,a_{d-1}\prod_{k\neq j}\,(e_j-e_k)\,.
\end{equation}
We determine the nonzero integer coefficients $a_{d-1}$ later.
Explicitly, the entries of the Jacobian matrix $\mathbf{f}'(\mathbf{e})$ are the partial derivatives $f_{jk}:=\partial_k f_j$, i.e.
\begin{align}
\label{fjk}
f_{jk}\,&=\, -e_j\,\prod_{\ell\neq j,k}\ (e_j-e_\ell)\,=\,-f_j/(e_j-e_k)\,,\quad\mathrm{for}\quad k\neq j\,,\ \mathrm{and}\\
\label{fjj}
f_{jj}\,&=\,\prod_{k\neq j}\ (e_j-e_k)\,+\, \sum_{k\neq j} e_j\prod_{\ell\neq j,k}\ (e_j-e_\ell)\,=\,f_j/e_j+\sum_{k\neq j}\ f_{j}/(e_j-e_k)\,,
\end{align}
for $0< j,k,\ell < d$.
Each of the $(d-1)\times (d-1)$ entries is a homogeneous polynomial in $\mathbf{e}=(e_1,\ldots,e_{d-1})$ of degree $d-2$.
In particular, the polynomial $D_{d-1}$ is homogeneous of degree $(d-2)(d-1)$ with integer coefficients (or else identically zero).
If $e_k=e_{k'}$, the columns $k$ and $k'$ of the Jacobian matrix coincide.
Therefore $e_k-e_{k'}$ is a factor of the polynomial $D_{d-1}$.
Next we note equivariance \begin{equation}
\label{fequi}
\mathbf{f}\sigma=\sigma\mathbf{f}\,,
\end{equation} 
under any permutation $\sigma\in S_{d-1}$ of indices.
Therefore, the determinant $D_{d-1}$ is invariant under the action of $\sigma$.
Let $\sigma$ be the transposition $k\leftrightarrow k'$.
Then $\sigma$-invariance of $D_{d-1}$ implies, that not only $e_k-e_{k'}$, but even  $(e_k-e_{k'})^2$ must be a factor of $D_{d-1}$.
In other words, $D_{d-1}$ takes the form \eqref{Dd1a}, for some polynomial factor $a_{d-1}$.
However, the product on the right side of \eqref{Dd1a} is of the same polynomial degree $(d-2)(d-1)$ in $\mathbf{e}$ as the determinant $D_{d-1}$ itself.
Hence the coefficient $a_{d-1}$ must be a constant.

To show $a_d=d!$\,, as claimed in \eqref{Dd1fact}, we establish the recursion
\begin{equation}
\label{add1}
a_d=d\, a_{d-1}
\end{equation}
for the coefficients in \eqref{Dd1a}.
For $d=2$, we obtain $\mathbf{f}=f_1=e_1$ and $D_1=a_1=1$; see \eqref{fj}.
This agrees with the empty product in claim \eqref{Dd1fact}.

We slightly refine our notation to keep track of degrees of the polynomials $f(w)$ in the two cases $d-1$ and $d$ of our induction from $d-1$ to $d$.
We use superscripts $f_j^{d-1},\ \mathbf{f}^{d-1}$ for the maps defined in \eqref{fj}--\eqref{fbEY} which lead to the Jacobian determinant $D_{d-1}$.
We write $f_j^d\,,\,\, \mathbf{f}^d$ for their analogues leading to the Jacobian determinant $D_d$\,. 
For superscripts $d$ only, an extra zero $e_d$\ appears.

Analogously to \eqref{Dd1a}, for example, we obtain the recursion
\begin{equation}
\label{Dd1}
D_d\,=\,a_d\,\,\prod_{k\neq j}\,(e_j-e_k)\,=\,(-1)^{d-1}\frac{a_d}{a_{d-1}}\cdot D_{d-1}\cdot(e_d^{2(d-1)}+\ldots)\,,
\end{equation}
this time for indices $0<j,k\leq d$.
On the right, we have singled out all terms from the first product which involve  $j=d$ or $k=d$, i.e. all factors of $D_d$ which involve $e_d$\,.
We have also expanded for the highest power $2(d-1)$ of the new zero $e_d$\,; dots indicate terms of lower order.
The prefactor $(-1)^{d-1}$ accounts for the $d-1$ factors $k=d$ in the product.
Substitution of the remaining product over distinct $0<j,k<d$ by $D_{d-1}/a_{d-1}$ invokes the induction hypothesis \eqref{Dd1fact}, which also asserts $a_{d-1}=(d-1)!\neq 0$.

Similarly, \eqref{fjk} implies the explicit recursions
\begin{equation}
\label{fjkd}
\begin{aligned}
f_{jk}^d\ &=& (e_j-e_d)f_{jk}^{d-1}\,, \qquad\qquad\qquad&\ \mathrm{for}&\  0<j,k<d\,, \\
f_{jd}^d\ &=& -f_j^{d-1}\,,\qquad\qquad\qquad&\ \mathrm{for}&\  0<j<d\,,\ \\
f_{dk}^d\ &=& -e_d\,\prod_{\ell\neq d,k}\,(e_d-e_\ell)\,=\,-e_d^{d-1}+\ldots\,,\ &\ \mathrm{for}&\  0<k<d\,.
\end{aligned}
\end{equation}
Dots indicate terms of lower order than $d-1$ in $e_d$.
Last but not least, \eqref{fjj} implies
\begin{equation}
\label{fddd}
f_{dd}^d\ =\ \prod_{k\neq d}\ (e_d-e_k)\,+\, \sum_{k\neq d} e_d\prod_{\ell\neq d,k}\ (e_d-e_\ell)\ =\ (1+(d-1))e_d^{d-1}+\ldots \,.\quad\quad 
\end{equation}
Combining \eqref{fjkd} and \eqref{fddd} we obtain an alternative $e_d$-expansion of the determinant:
\begin{equation}
\label{Dd2}
\begin{aligned}
D_d\ &=\ 
\begin{vmatrix}
 (e_1-e_d)f_{11}^{d-1}     & \ldots & (e_1-e_d)f_{1,d-1}^{d-1} & -f_1^{d-1}  & \\
 \vdots     & \ddots & \vdots & \vdots &  \\
 (e_{d-1}-e_d)f_{d-1,1}^{d-1}     & \ldots & (e_{d-1}-e_d)f_{d-1,d-1}^{d-1} & -f_{d-1}^{d-1} &  \\[0.5em]
-e_d^{d-1}+\ldots    & \ldots & -e_d^{d-1}+\ldots &d\,e_d^{d-1}+\ldots &
\end{vmatrix} \ =\ \\[1em]
\ &=\ d\,D_{d-1}\ e_d^{d-1}\,\prod_{0<j<d}(e_j-e_d)\,+\,\ldots\ =\ (-1)^{d-1}\,d\,D_{d-1}\ e_d^{2(d-1)}\,+\,\ldots\,.
\end{aligned}
\end{equation}
In the second line, we have first used that $f_j^{d-1}$ and $f_{jk}^{d-1}$ do not depend on $e_d$, by definition.
We have also expanded the determinant with respect to the last row.
Because the last column does not contain $e_d$\,, the leading term of the expansion \eqref{Dd2} originates from the right lower corner of the last row. 

Comparison of the leading coefficients in determinant expansions \eqref{Dd1}, \eqref{Dd2}, and subsequent cancellation of $D_{d-1}\neq 0$, prove the remaining recursion claim \eqref{add1}.
This proves claim \eqref{Dd1fact} and the lemma.
\hfill$\bowtie$

\bigskip

By the explicit evaluation \eqref{Dd1fact} of the Jacobian determinant $D_{d-1}$\,, we have just proved that all values $\mathbf{y}\in \mathbf{Y}$ of the derivatives map $\mathbf{f}: \mathbf{E}\rightarrow \mathbf{Y}$ are regular values.
The set $\mathbf{f}^{-1}(\mathbf{Y}')$ denotes those equilibrium configurations $\mathbf{e}\in \mathbf{E}$ which satisfy the nonvanishing constraint \eqref{sumseta}.
On the complementary set $\mathbf{f}^{-1}(\mathbf{Y}\setminus \mathbf{Y}')$, we encounter some vanishing sum $\sum_{j\in J}\ 1/y_j\,=\,0$.
Since each vanishing is of complex codimension 1 in $\mathbf{Y}$, we immediately obtain the following corollary.

\begin{cor}\label{finvY'}
The set $\mathbf{f}^{-1}(\mathbf{Y}')$, where nonvanishing condition \eqref{sumseta} holds true, is open and dense in $\mathbf{E}$.\\
More precisely the complementary set $\mathbf{f}^{-1}(\mathbf{Y}\setminus \mathbf{Y}')$ of violations is a locally finite union of submanifolds embedded in $\mathbf{E}$, each of complex codimension one.
\end{cor}

\emph{Proof of lemma \ref{fbclosed}\,.}\\
To show that $\mathbf{f}(\mathbf{E})\cap \mathbf{Y}'$ is closed in $\mathbf{Y}'$, we consider any sequence $\mathbf{e}\in \mathbf{E}$ such that $\mathbf{Y}'\ni \mathbf{y}=\mathbf{f}(\mathbf{e})\rightarrow \mathbf{y}^*\in \mathbf{Y}'$.
We suppress the index of that sequence.
We have to establish the existence of $\mathbf{e}^*\in \mathbf{E}$ such that $\mathbf{f}(\mathbf{e}^*)=\mathbf{y}^*$.
In case $\mathbf{e}^*=\lim \mathbf{e} \in \mathbf{E}$ exists, this follows from mere continuity of $f$.

Even after passing to subsequences, two problems remain.
The sequence of $\mathbf{e}$ might not possess any limit $\mathbf{e}^*$, viz.~$\mathbf{e}$ is unbounded.
Or, the limit $\mathbf{e}^*\notin \mathbf{E}$ possesses some identical components which have coalesced in the limiting process.
Recall that we have normalized $\mathbf{e}$ to $e_0=0$, in expression \eqref{fj} for $f_j=y_j\rightarrow y_j^*$\,.
If any $e_j$ is unbounded, therefore, boundedness of $y_j$ implies $e_j-e_k\rightarrow 0$ for one or several factors in \eqref{fj}.
We call this effect \emph{clustering}.
In case the limit $\mathbf{e}^*\notin \mathbf{E}$ exists but possesses some identical components, of course, clustering occurs as well.

Permutation equivariance \eqref{fequi} of $\mathbf{f}$ allows us to relabel that $j$ to $j=0$, so that normalized $e_0=0$ becomes part of the cluster.
Let $j\in J$ denote the indices of the cluster, including $j=0$, for the duration of this proof.
We reserve the notation $k\notin J$ for the remaining $d-|J|$ indices.

To prove lemma \ref{fbclosed} we will exclude clustering, by contradiction.
Suppose clustering occurs.
On the cluster set $J$, we can then invoke nondegeneracy assumption 
\begin{equation}
\label{sumY'}
 \sum_{j\in J}\ 1/y_j^*\,\neq\, 0\,, \textrm{ for } 2\leq |J|<d-1,
\end{equation}
as specified in definition \eqref{Y'} for $\mathbf{y}^*\in \mathbf{Y}'$.
Note that we have adapted the condition to the inclusion of $0$ in the clustering set $J$ here.
To see that $J$ qualifies, first note that $J$ is nonempty by construction.
We have to exclude the case $|J|=d$ of $\mathbf{e}^*=0$.
But then $\mathbf{y}^*=0$, which is excluded by $\mathbf{y}^*\in \mathbf{Y}\not\owns 0$.

In violation of \eqref{sumY'}, we now claim 
\begin{equation}
\label{sumJ0}
 \sum_{j\in J}\,1/y_j^*\,=\,0
\end{equation}
This contradiction will prove the lemma.

To establish claim \eqref{sumJ0}, we first recall
\begin{equation}
\label{ej0}
e_j\rightarrow 0\,,\qquad\textrm{for all } j\in J,
\end{equation}
By definition of $J$ there also exists $\delta>0$ such that $|e_k|\geq\delta>0$ for all $k\notin J$.
Rescaling $e_\ell$ to $e_\ell/\delta$ for all $0\leq \ell < d$, without loss of generality, we may therefore assume 
\begin{equation}
\label{ek1}
|e_k|\geq 1, \textrm{ for all } k\notin J.
\end{equation}

To prove claim \eqref{sumJ0}, we now introduce the auxiliary polynomial
\begin{equation}
\label{gpolyw}
g(w)\,:=\,\prod_{j\in J} (w-e_j),
\end{equation}
quite analogously to $f$, but extended over the cluster indices $j\in J$, only.
Introducing $1/h_j=g_j=g'(e_j)$ on the cluster $J$, analogously to $\eta_j\,,f_j$ in \eqref{partialfrac}, \eqref{fj}, the residue theorem for $1/g(w)$ implies
\begin{equation}
\label{sumg'h}
\sum_{j\in J}\,\frac{1}{g'(e_j)}\,=\,\sum_{j\in J}\,h_j\,=\,0\,,
\end{equation}
analogously to \eqref{sumf'eta}.
We now claim
\begin{align}
\label{chj}
\lvert\, c&\,h_j-\eta_j\rvert \rightarrow 0\,, \qquad \mathrm{for}\\
\label{cdef}
c &:=(-1)^{d-|J|}\prod_{k\notin J} \,e_k^{-1}\,
\end{align}
and all $j\in J$.
By definition \eqref{cdef}, the (not necessarily bounded) sequence $c \in\mathbb{C}$ is independent of any $j\in J$.
Since $\eta_j=1/y_j \rightarrow 1/y_j^*$ has been assumed to converge, for all $j$, claim \eqref{sumJ0} follows from fact \eqref{sumg'h} and claim \eqref{chj}.

To prove the remaining claim \eqref{chj}, for any $j\in J$, we decompose and expand
\begin{align}
\label{gjq}
\frac{1}{\eta_j}\,=\, f_j\,=\,& g_j \cdot q \, =\,h_j^{-1} \cdot q\,,\qquad \mathrm{where}   \\
\label{qdef}
q\,:=\,& \prod_{k\notin J}\,(e_j-e_k)\,.  
\end{align}
Since all $e_j\rightarrow 0$ by \eqref{ej0}, and all $|e_k|\geq 1$ by \eqref{ek1}, definitions \eqref{cdef} and \eqref{qdef} imply $c\,q\rightarrow 1$.
Since $\eta_j\rightarrow\eta_j^*\neq0$ also converge, back substitution of \eqref{gjq} in \eqref{chj} evaluates
\begin{equation}
\label{cheta}
\lvert\, c\,h_j-\eta_j\rvert=\lvert\eta_j\rvert \cdot \lvert\, c\,h_j/\eta_j-1 \rvert =\lvert\eta_j\rvert \cdot \lvert\, c\,q-1 \rvert \rightarrow 0.
\end{equation}
Backtracking, this proves claim \eqref{chj}, contradiction \eqref{sumJ0} to \eqref{sumY'}, and the lemma. \hfill$\bowtie$.

With the proofs of lemmata \ref{fbopen} and \ref{fbclosed}, the proofs of theorems \ref{fbthm} and \ref{f'thm} are now also complete.

\section{Proof of Riemann coverings}\label{RPf}

In this section we prove parts \emph{(i)--(iv)} of theorem \ref{Rthm}.
We make free and extensive use of the notation, terminology, and concepts summarized in appendix section \ref{Apx}.
We repeat that this is classical material, and we do not claim any originality here.

\emph{Proof of theorem \ref{Rthm}(i).}\\
Expansion \eqref{sovz} presents an explicit holomorphic germ for the general solutions of \eqref{sovw}, starting with $t_0=0$ at $w_0=\infty\in\widehat{\mathbb{C}}_d$, alias at $z_0=1/w_0=0$.
Maximal analytic continuation along integration paths defines $t$ as the integral of the holomorphic differential form $\omega=dw/f(w)$ along curves $\gamma=\gamma_w$ from $w_0=\infty$ to $w$ in $\widehat{\mathbb{C}}_d$:
\begin{equation}
\label{tomg}
t=\int_{\gamma}\, \omega \,;
\end{equation}
see \eqref{tomx}.
We define the Riemann surface $\mathcal{R}$ by maximal continuation of holomorphic germs for $t$ along paths $\gamma$. 

By construction, any two holomorphic germs $t_1(w;w_0)$ and $t_2(w;w_0)$ for $t$ at the same $w_0\in\widehat{\mathbb{C}}_d$ can only differ by the integral \eqref{tomg} along a closed curve $\gamma$ through $w_0$.
By free homotopy in $\widehat{\mathbb{C}}_d$\,, this is equivalent to a closed curve $\gamma$ through $w=\infty$.
In other words, the difference
\begin{equation}
\label{t1t2}
T:=t_1(w;w_0)-t_2(w;w_0)\in\mathrm{range}\, \mathcal{P}
\end{equation}
does not depend on $w$, locally, and is a period of the period map $\mathcal{P}$.
See section \ref{PeriodT}, and in particular \eqref{T}, \eqref{alltw}, \eqref{Trange}, \eqref{Tjres}.

In particular, we may uniquely represent each local germ $t=t(w;w_0)$ by $w_0$ and its constant term $t_0=t(w_0;w_0)$.
This allows us to represent the Riemann surface $\mathcal{R}$ as a set of pairs $(w_0,t_0)$; see \eqref{Rdef}.
However, we have to keep in mind that the topology on each fiber $t_0+\mathrm{range}\,\mathcal{P}$ is discrete, due to the underlying process of analytic continuation in $\mathcal{R}$ between any two elements of the same fiber.

The construction also proves that the projection
\begin{equation}
\label{qw}
\mathbf{q}_w: \mathcal{R}\rightarrow\widehat{\mathbb{C}}_d\,.
\end{equation}
is an unlimited, unbranched, normal covering with deck group
\begin{equation}
\label{deckqw}
\mathrm{deck}(\mathbf{q}_w) \,=\, \mathrm{range}\, \mathcal{P}.
\end{equation}
The deck group acts by addition on the fibers 
\begin{equation}
\label{Rw}
\mathcal{R}_{w_0} := \mathbf{q}_w^{-1}(w_0) = \{t_0\in\mathbb{C}\,|\,(w_0,t_0)\in\mathcal{R}\} = t_0+\mathrm{range}\,\mathcal{P}\,,
\end{equation}
for some $t_0$\,.
This proves claims \emph{(i)} of theorem \ref{Rthm}. \hfill $\bowtie$

\emph{Proof of theorem \ref{Rthm}(ii).}\\
For $d\geq3$, the universal cover $\widetilde{X}$ of the Riemann surface $X:=\widehat{\mathbb{C}}_d$ is the upper half plane $\mathbb{H}$; see \eqref{HcoversCd}.
In other words, the punctured Riemann sphere $\widehat{\mathbb{C}}_d$ is hyperbolic.
For claim \eqref{deckHRCd}, see section \ref{Apx} and in particular \eqref{pi1Cd}.

By the general construction \eqref{XYX} of section \ref{Apx}, this implies the unbranched, unlimited, normal coverings
\begin{equation}
\label{HRCdPf}
\mathbb{H}\xrightarrow{\mathbf{p}}\mathcal{R}\xrightarrow{\mathbf{q}_w}\widehat{\mathbb{C}}_d\,,
\end{equation}
as claimed in \eqref{HRCd}.
In particular, the Riemann surface $\mathcal{R}$ is also hyperbolic.

Closed loops $\gamma$ through $(w_0,t_0)=(\infty,0)\in\mathcal{R}$ are precisely those loops in $\widehat{\mathbb{C}}_d$ on which the period map $\mathcal{P}(\gamma)=0$ vanishes.
In other words, $\pi_1(\mathcal{R})=\ker\mathcal{P}$.
The remaining claims \eqref{deckHRCd}--\eqref{deckRCd} now follow from assertion \eqref{XYXdeck} and the homomorphism theorem for the period map $\mathcal{P}: \pi_1(\widehat{\mathbb{C}}_d) \rightarrow \mathbb{C}$; see \eqref{T}.
This proves claims \emph{(ii)} of theorem \ref{Rthm}. \hfill $\bowtie$

\emph{Proof of theorem \ref{Rthm}(iii).}\\
To study the holomorphic projection $\mathbf{q}_t:\mathcal{R}\rightarrow\mathbb{C}$ of $(w,t)$ onto the integral $t$ of the differential form $\omega=dw/f(w)$, we first recall that $\mathcal{R}$ consists of sheets which are graphs $w\mapsto t=t_0+t(w;w_0)$, locally near any $w_0\in\widehat{\mathbb{C}}_d$\,.

Concerning branching of $\mathbf{q}_t$\,, consider finite $w\in\widehat{\mathbb{C}}_d\setminus\{\infty\}$ first.
Then $dt=dw/f(w)$ is nonzero.
By the implicit function theorem, this implies that the maps $w\mapsto t$ via any sheet are locally biholomorphic, differing only by periods in $\mathrm{range}\,\mathcal{P}$.
In particular $\mathbf{q}_t$ is unramified, at any finite $w\in\widehat{\mathbb{C}}_d\setminus\{\infty\}$.

At $w=\infty\in\widehat{\mathbb{C}}_d$\,, alias at $z=1/w=0$, we instead recall the local expansion 
\begin{equation}
\label{sovtz}
t=t_0+t(z;0)= -\tfrac{1}{d-1}z^{d-1}(1+\ldots)\,,
\end{equation} 
uniformly for any period $t_0=T\in\mathrm{range}\,\mathcal{P}$.
See \eqref{sovz}.
Therefore the projection $\mathbf{q}_t$ is branched of multiplicity $d-1$, at the branching points $t\in\mathrm{range}\,\mathcal{P}$ with ramification points $(w,t)\in\mathcal{R}$ at $w=\infty$.

Next, we study the closure $\mathbb{B}$ of the set of branch points of the projection $\mathbf{q}_t$\,, i.e.
\begin{equation}
\label{B}
\mathbb{B}:= \mathrm{clos}\,\mathrm{range}\,\mathcal{P} = 2\pi\mi\,\mathrm{clos} \langle \eta_1,\ldots,\eta_{d-1}\rangle_\mathbb{Z}\,;
\end{equation}
see \eqref{Trange}, \eqref{Tjres}.
The cases in claim \eqref{closT} enumerate the closed Abelian subgroups of $(\mathbb{R}^2,+)$.
Since $\mathbb{B}$ is such a closed Abelian subgroup, it has to appear in that list.

Note that only the discrete first two cases $\mathbb{Z}$ and $\mathbb{Z}\times \mathbb{Z}$ for $\mathbb{B}$ properly accommodate the discrete topology of the fibers of the projection $\mathbf{q}_w$ of $\mathcal{R}$ and therefore provide an embedding $\mathcal{R}\hookrightarrow\widehat{\mathbb{C}}_d\times\mathbb{C}$, by representation \eqref{Rdef}.
In particular, the projection $\mathbf{q}_t:\mathcal{R}\rightarrow\mathbb{C}$ becomes an unlimited branched covering.

The remaining three cases feature dense immersions $\mathcal{R}\hookrightarrow\widehat{\mathbb{C}}_d\times\mathbb{C}$ of sheets.
Indeed, each fiber carries the discrete topology in the Riemann surface $\mathcal{R}$, but becomes dense in at least one real direction, in the topology of $t\in\mathbb{C}$.
The accumulation points of the branching closure $\mathbb{B}$ also signify that the projection $\mathbf{q}_t:\mathcal{R}\rightarrow\mathbb{C}$ does not qualify as a branched covering.

It remains to show that all five cases actually do arise.
By theorem \ref{f'thm}, we can prescribe $\eta_j$ arbitrarily, as long as nonvanishing condition \eqref{sumseta} remains satisfied.
We proceed case by case.

We have already encountered the first discrete case $\mathbb{B}=\mathrm{range}\,\mathcal{P} = \pi\mi\,\mathbb{Z}$ in the quadratic examples $d=2$ of a single generator; see sections \ref{ODE2}, \ref{RRes},  \ref{Cyclotomic}, and figure \ref{fig1}.

The second discrete case $\mathbb{B}\cong\mathbb{Z}\times \mathbb{Z}$ arises for $d=3$, for example, whenever we prescribe two generators $\eta_1, \eta_2$ which are linearly independent over $\mathbb{R}$.

Next suppose $d=3$, but $\eta_1, \eta_2$ are $\mathbb{R}$-collinear.
Then linear independence over $\mathbb{Z}$ or, equivalently, over $\mathbb{Q}$ implies $\mathbb{B}\cong\mathbb{R}$.
This is the third case.
Linear $\mathbb{Z}$-dependence, of course, reverts to $\mathbb{B}\cong\mathbb{Z}$ .

To construct case 4, consider $d=4$ with $\mathbb{R}$-dependent, $\mathbb{Z}$-independent  $\eta_1,\eta_2$ and any $\eta_3$ which is $\mathbb{R}$-independent from $\eta_1,\eta_2$\,. 
Indeed, this leads to case 4 of $\mathbb{B}\cong\mathbb{R}\times \mathbb{Z}$.


To construct an example for the final case 5 of $\mathbb{B} = \mathbb{C}= \mathbb{R}^2$\,, again with $d=4$, we first choose $\mathbb{R}$-independent generators like $\eta_1=1,\ \eta_2=\mi$ of the lattice .
We choose $\eta_3=a+\mi b$ such that $1,a,b$ are $\mathbb{Z}$-independent.
In particular, $\eta_1,\eta_2,\eta_3$ are $\mathbb{Z}$-independent.
Let $\mathbb{B}_3$ denote the closure of the subgroup generated by $\eta_3\, \mathbb{Z}$ in the 2-torus $\mathbb{T}^2= \mathbb{R}^2\,/\mathbb{Z}^2$\,.
Note that $\mathbb{B}_3$ cannot be discrete, i.e. finite, because $\eta_1,\eta_2,\eta_3$ are $\mathbb{Z}$-independent.
To show $\mathbb{B}_3=\mathbb{T}^2$, we only have to exclude that $\mathbb{B}_3$ is a line with rational slope.
Passing to covering space, $\eta_3\in\mathbb{B}_3$ of rational slope implies $m'(a+m)+n'(b+n)=0$, for some integers $m,m',n,n'$ with $m',n'\neq 0$.
Hence rational slope contradicts $\mathbb{Z}$-independence of $1,a,b$.
This contradiction proves $\mathbb{B}_3=\mathbb{T}^2$, and hence $\mathbb{B}= \mathbb{R}^2$.

Of course we have only presented constructions involving minimal degrees $d$, for each case.
All instances trivially extend to any higher degree as well.
This proves claims \emph{(iii)}. 

\emph{Proof of theorem \ref{Rthm}(iv).}\\
We have to prove equivariance for the liftings, by $\mathbf{q}_w$ and $\mathbf{p}\circ\mathbf{q}_w$\,,  of the local flow $\Phi^t$ from $\widehat{\mathbb{C}}_d\setminus\{\infty\}$ to local flows on $\mathcal{R}':=\mathcal{R}\setminus(\{\infty\}\times\ker\mathcal{P})$ and $\mathbb{H}':=\mathbb{H}\setminus\mathbf{p}^{-1}(\mathcal{R}')$, respectively.
On the base space $\widehat{\mathbb{C}}_d\setminus\{\infty\}$, the actions of the deck groups $\mathrm{deck}(\mathbf{q}_w)$ and $\mathrm{deck}(\mathbf{p}\circ\mathbf{q}_w)$ on $\mathcal{R}'$ and $\mathbb{H}'$ amount to compositions of local time shifts of the local flow $\Phi^t$ along closed curves $\gamma$.
See \emph{(ii)} and \emph{(iii)}.
These time shifts add up to periods $T=\mathcal{P}(\gamma)$.
Since the local flow $\Phi^t$ in the base space commutes with such time shifts, the lifted flows commute with the deck transformations on $\mathcal{R}'$ and $\mathbb{H}'$.

This proves deck equivariance of the lifted local flows and completes the proof of theorem \ref{Rthm}.
\hfill $\bowtie$

\section{Equivalences and case counts}\label{ClassPf}

After some preparations on real-time periodic orbits, interior saddle-connections and geometric properties of real-time Poincaré compactifications, in sections \ref{ODEper}--\ref{StarQuad}, 
we prove correspondence theorem \ref{Classthm} in section \ref{OrbEqu}.
The equivalence of formulations \emph{(i)} and \emph{(ii)} follows from purely graph-theoretic considerations; see \cite{oeis}.
As a relation among the various tree portraits associated to Poincaré compactifications, we study generic real one-parameter bifurcations in the guise of bifurcation graphs; see section \ref{BifGraph}.
Section \ref{RealizethmComm} sketches a proof of realization theorem \ref{Realizethm}.
For terminology we urgently refer to the glossary at the end of section \ref{PoinRes}.

\subsection{Periodic orbits}\label{ODEper}

We start with a brief folklore of some relations among Brouwer degree,  residues, Poincaré linearization, and synchronous real-time periodicity for scalar holomorphic ODEs; see e.g. \cite{Zoladek}.

\begin{lem}\label{perlem}
Let $e$ denote equilibria $f(e)=0$ of the scalar complex entire ODE \eqref{ODEw}, and let $w(r)\in\gamma$ denote nonstationary periodic orbits $\gamma$ of minimal period $|p|$, in real time $r$.
We write $p<0$ to indicate negative, clockwise orientation of $w$.
Positive, anti-clockwise orientation has $p>0$.
Then the following holds true.
\begin{enumerate}[(i)]
  \item The local Brouwer degree $\deg(f,e,0)$ of any equilibrium $e$ coincides with the algebraic multiplicity of $e$ as a zero of $f$.
  \item Any periodic orbit $\gamma$ surrounds a single equilibrium $e$.
  \item The surrounded equilibrium $e$ is algebraically simple.
  \item The signed minimal period $p$ of $w$ is given by the residue $2\pi\mi/f'(e)$ of $1/f$ at $e$, with purely imaginary linearization $f'(e)$ at the \emph{Lyapunov center} $e$. 
   \item The interior of $\gamma$ is foliated by synchronously iso-periodic orbits of constant signed minimal period $p$, which are nested around the Lyapunov center $e$.
\end{enumerate}
\end{lem}

\begin{proof}
To prove claim \emph{(i)}, let $e=0$ be a zero of $f(w)=w^m(a+\ldots)$ with algebraic multiplicity $m>0$ and some complex coefficient $a\neq 0$. 
Introduce polar coordinates $w=\rho\exp(\mi\alpha)$ locally, with small $\rho>0$. 
Then the local Brouwer degree at $f(0)=0$ is given by the winding number $m>0$ around $w=0$ of $\alpha\mapsto f(w)=\rho^m(a\exp(m\mi\alpha)+\ldots)$, for $0\leq\alpha\leq 2\pi$.
Alternatively, the winding number coincides with the residue $m$ of $(\log f)'=f'/f$ at $w=0$.

To prove claims \emph{(ii)} and \emph{(iii)}, let $\Omega$ denote the interior of $\gamma$. 
The winding number of $f$ along any nonstationary periodic orbit $\gamma=\partial\Omega$ equals 1. 
Here we have to view the set $\gamma$ as a left oriented cycle, even if the trajectory runs clockwise.
Replacing $f$ by $-f$, in fact, does not change that winding number.
By (non-analytic) diffeotopy of $\partial\Omega$ to the unit circle with left-oriented unit tangent,
the winding number equals 1, and also coincides with the Brouwer degree $\deg(f,\Omega,0)$. 
That Brouwer degree adds up the positive multiplicity contributions of all interior equilibria $e\in\Omega$.
Invoking \emph{(i)} proves \emph{(ii)}, i.e. existence, uniqueness, and simplicity of $e$.

To prove \emph{(iv)}, let $e$ denote any simple zero of $f$.
Then $f(w)=(w-e)(f'(e)+\ldots)$ identifies $1/f'(e)$ as the residue of the reciprocal $1/f(w)$ at $e$.
Separation of variables, on the other hand, identifies the signed real period $p$ by integration along the time-oriented periodic orbit $w$ as
\begin{equation}
\label{wperiod}
p = \oint_w \ \frac{1}{f(w)}\,dw = 2\pi\mi /f'(e)\,.
\end{equation}
Here we have used \emph{(ii)} and the residue theorem, to evaluate the integral.
Since $p$ is real, $f'(e)\neq 0$ is purely imaginary.
This identifies the equilibrium $e$ as a Lyapunov center.
The positive or negative time orientation of the periodic orbit $w$, depending on the sign of $\Im f'(e)\neq 0$, provides the appropriate sign of $p$.
This proves claim \emph{(iv)}.

We prove claim \emph{(v)} for left-winding $p>0$, without loss of generality.
Real analyticity of $r\mapsto w(r)$ extends to a narrow strip $|s|<\eps$ of complex times $t=r+\mi s$. 
Since the real and imaginary ODE flows commute in the strip, as in \eqref{flow}, all orbits $r\mapsto w(r+\mi s)$ with $|s|<\eps$ share the same minimal period $p>0$.
Because $\dot{w}(r)\neq 0$, the map $t\mapsto w(t)$ is locally biholomorphic and preserves orientation.
In particular, the orbits $r\mapsto w(r+\mi s)$ with fixed imaginary parts $0<s<\eps$ are nested around $e$, in the interior $\Omega$ of the orbit $\gamma$.
(Local foliation of the exterior is provided by constant $-\eps<s<0$.)
This shows that the region foliated by periodic orbits is open in $\Omega$.
Since $\Omega$ is also bounded by $\partial\Omega=\gamma$, limits of periodic orbits in $\Omega$ are again periodic, or coincide with $e$.
Connectivity of the disk $\Omega$, and \emph{(iv)}, imply the synchronous iso-periodic foliation of $\Omega\setminus\{e\}$.
This proves the final claim \emph{(v)}, and the lemma.
\end{proof}

Of course it is interesting to explore the boundary of the connected component defined by complex continuation of all periodic orbits surrounding a Lyapunov center $e$, locally, in the Poincaré compactification of our polynomial setting.
Recall figure \ref{fig2} (c) for an example with two Lyapunov centers $e_1$ and $e_3$.
Following \cite{Sotomayor}, we obtain a finite cycle of saddle-saddle heteroclinic orbits.
Their location alternates between the disk interior $\mathbb{D}$ and the invariant boundary circle $\mathbb{S}^1$.
Indeed, all saddles are hyperbolic, located on $\mathbb{S}^1$, and of alternating stability within $\mathbb{S}^1$.
This excludes homoclinic orbits, in the Poincaré compactification.
If the heteroclinic cycle consists of just a single pair of interior/boundary saddle-connections, however, we obtain a homoclinic orbit in $\widehat{\mathbb{C}_d}$, where the boundary circle $\mathbb{S}^1$ is collapsed to $w=\infty$.
See \ref{fig2} (d).
We will see in the next section how all other cases will require linear degeneracies of real codimension at least 2 in the space of reciprocal linearizations $\eta_j$\,.

\subsection{Interior saddle-connections}\label{SadSad}

In lemma \ref{perlem} we have seen how the strong nondegeneracy condition \eqref{sumsetaRe} for $J=\{j\}$ excludes Lyapunov centers $e_j$\,.
We will now see how the same condition excludes interior saddle-connections.
See figure \ref{fig2}(c) for an illustration of the following lemma.

\begin{lem}\label{sad2lem1}
Let $\gamma:\mathbf{A}\leadsto\mathbf{B}$ denote an interior saddle-connection in the open disk $\mathbb{D}$ between boundary saddle equilibria $\mathbf{A}, \mathbf{B} \in \mathbb{S}^1$.
Let $J$ enumerate the interior source/sink equilibria $e_j$ which are on the opposite side of $e_0$ in $\mathbb{D}\setminus\gamma$. Then $J$ is nonempty and
\begin{equation}
\label{sumsetaRe0}
 \sum_{j\in J}\,\Re\,\eta_j\,=\,0
\end{equation}
\end{lem}

\begin{proof}
The Jordan arc $\gamma$ divides the open disk $\mathbb{D}$ into two connected components.
We label the component of $\mathbb{D}\setminus\gamma$ which does not contain $e_0$ as ``interior''.
Let $J$ collect the indices of equilibria in the interior.
Integration along $\gamma$ determines the finite real time $T$ from $\mathbf{A}$ to $\mathbf{B}$ as
\begin{equation}
\label{TAB}
0<T=\int_\gamma \omega=\pm 2\pi\mi\,\sum_{j\in J} \eta_j\,;
\end{equation}
see \eqref{tw}.
In particular $J\subseteq\{1,\ldots,d-1\}$ is nonempty.
Here we have used the residue theorem for $w\in\widehat{\mathbb{C}}_d$ to evaluate the integral.
Indeed, the loop $\gamma$ becomes closed in $\widehat{\mathbb{C}}_d$\,, because $\mathbf{A}=\mathbf{B}=\infty$ there; see figure \ref{fig2} (d).
The indeterminacy of sign results from the winding orientation of the loop $\gamma$ with respect to its designated interior.
Taking imaginary parts of \eqref{TAB} proves claim \eqref{sumsetaRe0}, and the lemma.
\end{proof}

For a partial converse to lemma \ref{sad2lem1}, which derives existence and uniqueness of a saddle-connection $\gamma$ from degeneracy \eqref{sumsetaRe0} under an additional condition, see lemma \ref{sad2lem2} below.

Lemma \ref{perlem} provides an example where an equilibrium $e=e_j$ is a Lyapunov center.
This is equivalent to a degeneracy \eqref{sumsetaRe0} where $J=\{j\}$ is a singleton. 
By \cite{Sotomayor}, the boundary of the family of iso-periodic orbits emanating from the Lyapunov center $e_j$ consists of a succession of saddle-saddle heteroclinic orbits $\Gamma$.
In our case of Poincaré compactification, this identifies a unique pair of adjacent saddles on the $\mathbb{S}^1$ boundary, and a unique interior saddle-connection  $\Gamma:\mathbf{A}\leadsto \mathbf{B}$ between them.

\subsection{Stars, quadrangles, tree portraits, and handshakes}\label{StarQuad}

We now replace nondegeneracy \eqref{sumseta} by the strong nondegeneracy assumption \eqref{sumsetaRe}.
By sections \ref{ODEper} and \ref{SadSad}, this excludes Lyapunov centers, periodic orbits, and interior saddle-connections in the Poincaré compactification \eqref{rho}, \eqref{alpha}.
We recall that equilibria $\mathbf{k}$ on the boundary circle $\alpha\in\mathbb{S}^1$correspond to saddles at $\alpha=\pi k/(d-1),\ 0\leq k < 2(d-1)$.
Alternatingly along the circle, each boundary saddle receives a unique red blow-up orbit from the open disk $\mathbb{D}$, for even $\mathbf{k}$, and sends a unique blue blow-down orbit back into it, for odd $\mathbf{k}$. 
See figure \ref{fig2} (a),(c).

By absence of Lyapunov centers, periodic orbits, and interior saddle-connections, any blue blow-down orbit emanating from a boundary saddle terminates at some interior blue sink.
This follows from the Poincaré-Bendixson theorem \cite{Hartman}.
Similarly, in reverse time, any red blow-up orbit towards a boundary saddle emanates from an interior red source.
Any other interior orbit in $\mathbb{D}$ has to be heteroclinic, from some red source to some blue sink.

To describe the resulting phase portraits, globally, we introduce stars and quadrangles.
Let $e\in\mathbb{D}$ denote any source or sink equilibrium.
Then the \emph{star} of $e$ is the set of all red blow-up, or blue blow-down, orbits between $e$ and any (boundary) saddles.
Recall that each boundary contributes a unique such red blow-up or blue blow-down interior orbit.
Each interior red source and blue sink, likewise, has to meet at least one such interior red blow-up or blue blow-down orbit.
Indeed, sources or sinks without any blow-up or blow-down orbit cannot exist, for degrees $d>1$: their domain of attraction in forward or backward time would consist of the whole open disk $\mathbb{D}$.
After local $C^0$ Grobman-Hartman linearization \cite{Hartman}, the stars consist of straight half-lines, locally.
A \emph{sector} of $e$ consists of all orbits which emanate or terminate between two angularly \hbox{adjacent half-lines}.

A \emph{regular quadrangle} is any 2-cell which is bounded by one red source $W_-$\,, one blue sink $W_+$\,, their two red and two blue blow-up and blow-down orbits defining a sector for each, and the two segments between their associated saddles on the boundary circle $\mathbb{S}^1$.
Note adjacency of the boundary saddles, along their boundary segments.
The interior of any quadrangle is filled with heteroclinic orbits from $W_-$ to $W_+$\,.

A \emph{leaf quadrangle} refers to the degenerate case where one of $W_\pm$\,, say the red source $W_-$\,, only possesses a single red blow-up to $\mathbb{S}^1$.
An analogous ``blue'' blow-down definition applies to sinks $W_+$\,.
We may safely skip the most degenerate case where both $W_\pm$ possess only a single sector.
By adjacency of their two associated boundary saddles on $\mathbb{S}^1$ that case reduces to the quadratic of section \ref{ODE2}. 

All quadrangles together, both of regular and of leaf type, partition the closed disk $\mathbf{D}$ into a \emph{quadrangulation}.

We can now define the planar \emph{tree portraits} $\mathcal{T}$ which schematically symbolize the global phase portraits  at the heart of correspondence theorem \ref{Classthm}.
Their \emph{vertices} are the source and sink equilibria $e_0,\ldots,e_{d-1}$\,.
Single \emph{edges} between pairs $W_\pm$ of $e_j$ represent the existence of heteroclinic orbits $W_-\leadsto W_+$\,.
We do not keep track of heteroclinic time direction, at this moment.
In other words, each quadrangle is represented by an undirected edge.
Since each quadrangle extends to the circle boundary, the resulting graph $\mathcal{T}$ is cycle-free, hence a tree.
Moreover, the quadrangle associated to any pair $W_\pm$ consists of all heteroclinic orbits $W_-\leadsto W_+$\,.
The tree $\mathcal{T}$ is connected because the quadrangulation is connected, as a partition of $\mathbf{D}$.
By construction, the tree portrait is unlabeled, unrooted, undirected, and planar.

Abstractly, undirected (connected) trees possess a unique alternating 2-coloring of vertices, up to interchanging the two colors.
Each edge meets both colors.
Of course, this allows us to direct edges, away from one distinguished color towards the other.
Vertices of our tree portraits $\mathcal{T}$ are then colored red and blue, alternatingly along paths.
The colors determine the heteroclinic time direction of the edges.
The two 2-colorings correspond to time reversal.
The tree $\mathcal{T}$ of colored sources and sinks, with heteroclinic edges directed by real time, has been called a \emph{connection graph}. 
For extensive studies of this concept, in the context of parabolic PDEs much more general than \eqref{PDEw}, see for example \cite{brfi88,brfi89,firoSFB, firoFusco} and the many earlier references there.
For  historical overviews and many further illustrations and applications, see for example \cite{Conley, Mischaikow,Yorke-a, Yorke-b}.

Quadrangulations are closely related to handshakes, i.e. to \emph{chord diagrams} of $d-1$ non\-intersecting chords of the closed unit disk $\mathbf{D}$, up to rotation.
To obtain handshakes, we split each quadrangle in two, by an arc $\gamma$ perpendicular to the one edge of the tree portrait in the quadrangle.
We extend $\gamma$ to the midpoints of the two circular boundary segments of the quadrangle.
Note that the \emph{splitting arc} $\gamma$ can be chosen transversely to the compactified vector field.
Since quadrangles only intersect along their interior heteroclinic boundaries, in $\mathbb{D}$, the splitting arcs are disjoint.
We obtain $d-1$ arcs, one for each edge of the connected tree portrait.
Homeomorphically, we can standardize the arc diagram to be represented by straight line chords between equidistant points $\beta_k:=\pi (k+\tfrac{1}{2})/(d-1),\ 0\leq k<2(d-1)$, on the unit circle $\mathbb{S}^1$.
The chord diagram defines a non-crossing handshake between the $2(d-1)$ midpoints $\beta_j$\,.

Conversely and dually, each chord diagram defines a quadrangulation.
Consider the decomposition of the disk $\mathbf{D}$ by the chord diagram.
We just have to mark each of the complementary $d$ faces by a vertex, and connect each vertex to each segment of its face on the circle boundary.
This results in a uniquely defined quadrangulation, once we omit the $d-1$ original chords.
The two interior vertices of source/sink type on the two interior boundaries of any resulting quadrangle define an edge between them, which splits the quadrangle perpendicularly to its removed chord.
Formally, this generates a tree portrait and shows the equivalence of claims \emph{(i)} and \emph{(ii)} in theorem \ref{Classthm}.

It is obvious by construction that $C^0$ orbit equivalent Poincaré compactifications generate equivalent tree portraits.
In the next section \ref{OrbEqu}, we prove that, conversely, compactified Poincaré flows with equivalent tree portraits are $C^0$ orbit equivalent.

\subsection{Orbit equivalence}\label{OrbEqu}

In sections \ref{ODEper} and \ref{SadSad}, we have seen how the strong nondegeneracy condition \eqref{sumsetaRe} implies absence of Lyapunov centers and interior saddle-connections.
We have also associated planar tree portraits to the Poincaré compactifications of ODEs \eqref{ODEw}.
In the present section we prove that ODEs with ``the same'', i.e. equivalent, tree portraits are in fact $C^0$ orbit equivalent.
This will prove correspondence theorem \ref{Classthm}\emph{(i)}.

\emph{Proof of theorem \ref{Classthm} (i)}.
Consider any two univariate polynomials $f_1, f_2$ of the same degree $d$, which satisfy strong nondegeneracy condition \eqref{sumsetaRe} and are described by tree portraits which are homeomorphic by orientation preserving homeomorphisms of their planar embeddings into the closed disk $\mathbb{D}$ of their Poincaré compactifications.
Comparing their quadrangulations allows us to construct orientation preserving homeomorphism which map quadrangles of $f_1$ to quadrangles of $f_2$\,, for corresponding vertex pairs $W_-\leadsto W_+$ in their tree graphs.
Here we may have to perform a global reversal of time first, to match the red/blue colorings associated to $f_1, f_2$.

We fix the homeomorphism on the boundary of each quadrangle.
It remains to extend that homeomorphism to an orbit equivalence, separately for each pair of corresponding regular quadrangles. We skip the analogous case of leaf quadrangles.

We follow the well-trodden path of proofs for structural stability of Morse-Smale systems; see \cite{PalisSmale, Palis, PalisdeMelo, Sotomayor}, with minor modifications.
First, we are dealing with the simplest planar Morse case, in absence of periodic orbits.
Second, however, we have to address the nontransverse boundary saddle-connections within the invariant circular boundary segments in $\mathbb{S}^1$.
Third, we have to preserve the prescribed homeomorphisms on the blow-up and blow-down boundaries of the two quadrangles, in the interior disk $\mathbb{D}$.

Standard Grobman-Hartman arguments linearize, locally at the source, sink and saddle equilibria \cite{Hartman}.
Next note that end points $\mathbf{A}_1, \mathbf{A}_2$ of a boundary circle segment connect such that
\begin{equation}
\label{WAAW}
W_-\leadsto \mathbf{A}_1 \leadsto  \mathbf{A}_2 \leadsto W_+
\end{equation}
is a concatenation of heteroclinic orbits, in the same real time direction.
We first construct our orbit equivalence mapping interior source/sink heteroclinic orbits $W_-\leadsto W_+$ onto each other, for nonlinearities $f_1, f_2$\,.
At real time $r=0$, we start from a homeomorphism of splitting arcs, transverse to the flow of each quadrangle; see section \ref{StarQuad}.
Trajectories near the circular boundary segments, however, will enter neighborhoods of linearization at the boundary saddles $\mathbf{A}_1, \mathbf{A}_2$\,.
The closed neighborhoods $U_1, U_2$ of linearization at these hyperbolic corner points can be chosen to consist of orbit segments: none of these neighborhoods is entered twice.
This allows us to define orbit preserving homeomorphisms in quadrangles up to, and including their consistent heteroclinic boundaries.
Since even simpler linearization arguments apply near $W_\pm$, this completes the proof of theorem \ref{Classthm}\emph{(i)}.
\hfill $\bowtie$.

\emph{Proof of theorem \ref{Classthm} (ii)}.
The equivalence of claims \emph{(i)} and \emph{(ii)} in theorem \ref{Classthm} has already been established at the end of section \ref{StarQuad}.
This proves theorem \ref{Classthm}.
\hfill $\bowtie$.

\subsection{Bifurcation graphs}\label{BifGraph}

In this section, we discuss some bifurcation aspects for univariate polynomials of any degree $d$ which violate the strong nondegeneracy condition \eqref{sumsetaRe}.
We keep the defining nondegeneracy condition \eqref{sumseta} of the connected space $\mathbf{Y}'$ of equilibrium configurations $\mathbf{e}\in \mathbf{f}^{-1}(\mathbf{Y}')$ intact.
For notation and background see section \ref{f'Pf}, and in particular \eqref{Y'} and corollary \ref{finvY'}.

To formalize our approach, we introduce \emph{bifurcation graphs} $\mathcal{B}_d$\,.
\emph{Vertices} of $\mathcal{B}_d$ correspond to orbit equivalence classes  of Poincaré compactifications for those polynomial ODEs \eqref{ODEw} which satisfy the strong nondegeneracy condition \eqref{sumsetaRe}.
Equivalently, we will view each vertex as a tree portrait $\mathcal{T}$ of a Poincaré compactification, by an unlabeled, unrooted, undirected, planar tree with $d$ vertices, up to homeomorphisms of the plane which preserve planar orientation (but may reverse orientation of the undirected edges).
Deviating slightly from previous conventions, we do distinguish vertices of $\mathcal{B}_d$ related by time reversal here.

Perhaps it is worth noting how the vertex set of $\mathcal{B}_d$ is nonempty, for $d\geq 1$.
When all equilibria $e_j$ are real, for example, the associated tree portrait $\mathcal{T}$, or source/sink connection graph, is a line with only two leaves.
In the bifurcation graph $\mathcal{B}_d$\,, this example provides a vertex with minimal absolute difference in the counts of sources and sinks.
That minimal difference count is 0, for even $d$, and 1 when $d$ is odd.
The cyclotomic cases $f(w)=w^d-1$ provide further examples, unless two Lyapunov centers $e_{d/4},e_{3d/4}$ arise because $d$ is an integer multiple of 4.
See figure \ref{fig2}, as always, and section \ref{Cyclotomic}.

Star-shaped tree portraits, with one central vertex of the tree and $d-1$ leaves, provide two more vertices ``$\star$'' of $\mathcal{B}_d$\,.
Here $\star$ denotes the planar star-tree $\mathcal{T}$ of $d\geq 2$ vertices, with one central vertex of edge-degree $d-1$, and $d-1$ attached terminal leaves.
For star-trees, the difference of source versus sink counts is maximal, i.e. $d-2$, in absolute value.
Time reversal actually identifies two star-trees: one where the central vertex is the only source, and the other where that vertex is the only sink.
The case $d=4$, for example, arises from the Poincaré compactification in figure \ref{fig2}(c) by a perturbation which pushes both Lyapunov centers $e_1, e_3$ to become sinks, like $e_2$ already is.
The central vertex $e_0$ remains the only source.
The general star-trees arise in the course of connectedness theorem \ref{Bdthm} below.

\emph{Edges} of the bifurcation graph $\mathcal{B}_d$ will keep track of violations of strong nondegeneracy \eqref{sumsetaRe}.
Violations of \eqref{sumsetaRe}, we recall, require reciprocal derivatives $\eta_j$ which satisfy
\begin{equation}
\label{sumsetabif}
 \sum_{j\in J}\,\Re\,\eta_j\,=\,0\,,
\end{equation}
for some nonempty subset $\emptyset\neq J \subsetneq\{0,\ldots,d-1\}$.
So far, we know that such violations are necessary conditions for the occurrence of interior saddle-connections, and, in some special cases, necessary and sufficient conditions for Lyapunov centers.
See lemmata \ref{perlem} and \ref{sad2lem1}, \eqref{sumsetaRe0}.
In the spirit of corollary \ref{finvY'}, the $\mathbf{e}\in \mathbf{Y}'$ which satisfy any bifurcation condition \eqref{sumsetabif} define a finite collection of submanifolds of $\mathbf{Y}'$, but this time of \emph{real} codimension one.
Let any local curve of equilibrium configurations $\mathbf{e}\in \mathbf{f}^{-1}(\mathbf{Y}')$ transversely cross a maximal stratum, of local real codimension one, of this variety.
Then we draw an \emph{edge} in the bifurcation graph between those vertices which represent the tree portraits of the Poincaré compactifications on either side of the separating codimension one condition \eqref{sumsetabif}.
We omit self-looping edges, in case those tree portraits are identical. 

In lemma \ref{sad2lem1}, we have seen how interior saddle-connections imply degeneracies \eqref{sumsetabif}. To clarify the role of edges in the bifurcation graph, we now establish a partial converse, which derives existence of a unique interior saddle-connection from disconnectedness of tree portraits $\mathcal{T}$.

\begin{lem}\label{sad2lem2}
Let $J$ and its complement $J^c\ni 0$ be a nonempty partition of the index set $j\in\{0,\ldots,d-1\}$ such that degeneracy \eqref{sumsetabif} occurs. 
Assume that each of the two tree portraits of interior source/sink connection graphs among $J$ and among $J^c$, separately, is connected.\\
Then there exists a unique interior saddle-connection $\gamma$, for the Poincaré compactification of \eqref{ODEw}, with $e_j$ on opposite sides of $\gamma$ for $j\in J$ and $j\in J^c$, respectively.
\end{lem}

\begin{proof}
To prove \emph{existence} of $\gamma$, we argue indirectly.
Assume absence of any interior saddle-connections.
By lemma \ref{perlem} and the Sotomayor remark, absence of interior saddle-connections implies absence of Lyapunov centers. 
In particular, the Poincaré compactification is described by a connected tree portrait $\mathcal{T}$.
Our analysis in section \ref{StarQuad} then implies that $\mathcal{T}$ is connected, even in presence of degeneracy \eqref{sumsetabif}.
The parts of $\mathcal{T}$ among vertices in $J$ and $J^c$ are also connected, each, by assumption.
Since the total tree portrait $\mathcal{T}$ is acyclic, it possesses a unique edge between some distinguished vertex $j\in J$ and its counterpart $j^c\in J^c$.
We call this edge a \emph{bridge} between the two parts $J$ and $J^c$.
The bridge defines a bridging quadrangle between its source/sink pair.

As in section \ref{StarQuad}, let the arc $\gamma$ split the bridging quadrangle, transversely to the (compactified) vector field $f$ in it.
Contrary to what we seek to prove, indirectly, $\gamma$ cannot be an orbit.
In $\widehat{\mathbb{C}}^d$, the splitting arc $\gamma$ becomes a closed loop through $w=\infty$.
For the integration of $\omega=dw/f(w)$ our construction of $\gamma$, the residue theorem, and degeneracy assumption \eqref{sumsetabif} therefore imply
\begin{equation}
\label{oint=0}
\Im \oint_\gamma\, \omega \,=\, \pm\Im \left( 2\pi\mi\, \sum_{j\in J}\, \eta_j\right) \,=\, \pm2\pi\, \sum_{j\in J} \Re\, \eta_j \,=\,0\,.
\end{equation}
Transversality of $\gamma$ to $f$, on the other hand, implies a constant sign of the integrand
\begin{equation}
\label{ointnot0}
\lvert f\rvert^2\, \Im \omega = \lvert f\rvert^2\, \Im (dw/f) =  \Im (\bar{f}\cdot dw) = \Re (-\mi\bar{f}\cdot dw)=\langle\mi f, dw\rangle \neq 0\,.
\end{equation}
Indeed, the scalar product on the right hand side keeps track of the nonvanishing sine of the angle between the tangent direction $dw$ of $\gamma$ and the supposedly transverse vector field $f$.
The contradiction \eqref{oint=0}, \eqref{ointnot0} proves the existence of an interior saddle-connection $\gamma$, along the vector field $f$, rather than the existence of a splitting arc $\gamma$ transverse to it.

We prove \emph{uniqueness} of the interior saddle-connection, next.
By assumption, the separate tree portraits of $J$ and $J^c$ are each connected.
Therefore, the blow-up and blow-down orbits of their sinks and sources connect to a total of $2(|J|-1)+2(|J^c|-1)=2(d-2)$ boundary saddles on $\mathbb{S}^1$.
If Lyapunov centers are absent, this leaves precisely one pair of boundary saddles for interior saddle-connections $\gamma$. This proves uniqueness of $\gamma$.

Consider a Lyapunov center $e$, at last, which is neither  source nor sink. 
Then our assumptions imply that $J=\{j\}$ or $J^c=\{0\}$ is a singleton, without any connecting edges in its part of the tree portrait $\mathcal{T}$.
Moreover, the Lyapunov center is unique, for degrees $d\geq 3$.
The interior saddle-connection $\gamma$ then coincides with the unique interior saddle-connection which bounds the nest of periodic orbits around $e$.
For the trivial quadratic case $d=2$, see section \ref{ODE2} in imaginary time direction, and the orange orbits in figure \ref{fig1}.

This proves the lemma.
\end{proof}

See \cite{Sotomayor} for an interesting exploration of such transitions, in the much broader context of generic one-parameter families of planar Morse-Smale flows.
Although our constructions of quadrangulations, tree portraits, and handshakes in section \ref{StarQuad} make it clear how all trees can be realized in the class of all general planar ODE nonlinearities $f$, this is not obvious in our restricted class of Poincaré compactifications of complex scalar polynomial ODEs \eqref{ODEw}.

Instead, we only show connectivity of the bifurcation graph $\mathcal{B}_d$\,.
Specifically, we establish paths from any vertex in $\mathcal{B}_d$ to the two star-trees.

\begin{thm}\label{Bdthm}
For any $d\geq 2$, the bifurcation graph $\mathcal{B}_d$ is connected.
More precisely, there exists a path from any vertex $\mathcal{T}$ of $\mathcal{B}_d$ to each of the two $\star$ vertices.
\end{thm}

\begin{proof}
Let the tree portrait $\mathcal{T}$ denote any vertex of the bifurcation graph $\mathcal{B}_d$\,.
Let $\star$ denote the star-tree with a source as its central vertex.
Since $\mathbf{Y}'$ is connected, we can choose a generic path of the vector $\boldsymbol{\eta}=(\eta_1,\ldots,\eta_{d-1})$ of reciprocal derivatives, from its value at $\mathcal{T}$ to its value at $\star$\,.
The nonvanishing determinants \eqref{Dd1fact} in the proof of lemma \ref{fbopen} and corollary \ref{finvY'} allow us to pull back the path of $\boldsymbol{\eta}$ to a path of $\mathbf{e}$.

Of course, the path of $\boldsymbol{\eta}$ cannot avoid Lyapunov centers, unless $\mathcal{T}=\star$\,.
Indeed, the total number of sources must decrease to one, along the path of $\boldsymbol{\eta}$.
Along the way, unavoidably, other bifurcations of saddle-saddle type might occur, by violation \eqref{sumsetabif} of strong nondegeneracy condition \eqref{sumsetaRe}.
Because the violation is of real codimension one, however, the violations will only occur at isolated points, generically, where the path $\boldsymbol{\eta}$ crosses strata of maximal dimension transversely.
By definition, each such crossing amounts to traversing an edge in the bifurcation graph $\mathcal{B}_d$\,.
This shows that any vertex $\mathcal{T}$ connects to the star-tree with source central vertex, in each bifurcation graph $\mathcal{B}_d$\,.
Increasing the number of sources to $d-1$, in contrast, connects to the other star-tree with sink central vertex.

In particular, these constructions establish the existence of Poincaré compactifications of either star-tree type.
Moreover they show that each bifurcation graph $\mathcal{B}_d$ is connected, proving the theorem.
\end{proof}

Although small degrees $d$ are now accessible to hands-on inspection, the general structure of the bifurcation graphs $\mathcal{B}_d$ is not known.
For  a specific example, consider time rays $r\mapsto r\exp(\mi\theta)$ which are inclined by a fixed angle parameter $\theta$ against the real time axis of $\theta=0$.
At $\theta=\pi$, real time is reversed.
The vector $\boldsymbol{\eta}$ of reciprocal derivatives at equilibria $e_j$ co-rotates as $\exp(\mi\theta)\boldsymbol{\eta}$\,.
For generic $\boldsymbol{\eta}$, it would be interesting to detail the resulting cycles in the bifurcation graph $\mathcal{B}_d$\,, which are generated by the edges associated to the co-rotated degeneracies \eqref{sumsetabif}.
Even the nongeneric cyclotomic cases are of interest.

\subsection{On realization theorem \ref{Realizethm}}\label{RealizethmComm}

Beyond correspondence theorem \ref{Classthm} and connectedness theorem \ref{Bdthm}, theorem \ref{Realizethm} claims realization of \emph{all} planar trees $\mathcal{T}$, as source/sink connection graphs of Poincaré compactifications for ODE \eqref{ODEw}.
The techniques developed above suggest induction over the degree $d$ of the polynomial $f$, i.e. over the number $d$ of vertices in the planar tree portraits $\mathcal{T}^d$ of the source/sink connection graphs.
The cases $d=2,3$ of figures \ref{fig1} and \ref{fig2} (a) are trivial, because the corresponding trees are unique.
Indeed, theorem \ref{Classcor} then asserts the counts $A_{d-1}=1$.

In general, any tree possesses at least two leaves.
Plucking any leaf from a planar tree $\mathcal{T}^{d+1}$ therefore leaves us with a planar tree $\mathcal{T}^d$.
Conversely, if we manage to grow a leaf to an arbitrary planar tree $\mathcal{T}^d$, in any of the $2(d-1)$ positions which $\mathcal{T}^d$ has to offer, then we obtain (perhaps redundantly) any planar tree $\mathcal{T}^{d+1}$.
Since Lyapunov centers $J=\{j\}$ can (only) appear at leaves of the tree, this particular variant of degeneracy \eqref{sumsetabif} is the strategy which we will sketch next.
See section \ref{BifGraph} and particularly the proof of lemma \ref{sad2lem2}.
For complete details, we refer to \cite{FiedlerYamaguti}.

\emph{Proof of theorem \ref{Realizethm}, sketch.}
To be specific, suppose the planar connection graph $\mathcal{T}^d$ arises from the Poincaré compactification of $\dot w=g(w)=(w-e_0)\cdot\ldots\cdot(w-e_{d-1})$, under the usual nondegeneracy assumption \eqref{sumsetaRe}.
We now consider the homotopy
\begin{equation}
\label{ftheta}
\dot w = f^\theta(w):=g(w)\cdot(1-\eps \exp(\mi\theta)w).
\end{equation}
for $0\leq\theta< 2\pi$. 
Although we could rescale $w$ to render $f^\theta$ univariate, we prefer the more direct analysis of \eqref{ftheta}, for clarity of presentation.
We will choose $\eps>0$ small.
Evidently, the polynomial $f^\theta$ of degree $d+1$ inherits the equilibrium zeros $e_0,\ldots,e_{d-1}$ from $g$, i.e. the vertices of $\mathcal{T}^d$.
The additional zero $e_d^\theta=1/\eps \exp(-\mi\theta) $ tracks a large circle around the inherited vertices, clockwise.

For $j<d$, the reciprocal linearizations $\eta_j=1/f'(e_j)=1/g'(e_j)+\ldots$ are also inherited, up to a perturbation of order $\eps$.
Moreover, the structurally stable connection graph $\mathcal{T}^d$ is contained in $\mathcal{T}^{d+1}$.

At the large equilibrium $e_d^\theta$, we obtain the linearization
\begin{equation}
\label{etatheta}
1/\eta_d=(f^\theta)'(e_d^\theta)=-\eps \exp(\mi\theta) g(e_d^\theta)=-\eps^{-(d-1)}(\exp(-\mi(d-1)\theta)+\ldots).
\end{equation}
Therefore, interior saddle connections cannot arise, except possibly via degeneracy \eqref{sumsetabif} at $J=\{d\}$, i.e. via Lyapunov centers $e_d^\theta$\,.
See  our remark at the end of section \ref{ODEper}, and lemmata \ref{sad2lem1}, \ref{sad2lem2}. 

Hopf transitions at Lyapunov centers $e_d^\theta$ occur at $\theta=\theta_k=\pi (\tfrac{1}{2}+k)/(d-1)+\ldots$\,, for $0\leq k< 2(d-1)$. 
The $2(d-1)$ Lyapunov centers $e_d^\theta$\,, at $\theta=\theta_k$\,,  are accompanied by interior saddle connections $\Gamma_k:\mathbf{A}_k\leadsto \mathbf{B}_k$ between adjacent saddles $\mathbf{A}_k\,,\, \mathbf{B}_k$ on the Poincaré boundary $\mathbb{S}^1$.
See figure \ref{fig2} (c) for a double example.
The local unfoldings of each Lyapunov center alternate between source and sink type as $\theta$ increases through $\theta=\theta_k$\,.
The linearization $(f^\theta)'(e_d^\theta)$ in eqref{etatheta} rotates around zero, clockwise.
In particular, the rotation directions of the periodic orbits around successive Lyapunov centers $e_d^{\theta_k}$ are of alternating sign.

As $\theta$ and $k$ increase, the adjacent saddles $\mathbf{A}_k\,,\, \mathbf{B}_k$ therefore ratchet through $2(d-1)$ such pairs on the boundary $\mathbb{S}^1$, clockwise and cyclically.
The pairs successively attach $e_d^\theta$ as a leaf of $\mathcal{T}^{d+1}$ to any of the $2(d-1)$ possible positions which $\mathcal{T}^d$ has to offer. 
Since the planar tree $\mathcal{T}^d$ was arbitrary, the parameter $\theta$ realizes all planar trees $\mathcal{T}^{d+1}$; again see \cite{FiedlerYamaguti} for full details.

This completes our sketch of a proof of realization theorem \ref{Realizethm}. \hfill $\bowtie$

\section{Discussion}\label{Dis}

\subsection{Overview}\label{Over}

We conclude with a cursory collection of some further scalar and higher-dimensional aspects.
In section \ref{ODEhet}, we consider complex scalar ODEs $\dot{w}=f(w)$ with complex entire, rather than just polynomial, nonlinearities $f$.
We show how heteroclinic and homoclinic orbits between not necessarily hyperbolic equilibria $W_\pm$ still have to blow up in finite time.
In higher dimensions, we also recall some classical results on non-entire solutions to complex entire ODEs.

Section \ref{1000} follows \cite{FiedlerClaudia} to describe the potential relevance of complex entire homoclinic orbits, even in ODEs, for the elusive phenomenon of ultra-exponentially small homoclinic splittings and ultra-invisible chaos.
The existence of such homoclinic orbits remains a 1,000 \euro\ open question.

Our final section \ref{EntirePer} recalls a reversible scalar ODE example of second order.
The example features a non-entire homoclinic orbit which coexists with a complex entire periodic orbit (a mere cosine, in fact).
As a consequence, ultra-exponential homoclinic splitting may fail, but ultra-exponentially  sharp Arnold tongues emerge under suitable, rapidly periodic forcing.

\subsection{Blow-up in complex entire ODEs}\label{ODEhet}

For complex entire, rather than just polynomial, ODEs $\dot{w}=f(w)$, let us consider $\Gamma:W_-\leadsto W_+$ which are heteroclinic between any two complex equilibria $f(W_\pm)=0$, in real time $t=r$.
We do not assume $W_\pm$ to be simple zeros of $f$, and we explicitly allow $\Gamma$ to be homoclinic, i.e. $W_+=W_-$\,.
We claim $\Gamma(t)$ itself cannot be entire, in complex time $t$.
 
If $W_\pm$ are hyperbolic, i.e.~for $\Re f'(W_\pm)\neq 0$, the claim follows from \cite{FiedlerClaudia}.
In the general, possibly nonhyperbolic, scalar case, again, we argue indirectly and assume $\Gamma$ to be complex entire.
To reach a contradiction, we first show that $\Gamma$ then possesses a complex period $p$.

Indeed, complex periodicity follows from the great Picard theorem, as follows.
Let $\Gamma_0:=\Gamma(0)$.
Since heteroclinic orbits $\Gamma(t)$ cannot be polynomial, the holomorphic map 
$0\neq t \mapsto\Gamma(1/t)$ possesses an essential singularity at $t=0$.
The great Picard theorem then implies that $\Gamma(t)$ attains any complex value infinitely many times, of course with the single equilibrium exception $W_+=W_-$\,.
In particular, the (supposedly entire and supposedly heteroclinic) orbit $\Gamma(t)$ must then be homoclinic to $W=W_\pm$ in real time $t=r$, with a complex period $p$ of nonvanishing imaginary part.

We may therefore rescale time to normalize the complex period to $p=\pi\mi$, at the expense of ``real'' time rays $r\exp(\mi\theta),\ r\geq0$ now progressing at some fixed angle $0\leq|\theta|<\pi/2$ to the imaginary time axis.
By lemma \ref{perlem}\emph{(ii)}, a resulting $\pi$-periodic orbit of $\Gamma$ in rescaled imaginary time $t=\mi s$ surrounds a unique equilibrium, say $f(0)=0$.
Lemma \ref{perlem}\emph{(iv)} identifies $e=0$ as a nondegenerate Lyapunov center, in rotated imaginary time.
Therefore $W_\pm=e=0$ is a source or sink in original time.
This contradiction to homoclinicity of $\Gamma$ proves that scalar heteroclinic orbits $\Gamma$ cannot be entire.

Let us briefly comment on the case of ODEs $\dot{w}=f(w)\in\mathbb{C}^N$ in several complex variables.
Consider heteroclinic orbits $\Gamma:W_-\leadsto W_+$ between hyperbolic equilibria $W_\pm$\,.
For recent blow-up results in the spirit of section \ref{Para}, we can refer to \cite{FiedlerClaudia}, theorem 1, and the references and discussion there.
Beyond hyperbolicity, the unstable part of the spectrum of $f'(W_-)$ is assumed to be real and, among itself, non-resonant.
At the target equilibrium $W_+$, analogously, real and non-resonant spectrum is assumed in the stable part.
Non-resonance, in the sense of Poincaré, forbids eigenvalues to arise as integer combinations in the respective spectral parts.
The significance of non-resonance lies in \emph{Poincaré's theorem on locally holomorphic linearization}; see \cite{Ilyashenko} for a complete proof.
For a scalar example of Poincaré linearization, see lemma \ref{perlem}\emph{(v)}.

The argument in \cite{FiedlerClaudia} can be summarized as follows.
Poincaré linearization with non-resonant real eigenvalues implies local quasi-periodicity, in imaginary time $s$.
For complex entire homoclinic or heteroclinic orbits $\Gamma=\Gamma(r+\mi s)$, we can then invoke standard results on almost-periodicity \cite{Bohr, Besicovitch, Corduneanu}, and compare quasi-periodic Fourier coefficients in $s$, globally for large positive and negative real times $r$, to reach a contradiction.
This establishes blow-up in complex time.

In the quadratic parabolic PDE \eqref{PDEw} of section \ref{Para} and \cite{fiestu24}, infinite-dimensional non-resonance at $W^+$ cannot be brought to bear, because Poincaré linearization is unavailable in infinite dimensions.

A much less explicit foreboding of the above ODE results, for unstable manifolds of complex entire diffeomorphisms in $\mathbb{C}^2$, can be attributed to Ushiki \cite{Ushiki-2}.
In higher finite dimensions $\mathbb{C}^N$, see also \cite{Ushiki-N} for the special heteroclinic case of unstable dimensions 1, at the source $W_-$\,, and stable dimension 1 at the target $W_+$\,.

For the special case of scalar second order pendulum equations 
\begin{equation}
\label{pendulum}
\ddot w + \mathbf{h}(w) = 0
\end{equation}
with nonlinear complex entire $\mathbf{h}$, already Rellich \cite{Rellich} in 1940 has most elegantly noticed the complete absence of \emph{any} non-constant complex entire solutions $w(t)$; see also \cite{Wittich1}.
Wittich has generalized his negative answer to scalar nonautonomous $m$-th order equations which are a polynomial $\mathbf{p}$ in $t,w,\dot w,\ldots,w^{(m)}$ plus a complex entire nonlinearity $\mathbf{h}(w)$:
\begin{equation}
\label{pendulumm}
\mathbf{p}(w^{(m)},\ldots,w,t)+\mathbf{h}(w)=0\,.
\end{equation}
See \cite{Wittich2}.
His setting includes the purely polynomial case $\mathbf{h}=0$.
The only complex entire solutions $w(t)$ are then polynomial in $t$, and hence neither homo- or heteroclinic, nor non-constant periodic.
The Wittich result requires an elementary additional condition on $\mathbf{p}$. 
That condition is violated in our autonomous, time reversible, second order example \eqref{oderev} below, which features a periodic cosine solution \eqref{revper}.

\subsection{Ultra-invisible chaos: a 1000 \euro\ question}\label{1000}

Even in pure ODE settings, complex time extensions of real analytic heteroclinic and homoclinic orbits $\Gamma$ hold interesting promise and applied interest.
We briefly sketch some recent progress.
For details and many further references, we refer to \cite{FiedlerScheurle, FiedlerClaudia}.
To be specific, consider rapid $\eps$-periodic forcings $g$ of autonomous ODE systems $f$, as in
\begin{equation}
\label{odefg}
\dot{w}(t)=f(\lambda,w(t)) + \eps^p g(\lambda,\eps,t/\eps,w(t))\,.
\end{equation}
The functions $f,g$ are assumed to be analytic for small real $\lambda,\,\eps$, all $w\in\mathbb{C}^N$, and real for real $w$. 
In the real variable $\beta=t/\eps$ of the rapid forcing $g=g(\lambda,\eps,\beta,w)$, we only assume smoothness and real period 1.
Throughout, the origin $W=0$ is assumed to be a hyperbolic equilibrium.

\begin{figure}[t]
\centering \includegraphics[width=0.6\textwidth]{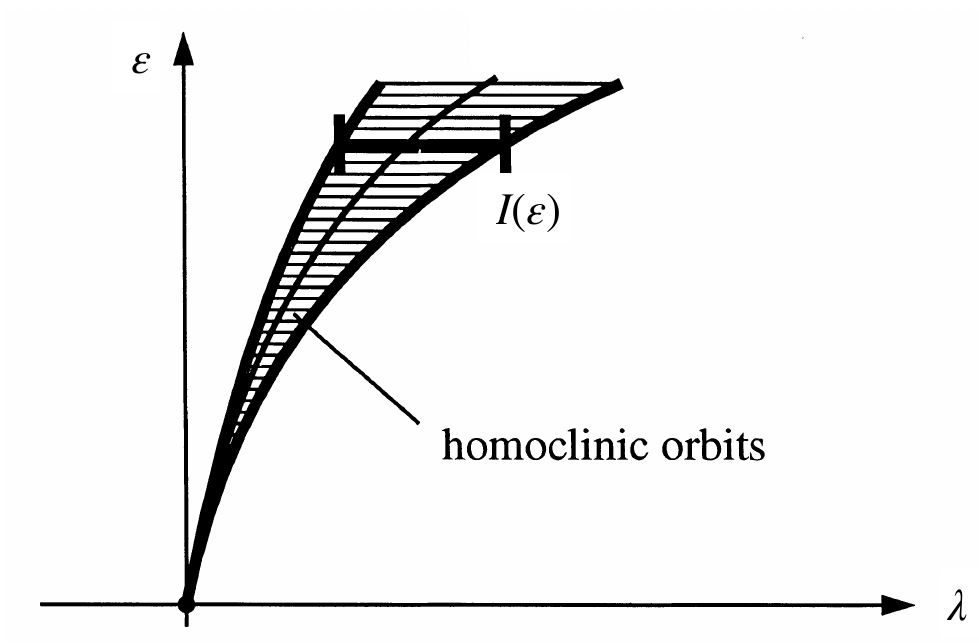}
\caption{\emph{
Schematic splitting region for a homoclinic orbit $\Gamma$ of the ODE-flow \eqref{odef} (hashed), under discretization with step size $\eps$ or, equivalently, under $\eps$-periodic non-autonomous forcing \eqref{odefg}; see \cite{FiedlerScheurle}.
At fixed levels of $\eps$, the horizontal splitting intervals $\lambda\in I(\eps)$ mark parameters $\lambda$ for which single-round homoclinic orbits occur near $\Gamma$.
For analyticity of $\Gamma(r+\mi s)$ in a strip $|s|\leq a$,  the width $\ell(\eps)$ of the splitting interval $I(\eps)$ is exponentially small in $\eps$ with exponent $a$; see \eqref{split}.
For better visibility, the horizontal width $\ell(\eps)$ of the exponentially flat splitting region has therefore been much exaggerated, in our schematic illustration.
We call the dynamics in the resulting chaotic region ``invisible chaos''.
For complex entire homoclinic orbits $\Gamma(t)$, the exponent $a$ could be chosen arbitrarily large.
Such ultra-exponentially small splittings would therefore lead to ultra-invisible chaos.\\
}}
\label{fig3}
\end{figure}

The time-$\eps$ stroboscope map of \eqref{odefg} may, equivalently, be seen as a one-step discretization scheme, of step size $\eps$ and order $p$, for the autonomous ODE
\begin{equation}
\label{odef}
\dot{w}=f(\lambda,w).
\end{equation}
This topic runs under the name of \emph{backward error analysis}; see for example \cite{MoserPhysics, Reich, WulffOliver} and the more detailed references in \cite{FiedlerClaudia}.

Consider a nondegenerate homoclinic orbit $\Gamma:0\leadsto 0$ of \eqref{odef}.
Here \emph{nondegeneracy} requires that real stable and unstable manifolds of the hyperbolic equilibrium $W=0$ cross each other at nonvanishing speed in $\mathbb{R}^N$, as the scalar real parameter $\lambda$ increases through $\lambda=0$.
Under certain further nondegeneracy conditions on $f$ and the perturbation $g$, detailed in section 5 of \cite{FiedlerScheurle}, this leads to \emph{invisible chaos}: the chaotic region accompanying transverse homoclinics is exponentially thin, both, in parameter space $(\lambda,\eps)$ and along the homoclinic loop.
More precisely, let 
\begin{equation}
\label{astrip}
r\in\mathbb{R},\ |s|\leq a
\end{equation} 
denote a horizontal complex strip where the complex time extension $\Gamma(r+\mi s)$ of the real analytic homoclinic orbit $\Gamma(r)$ of \eqref{odef} remains complex analytic.
Then the width $\ell(\eps)$ of the splitting interval $I(\eps)$, where transverse homoclinicity and accompanying shift dynamics prevail, satisfies an exponential splitting estimate
\begin{equation}
\label{split}
\ell(\eps) \leq C_a \exp(- a/\eps),
\end{equation}
for some constant $C_a>0$ and all small $\eps \searrow 0$.
Remarkably, the splitting region is small of infinite order $\exp(-a/\eps)$ in $\eps$, even under forcings \eqref{odefg} or discretizations of finite order $\eps^p$.
See figure \ref{fig3} for illustration.
For thorough ODE surveys on this topic, going back as far as Poincaré \cite{Poincare3body} and starting technically with \cite{Neishtadt}, see \cite{Gelfreich01,Gelfreich02}.
For PDE extensions see \cite{MatthiesDiss,Matthieshom,MatthiesScheel} and the references there.

The significance of \emph{complex entire homoclinic orbits} $\Gamma$ is now evident.
In the complex entire case, an upper estimate \eqref{split} would hold for any strip \eqref{astrip}, i.e.~for any $a$. 
In fact, we could replace \eqref{split} by an estimate
\begin{equation}
\label{ultrasplit}
-\log \ell(\eps) \geq c(1/\eps) >0 ,
\end{equation}
with a suitable convex function c of unbounded positive slope.
We would call this phenomenon \emph{ultra-exponentially small splitting of separatrices}.
Under further nondegeneracy, such splitting would be accompanied by \emph{ultra-invisible chaos}.

Alas, can it actually happen?
In \cite{FiedlerClaudia}, the author has personally offered a
\medskip\\
\centerline{\textbf{1,000\,\euro\ reward}}
\medskip
for settling this question.

More precisely, a \emph{positive answer} would require any explicit example of a real nonstationary homoclinic orbit $\Gamma(t)$ to a hyperbolic equilibrium, for a complex entire ODE \eqref{odef} on $X=\mathbb{C}^N$ with $f$ real for real arguments, such that  $\Gamma(t)$ is complex entire for $t\in\mathbb{C}$. 
Such an example would initiate many challenges, e.g., for numerical explorations.
Eventually, it might lead towards a whole new theory to address ultra-exponentially small splitting behavior under discretization.

A \emph{negative answer} would prove that such complex entire homoclinic orbits $\Gamma(t)$ cannot exist.
In section \ref{ODEhet}, we have already collected some very partial and incomplete results in that direction. 
However, cases involving complex conjugate, or even just algebraically degenerate real, eigenvalues remain wide open.
This includes the celebrated Shilnikov homoclinic orbits \cite{Shilnikov}.

In order to protect well-established colleagues against their own, potentially intensely distracting, financial interests, the prize will be awarded to the first solution by anyone \emph{without} a permanent position in academia.
Priority is defined by submission time stamp at arxiv.org or equivalent repositories.
Subsequent confirmation by regular refereed publication is required.

\subsection{Ultra-sharp Arnold tongues}\label{EntirePer}

Consider the time reversible, second-order, scalar ODE
\begin{equation}
\label{oderev}
\ddot w +\dot w^2 +w^2-3w=0.
\end{equation}
By integrability, the homoclinic orbit $w(t)=\Gamma(t)$ to the hyperbolic trivial equilibrium $w=0$ satisfies
\begin{equation}
\label{revint}
\tfrac{1}{2}\dot{w}^2=\exp(-2w)-1+2w\,-\,\tfrac{1}{2}w^2\,.
\end{equation}
In particular, $\Gamma(t)$ encounters a non-meromorphic singularity in complex time.
Indeed, theorem 1.1 in \cite{FiedlerClaudia} implies that $\Gamma$ cannot be complex entire.
The term $\exp(-2w)$ in ODE \eqref{revint}, on the other hand, prevents meromorphic singularities.
If the singularity were isolated, it would therefore have to be essential.

The example \eqref{oderev} violates the condition of  \cite{Wittich2} which prevents non-polynomial entire solutions.
In the language of Wittich, each of the \emph{two} terms $\dot w^2$ and $w^2$ is ``of maximal dimension two'', in fact, where only a \emph{single} such term is permitted.

Notably, the nonlinear ODE \eqref{oderev} does possess a \emph{complex entire $2\pi$-periodic orbit} $w$ given by the harmonic function
\begin{equation}
\label{revper}
w(t) := 2+\sqrt{2}\cos t\,.
\end{equation}
Under discretizations, alias rapid forcings \eqref{odefg}, such complex entire periodic orbits provide devil's staircases with ultra-exponentially thin resonance plateaus, and correspondingly ultra-sharp Arnold tongues.
See \cite{ArnoldODE,Devaney,GuckenheimerHolmes} for a general background.
Some further references and details are discussed in section 7 of \cite{FiedlerClaudia}.

Alas, our results have only provided an ultra-exponential upper estimate \eqref{ultrasplit} for Arnold tongues associated to complex entire periodic orbits. 
This raises the question: how ultra-sharp are they, really?
Complementary lower estimates of ultra-sharp Arnold tongues are not known.
For homoclinic, rather than periodic, splittings we even lack any 1,000\ \euro\ ODE example with an ultra-exponential upper estimate.

\newpage
\section{Appendix on Riemann surfaces}\label{Apx}

For convenience we collect some definitions, notations, and facts on Riemann surfaces, following the lines of \cite{Forster, Hartshorne, Jost, Lamotke}.

\emph{Riemann surfaces} are path-connected complex one-dimensional manifolds, with a maximal atlas and holomorphic coordinate changes.
The grand \emph{Riemann mapping theorem}, for example, asserts that (up to biholomorphic equivalence) any simply connected Riemann surface realizes one of the following three cases: the compact Riemann sphere $\widehat{\mathbb{C}}$, the complex plane $\mathbb{C}$, or the open complex upper half plane $\mathbb{H}$.

In general, let $X,Y$ be any two Riemann surfaces.
Holomorphic maps between them are defined via holomorphy in local coordinates.
Any nonconstant holomorphic map
\begin{equation}
\label{pYX}
\mathbf{p}: \ Y\rightarrow X
\end{equation}
is called a \emph{covering}.
I.e., $\mathbf{p}$ is open, and the \emph{fiber} $\mathbf{p}^{-1}(x)\subset Y$ is discrete, for each $x\in X$.
The covering is called \emph{unlimited}, if every base point $x\in X$ possesses a neighborhood $U$ such that the restriction $\mathbf{p}:V_j\rightarrow U$ is holomorphic and surjective (onto), on each connected component $V_j$ of $\mathbf{p}^{-1}(U)$ in $Y$.
(Coverings are often required to be unlimited, by definition.)
The unlimited covering is called \emph{unbranched}, if the restrictions are biholomorphic.
Unlimited unbranched coverings allow the lifting of curves in $X$ to curves in $Y$.
Closed curves in $X$ need not lift to closed curves in $Y$.
The components $V_j\subseteq Y$ are often called \emph{sheets} of $Y$ over $U\subseteq X$.
Some also think of the sheets as the ``decks'' of a parking garage.

In contrast, a point $y_0\in Y$ is called a \emph{ramification point} of $\mathbf{p}$, if $y_0$ does not possess any open neighborhood $V$ in $Y$ such that the restriction of $\mathbf{p}$ to $V$ is injective (one-to-one).
The trace point $x_0=\mathbf{p}(y_0)$ of any ramification point $y_0$ is called a \emph{branch point} of $\mathbf{p}$.
Coverings $\mathbf{p}$ without branch points are called \emph{unbranched}.
In suitable holomorphic local coordinates near a ramification point $y_0$, over a branch point $x_0$, the covering $\mathbf{p}$ can be written in normal form as
\begin{equation}
\label{branchm}
x=\mathbf{p}(y)=x_0+(y-y_0)^m\,,
\end{equation}
for some integer $m\geq 2$, which is called the \emph{branching multiplicity} or \emph{branching degree} of $\mathbf{p}$ at $y_0$.
Locally, for $y_0\neq y\in Y$ near the ramification point $y_0$, the covering $\mathbf{p}$ is therefore $m$-to-one.
Expansion \eqref{sovz}, for example, exhibits local branching of multiplicity $m=d-1$ at the ramification point $w=\infty\in\widehat{\mathbb{C}}_d $ over the branch point $t=0$, for $d\geq 3$.
\emph{Branched coverings} are unlimited coverings for which the set of branching points is discrete in $X$.

Assume the covering $\mathbf{p}: Y\rightarrow X$ is unbranched and unlimited.
\emph{Deck transformations} are homeomorphisms $\varphi$ of $Y$ which satisfy 
\begin{equation}
\label{decktrf}
\mathbf{p}\circ\varphi = \mathbf{p}.
\end{equation}
In other words, deck transformations $\varphi$ preserve each individual fiber $\mathbf{p}^{-1}(x)$, but shuffle the decks or sheets $V_j$ of the covering $\mathbf{p}$.
By composition, deck transformations form the \emph{deck group} $\mathrm{deck}(\mathbf{p})$, sometimes also abbreviated as $\mathrm{deck}(Y/X)$.
A posteriori, deck transformations turn out to be biholomorphic fixed-point free automorphisms of $Y$.
They are in fact liftings of the covering $\mathbf{p}$ itself.
The action of the deck group on $y$ is discontinuous.

An unbranched, unlimited covering $\mathbf{p}: Y\rightarrow X$ is called \emph{normal} or \emph{Galois}, if the deck group acts transitively on each fiber.
In other words, for each pair $y_1,y_2$ in the same fiber $\mathbf{p}(y_1)=\mathbf{p}(y_2)=x$, there exists a deck transformation $\varphi$ such that $\varphi(y_1)=y_2$\,.
Since the deck group acts freely on each fiber, the transformation $\varphi$ is also unique.
We can therefore identify the fiber $\mathbf{p}^{-1}(x)$ with $\mathrm{deck}(\mathbf{p})$.
More precisely, this makes the covering $\mathbf{p}$ a principal fiber bundle over $X$ with structure group $\mathrm{deck}(\mathbf{p})$.

The \emph{universal cover} $\widetilde{X}$ of a Riemann surface $X$, for example, is itself a Riemann surface, with an unlimited, unbranched, normal covering $\mathbf{p}: \widetilde{X}\rightarrow X$.
Since the universal cover $\widetilde{X}$ is simply connected, the grand Riemann mapping theorem identifies $\widetilde{X}$ as one of: the compact Riemann sphere $\widehat{\mathbb{C}}$, the complex plane $\mathbb{C}$, or the open complex upper half plane $\mathbb{H}$.
Accordingly, the original Riemann surface $X$ is then called \emph{elliptic}, \emph{parabolic}, or \emph{hyperbolic}.

The deck group $\mathrm{deck}(\widetilde{X}/X)$ is isomorphic to the \emph{fundamental group} $\pi_1(X)$.
Following Poincaré, the deck groups are called \emph{Fuchsian groups} in the hyperbolic case \cite{Katok}.
Compactness of $X$ is then equivalent to compactness of the fundamental domain of the discontinuous group action on $\widetilde{X}$.

For illustration, let us recall the structure of the punctured spheres $w\in\widehat{\mathbb{C}}_d$\,, as a Riemann surface.
In the simplest case $d=2$, we obtain the parabolic cylinder with simply connected universal cover
\begin{equation}
\label{CcoversC2}
\mathbb{C} \rightarrow \widehat{\mathbb{C}}_2\,.
\end{equation}
Indeed, a biholomorphic fractional linear transformation of the Riemann sphere allows us to consider the case $\widehat{\mathbb{C}}_2\cong \widehat{\mathbb{C}}\setminus \{0,\infty\}\cong\mathbb{C}\setminus\{0\}$, without loss of generality. 
Then the exponential function provides the desired covering projection \eqref{CcoversC2}.
The associated fundamental group is $\pi_1( \widehat{\mathbb{C}}_2)\cong 2\pi\mi\,\mathbb{Z}\cong\mathbb{Z}$.

For $d\geq 3$, we claim that the base $X=\widehat{\mathbb{C}}_d$ is a hyperbolic Riemann surface, i.e. the simply connected universal cover is
\begin{equation}
\label{HcoversCd}
\mathbb{H} \rightarrow \widehat{\mathbb{C}}_d\,.
\end{equation}
Indeed, the cover candidate $\widetilde{X}=\widehat{\mathbb{C}}$ is out because $\widehat{\mathbb{C}}_d$ is noncompact.
Next we recall that any fiber, 
alias the deck group of the universal covering \eqref{HcoversCd}, 
coincides with the fundamental group of $\widehat{\mathbb{C}}_d$\,,
\begin{equation}
\label{pi1Cd}
\pi_1(\widehat{\mathbb{C}}_d)\cong\mathbb{F}_{d-1}\,,
\end{equation}
i.e. with the free group $\mathbb{F}_{d-1}$ of $d-1$ generators.
The generators are small closed loops encircling any of the punctures $e_j\,,\  0<j<d$\,.
The loop around $e_0$ can be omitted.
Indeed, consider any closed Jordan curve $\gamma\subset\mathbb{C}$ which is large enough to contain all zeros $e_j\,, 0\leq j<d$ of $f$ in its interior, including $e_0$\,. 
Contracting $\gamma$ away to $w=\infty\in\widehat{\mathbb{C}}_d$ then provides the only relation $\gamma_0\cdot\ldots\cdot\gamma_{d-1}=\mathrm{id}$ in $\pi_1$\,.

In particular, the second cover candidate $\widetilde{X}=\mathbb{C}$ is also out because the free group $\mathbb{F}_{d-1}$ cannot be a discrete subgroup of the automorphism group $\{z\mapsto aw+b\}$ of $\mathbb{C}$, unless $d=2$.
The Riemann mapping theorem therefore proves claim \eqref{HcoversCd}.

By universality of the universal cover $\widetilde{X}$, any unbranched, unlimited covering $\mathbf{q}: Y\rightarrow X$ gives rise to a sequence of unbranched, unlimited coverings
\begin{equation}
\label{XYX}
\widetilde{X}\xrightarrow{\mathbf{p}} Y\xrightarrow{\mathbf{q}} X\,,
\end{equation}
such that $\mathbf{p}\circ\mathbf{q}: \widetilde{X}\rightarrow X$ is the universal covering of $X$, and $\mathbf{p}$ is the universal covering of $Y$.
The covering $\mathbf{q}$ is normal, if and only if the fundamental group $\mathrm{deck}(\mathbf{p})\cong\pi_1(Y)$ is a normal subgroup of $\mathrm{deck}(\mathbf{p}\circ \mathbf{q})\cong\pi_1(X)$. 
In that case,
\begin{equation}
\label{XYXdeck}
\mathrm{deck}(\mathbf{q})\,\cong\,\mathrm{deck}(\mathbf{p}\circ \mathbf{q})/\mathrm{deck}(\mathbf{p})\,\cong\,\pi_1(X)/\pi_1(Y)\,;
\end{equation}
see \cite{Forster}, theorem 5.9 and exercise 5.2.

Given any holomorphic (and in particular, closed) 1-form $\omega$ on $X$, the integral 
\begin{equation}
\label{tomx}
t=\int_{x_0}^x\, \omega
\end{equation}
depends on the homotopy class in $X$ of any integration path from $x_0$ to $x$; see \eqref{tw}.
The process of maximal analytic continuation, however, defines a Riemann surface $Y$ of holomorphic germs for $t$, and an unbranched holomorphic cover $\mathbf{p}: Y\rightarrow X$, such that the integral $t:Y\rightarrow\mathbb{C}$ becomes holomorphic on $Y$.
Since integration \eqref{tomx} can be defined along any path in $X$, the covering $\mathbf{p}$ is also unlimited.

\emph{Sheaves} provide an abstract framework for describing analytic continuation, e.g. in the Abelian group of holomorphic germs \cite{Forster, Hartshorne, Lamotke}.
In our present setting, however, we decided to stay with the integration case at hand, and not overburden the reader with additional abstract language.



\newpage

\bigskip

\begin{thebibliography}{9999)999}

{\footnotesize{

\bibitem[ADLY]{Yorke-b}
C.~Adwani, R.~De Leo, and J.A.~Yorke.
What is the graph of a dynamical system? 
(2018) \url{https://arxiv.org/abs/2410.05520}

\bibitem[Arn88]{ArnoldODE}
V.I.~Arnold.
\emph{Geometrical Methods in the Theory of Ordinary Differential Equations.}
Springer-Verlag, Berlin 1988.

\bibitem[BG94]{MoserPhysics}
G. Benettin and A. Giorgilli.
On the Hamiltonian interpolation of near-to-the identity symplectic mappings with application to symplectic integration algorithms.
\emph{J. Stat. Physics} \textbf{74} (1994), 1117--1143.

\bibitem[Bes32]{Besicovitch}
A.S. Besicovitch.
\emph{Almost Periodic Functions.}
Cambridge University Press 1932.

\bibitem[Bo32]{Bohr}
H. Bohr.
\emph{Fastperiodische Funktionen.}
Springer-Verlag, Berlin 1932.

\bibitem[BrFie88]{brfi88}
P.~Brunovsk\'y and B.~Fiedler.
 Connecting orbits in scalar reaction diffusion equations.
 \emph{Dynamics Reported} \textbf{1} (1988), 57--89.

\bibitem[BrFie89]{brfi89}
P.~Brunovsk\'y and B.~Fiedler.
Connecting orbits in scalar reaction diffusion equations {II}: The complete solution.
 \emph{J.~Diff.~Eqns.} \textbf{81} (1989), 106--135.

\bibitem[ChIn74]{chin74}
N. Chafee and E. Infante.
A bifurcation problem for a nonlinear parabolic equation.
\emph{J. Applic.~Analysis} \textbf{4} (1974), 17--37.

\bibitem[COS16]{COS}
C.-H. Cho, H. Okamoto, M. Sh\={o}ji. 
A blow-up problem for a nonlinear heat equation in the complex plane of time. 
\emph{Japan J. Ind. Appl. Math.} \textbf{33} (2016), 145--166.

\bibitem[Con78]{Conley}
C.C.~Conley.
\emph{Isolated Invariant Sets and the Morse Index.} 
AMS, Providence R.I. 1978. 

\bibitem[Cor89]{Corduneanu}
C. Corduneanu.
\emph{Almost Periodic Functions.} With the collaboration of N. Gheorghiu and V. Barbu.
Chelsea Publishing Company, New York 1989. 

\bibitem[DLY24a]{Yorke-a}
R. De Leo and J.A. Yorke.
Streams and graphs of dynamical systems.
\emph{Qual. Theory Dyn. Syst.} \textbf{24} (2024);
\url{https://doi.org/10.1007/s12346-024-01112-x}

\bibitem[Dev22]{Devaney}
R.L. Devaney.
\emph{An Introduction to Chaotic Dynamical Systems.}
CRC Press, Boca Raton FL 2022. 

\bibitem[FKW23]{Fasondini23}
M. Fasondini, J.R. King, and J.A.C. Weideman.
Blow up in a periodic semilinear heat equation.
\emph{Physica D} \textbf{446} (2023); 
\url{https://doi.org/10.1016/j.physd.2023.133660} 

\bibitem[FKW24]{Fasondini24}
M. Fasondini, J.R. King, and J.A.C. Weideman.
Complex-plane singularity dynamics for blow up in a nonlinear heat equation: analysis and computation. 
\emph{Nonlinearity} \textbf{37} (2024);
\url{https://doi.org/10.1088/1361-6544/ad700b}

\bibitem[Fie02]{fi02}
B. Fiedler (ed.).  \emph{Handbook of Dynamical
Systems} \textbf{2}. Elsevier, Amsterdam 2002.

\bibitem[Fie23]{FiedlerClaudia}
B.~Fiedler.
Real chaos and complex time. (2023); 
\url{https://arxiv.org/abs/2310.08136}

\bibitem[Fie25]{FiedlerYamaguti}
B.~Fiedler.
Real-time blow-up and connection graphs of rational vector fields on the Riemann sphere.
In preparation (2025).

\bibitem[FieRo23]{firoSFB}
B.~Fiedler and C.~Rocha.
Design of Sturm global attractors 1: Meanders with three noses, and reversibility.
\emph{Chaos} \textbf{33}, 083127 (2023); \url{https://doi.org/10.1063/5.0147634}

\bibitem[FieRo24]{firoFusco}
B.~Fiedler and C.~Rocha.
Design of Sturm global attractors 2: Time-reversible Chafee-Infante lattices of 3-nose meanders.
\emph{São Paulo J. Math. Sciences.} (2024);
\url{https://doi.org/10.1007/s40863-023-00385-5}

\bibitem[FieSch96]{FiedlerScheurle}
B. Fiedler and J. Scheurle.
\emph{Discretization of Homoclinic Orbits, Rapid Forcing and “Invisible” Chaos.}
Mem. Am. Math. Soc. \textbf{570}, Providence R.I. 1996. 

\bibitem[FieStu25]{fiestu24}
B.~Fiedler and H. Stuke.
Real eternal PDE solutions are not complex entire: a quadratic parabolic example.
\emph{J. Ell. Par. Eqs.} (2025), 53 pp;
\url{doi.org/10.1007/s41808-024-00309-0}

\bibitem[For81]{Forster}
O. Forster.
\emph{Lectures on Riemann Surfaces.}
Springer-Verlag, New York 1981.

\bibitem[GL01]{Gelfreich01}
V.G. Gelfreich and V.F. Lazutkin. 
Splitting of separatrices: perturbation theory and exponential smallness. 
\emph{Russ. Math. Surv.} \textbf{56} (2001), 499--558.

\bibitem[Gel02]{Gelfreich02}
V.G. Gelfreich.
Numerics and exponential smallness. 
In \cite{fi02} (2002), 265--312. 

\bibitem[GH83]{GuckenheimerHolmes}
J. Guckenheimer and P. Holmes.
\emph{Nonlinear Oscillations, Dynamical Systems, and Bifurcations of Vector Fields.}
Springer-Verlag, New York 1983. 

\bibitem[GNSY13]{Yanagida}
J.-S. Guo, H. Ninomiya, M. Shimojo, and E. Yanagida.
Convergence and blow-up of solutions for a complex-valued heat equation with a quadratic nonlinearity. 
\emph{Trans. Am. Math. Soc.} \textbf{365} (2013), 2447--2467. 

\bibitem[Hartm02]{Hartman}
Ph. Hartman.
\emph{Ordinary Differential Equations.}
SIAM, Providence RI 2002.

\bibitem[Harts77]{Hartshorne}
R. Hartshorne.
\emph{Algebraic Geometry.}
Springer-Verlag, New York 1977.

\bibitem[IY08]{Ilyashenko}
Y. Ilyashenko and S. Yakovenko.
\emph{Lectures on Analytic Differential Equations.}
AMS, Providence RI 2008.

\bibitem[Jaq21]{Jaquetteqp}
J. Jaquette.
Quasiperiodicity and blowup in integrable subsystems of nonconservative nonlinear Schrödinger equations. 
\emph{J. Dyn. Differ. Eqs.} \textbf{36} (2024), 1--25. 
\url{https://doi.org/10.1007/s10884-021-10112-3}

\bibitem[JLT22a]{JaquetteHet}
J.~Jaquette, J.-P.~Lessard, A.~Takayasu.
Global dynamics in nonconservative nonlinear Schrödinger equations. 
\emph{Adv. Math.} \textbf{398} (2022), 108234. 

\bibitem[JLT22b]{JaquetteStuke}
J.~Jaquette, J.-P.~Lessard, A.~Takayasu.
Singularities and heteroclinic connections in complex-valued evolutionary equations with a quadratic nonlinearity. 
\emph{Commun. Nonlinear Sci. Numer. Simul.} \textbf{107} (2022) 106188.

\bibitem[Jo06]{Jost}
J. Jost.
\emph{Compact Riemann Surfaces. An Introduction to Contemporary Mathematics.}
Springer-Verlag, Berlin 2006.

\bibitem[KMM04]{Mischaikow}
T.~Kaczynski, K.~Mischaikow, M.~Mrozek.
\emph{Computational Homology.} 
Springer-Verlag, New York 2004. 

\bibitem[Kat92]{Katok}
S. Katok. 
\emph{Fuchsian Groups.} 
Univ. of Chicago Press, 1992.

\bibitem[KSK17]{Kevrekidis}
P.G. Kevrekidis, C.I. Siettos, and  Y.G. Kevrekidis.
To infinity and some glimpses of beyond.
\emph{Nature Comm.} \textbf{8} (2017), 1562; \url{https://doi.org/10.1038/s41467-017-01502-7}

\bibitem[Ku1440]{Cusanus}
Nikolaus von Kues.
\emph{De docta ignorantia.}
Bernkastel-Kues a.d. Mosel, 1440.
 
\bibitem[Lam09]{Lamotke}
K. Lamotke.
\emph{Riemannsche Flächen.}
Springer-Verlag, Heidelberg 2009.

\bibitem[Lap23]{LappicyBlowup}
P. Lappicy and E. Beatriz. 
An energy formula for fully nonlinear degenerate parabolic equations in one spatial dimension.
\emph{Math. Ann.} (2023); \url{https://doi.org/10.1007/s00208-023-02740-5}

\bibitem[LS08]{LiSinai}
D. Li and Y.G. Sinai.
Blow ups of complex solutions of the 3D Navier-Stokes system and renormalization group method.
\emph{J. Eur. Math. Soc.} \textbf{10} (2008), 267--313;
\url{https://doi.org/10.4171/JEMS/111}

\bibitem[LJR83]{Interpol1}
G.G.~Lorentz, K.Jetter, S.D.~Riemenschneider.
\emph{Birkhoff Interpolation.}
Enc. Math. Applic. \textbf{19}, Addison-Wesley,  Reading, Massachusetts1983. 

\bibitem[MH99]{Marsden}
J.E. Marsden and M.J. Hoffman.
\emph{Basic Complex Analysis.}
Freeman, New York 1999.

\bibitem[Mas82]{Masuda1}
K. Masuda.
Blow-up of solutions of some nonlinear diffusion equations.  
\emph{North-Holland Math. Stud.} \textbf{81} (1982), 119--131. 

\bibitem[Mas84]{Masuda2}
K. Masuda.
Analytic solutions of some nonlinear diffusion equations. 
\emph{Math. Z.} \textbf{187} (1984), 61--73. 

\bibitem[Mat01]{MatthiesDiss}
K. Matthies.
Time-averaging under fast periodic forcing of parabolic partial differential equations: Exponential estimates.
\emph{J. Differ. Eqs.} \textbf{174} (2001), 133--180. 

\bibitem[Mat03a]{Matthieshom}
K. Matthies.
Exponentially small splitting of homoclinic orbits of parabolic differential equations under periodic forcing.
\emph{Discr. Contin. Dyn. Syst.} \textbf{9}  (2003), 585--602. 

\bibitem[MS03]{MatthiesScheel}
K. Matthies and A. Scheel.
Exponential averaging for Hamiltonian evolution equations.
\emph{Trans. Am. Math. Soc.} \textbf{355} (2003), 747--773. 

\bibitem[MP00]{Interpol2}
H.N. Mhaskar and D.V. Pai. 
\emph{Fundamentals of Approximation Theory.} 
Narosa, India 2000.

\bibitem[Nei84]{Neishtadt}
A.I. Neishtadt. 
On the separation of motions in systems with rapidly rotating phase. 
\emph{J. Appl. Math. Mech.} \textbf{48} (1984), 134--139. 

\bibitem[Pa69]{Palis}
J.~Palis.
 On Morse-Smale dynamical systems.
 \emph{Topology} \textbf{8} (1969), 385--404.

\bibitem[PaSm70]{PalisSmale}
J.~Palis and S.~Smale.
 Structural stability theorems.
 In \emph{Global Analysis},  S. Chern, S. Smale (eds.). Proc. Symp. in
 Pure Math.~vol.~XIV.~AMS, Providence 1970. 

\bibitem[PdM82]{PalisdeMelo}
J.~Palis and W.~de~Melo.
\emph{Geometric Theory of Dynamical Systems. An Introduction.}
Springer-Verlag, New York 1983. 

\bibitem[Poi52]{Poincare3body}
H. Poincaré.
Sur le problème des trois corps et les équations de la dynamique. 
\emph{{\OE}uvres de Henri Poincaré}, Tome VII, Gauthier-Villars Paris 1952, 262--469.

\bibitem[Rei99]{Reich}
S. Reich. 
Backward error analysis for numerical integrators. 
\emph{SIAM J. Numer. Analysis} \textbf{36} (1999), 1549--1570.

\bibitem[Rel40]{Rellich}
F. Rellich.
Elliptische Funktionen und die ganzen Lösungen von $y''=f(y)$.
\emph{Math. Z.} \textbf{47} (1940), 153--160.
\url{https://doi.org/10.1007/bf01180954}

\bibitem[Shi03]{Interpol3}
Y. G. Shi. 
\emph{Theory of Birkhoff Interpolation.} 
Nova Science Publishers, New York 2003.

\bibitem[SSTC01]{Shilnikov}
L.P. Shilnikov, A.L. Shilnikov, D.V. Turaev, and L.O. Chua.
\emph{Methods of Qualitative Theory in Nonlinear Dynamics. II.}
World Scientific, Singapore 2001. 

\bibitem[Slo24]{oeis}
L.P.A. Sloane (ed.).
The Online encyclopedia of integer sequences.
A002995 (2024). \url{https://oeis.org/A002995}

\bibitem[Sot73]{Sotomayor}
J. Sotomayor.
Generic one-parameter families of vector fields on two-dimensional manifolds.
\emph{Publ. Math. I.H.E.S.} \textbf{43} (1973), 5--46.

\bibitem[Sto00]{Treecount2}
A. Stoimenow.
On the number of chord diagrams.
\emph{Disc. Math.} \textbf{218} (2000), 209--233.

\bibitem[Stu17]{Stukediss}
H. Stuke.
\emph{Blow-up in Complex Time.}
Dissertation Thesis, Freie Universität Berlin 2017. 
\url{http://dx.doi.org/10.17169/refubium-11743}

\bibitem[Stu18]{Stukearxiv}
H. Stuke.
\emph{Complex time blow-up of the nonlinear heat equation.}
(2018); \url{https://arxiv.org/abs/1812.10707} 

\bibitem[TLJO22]{JaquetteMasuda}
A.~Takayasu, J.-P.~Lessard, J.~Jaquette, and H.~Okamoto.
Rigorous numerics for nonlinear heat equations in the complex plane of time.
\emph{Numer. Math.} \textbf{151} (2022), 693--752. 

\bibitem[Ush80]{Ushiki-2}
S. Ushiki. 
On unstable manifolds of analytic diffeomorphisms of the plane.
\emph{RIMS Kyoto Kokyuroku} \textbf{403} (1980), 1--7.

\bibitem[Ush81]{Ushiki-N}
S. Ushiki. 
Unstable manifolds of analytic dynamical systems.
\emph{J. Math. Kyoto Univ.} \textbf{21} (1981), 763--785.

\bibitem[Wal72]{Treecount1} 
D.W. Walkup.
The number of plane trees.
\emph{Mathematika} \textbf{19} (1972), 200--204.

\bibitem[Wit41]{Wittich1}
H. Wittich.
Ganze Lösungen der Differentialgleichung $w''=f(w)$.
\emph{Math. Z.} \textbf{47} (1941), 422--426.
\url{https://doi.org/10.1007/BF01180973}

\bibitem[Wit50]{Wittich2}
H. Wittich.
Ganze transzendente Lösungen algebraischer Differentialgleichungen.
\emph{Math. Ann.} \textbf{122} (1950), 37--46.
\url{https://doi.org/10.1007/BF01342967}

\bibitem[WO16]{WulffOliver}
C. Wulff and M. Oliver.
Exponentially accurate Hamiltonian embeddings of symplectic A-stable Runge-Kutta methods for Hamiltonian semilinear evolution equations. 
\emph{Proc. R. Soc. Edinb. A, Math.} \textbf{146} (2016), 1265--1301.

\bibitem[\.{Z}o\l06]{Zoladek}
H. \.{Z}o\l\c{a}dek.
\emph{The Monodromy Group.}
Birkhäuser, Basel 2006.

}}

\end{thebibliography}
\end{document}